\documentclass[11pt]{article}%
\usepackage{amsmath}
\usepackage{graphicx}%
\usepackage{amsfonts}%
\usepackage{amssymb}

\newcommand{\eqnb}{\begin{equation}}
\newcommand{\eqne}{\end{equation}}

\newtheorem{The}{Theorem}

\newtheorem{Cor}[The]{Corollary}
\newtheorem{Lem}{Lemma}
\newtheorem{Pro}{Proposition}
\newtheorem{Rem}{Remark}

\begin{document}

\title{A Complete Algebraic Solution to the Optimal Dynamic Rationing Policy in the
Stock-Rationing Queue \\with Two Demand Classes}
\author{Quan-Lin Li$^{a}$, Yi-Meng Li$^{b}$, Jing-Yu Ma$^{b}$, Heng-Li Liu$^{b}$\\$^{a}$School of Economics and Management\\Beijing University of Technology, Beijing 100124, China\\$^{b}$School of Economics and Management\\Yanshan University, Qinhuangdao 066004, China}
\maketitle

\begin{abstract}
In this paper, we study a stock-rationing queue with two demand classes by
means of the sensitivity-based optimization, and develop a complete algebraic
solution to the optimal dynamic rationing policy. We show that the optimal
dynamic rationing policy must be of transformational threshold type. Based on
this finding, we can refine three sufficient conditions under each of which
the optimal dynamic rationing policy is of threshold type (i.e., critical
rationing level). To do this, we use the performance difference equation to
characterize the monotonicity and optimality of the long-run average profit of
this system, and thus establish some new structural properties of the optimal
dynamic rationing policy by observing any given reference policy. Finally, we
use numerical experiments to demonstrate our theoretical results of the
optimal dynamic rationing policy. We believe that the methodology and results
developed in this paper can shed light on the study of stock-rationing queues
and open a series of potentially promising research.


\end{abstract}

\section{Introduction}

In this paper, we consider a stock-rationing queueing problem of a warehouse
with one type of products and two classes of demands, which may be viewed as
coming from retailers with two different priority levels. Now, such a
stock-rationing warehouse system becomes more and more important in many large
cities under the current COVID-19 environment. For example, Beijing has seven
super-large warehouses, which always supply various daily necessities (e.g.,
vegetables, meat, eggs, seafood and so on) to more than 40 million people
every day. In the warehouses, each type of daily necessities are supplied by
lots of different companies in China and the other countries, which lead to
that the successive supply stream of each type of products can be well
described as a Poisson process. In addition, the two retailers may be regarded
as a large supermarket group and another community retail store group.
Typically, the large supermarket group\ has a higher supply priority than the
community retail store group. When the COVID-19 at Beijing is at a serious
warning, the stock-rationing management of the warehouses and their optimal
stock-rationing policy plays a key role in strengthening the fine management
of the warehouses such that every family at Beijing can have a very
comprehensive life guarantee.

From the perspective of practical applications, such a stock-rationing
queueing problem with multiple demand classes can always be encountered in
many different real areas, for example, assemble-to-order systems,
make-to-stock queues and multiechelon inventory systems by Ha \cite{Ha:1997a};
multi-echelon supply chains by Raghavan and Roy \cite{Rag:2005} and Huang and
Iravani \cite{Hua:2006}; manufacturing by Zhao et al. \cite{Zha:2005};
military operations by Kaplan \cite{Kap:1969}; airline by Lee and Hersh
\cite{Lee:1993}; maritime by John \cite{Joh:1994}; hotel by Bitran and
Mondschein \cite{Bit:1995}; rental business by Papier and Thonemann
\cite{Pap:2010, Pap:2011} and Jain et al. \cite{Jai:2015}; health care by
Papastavrou et al. \cite{Pap:2014}; and so forth. This shows that the
stock-rationing queues with multiple demand classes are not only necessary and
important in many practical applications, but also have their own theoretical interest.

In the stock-rationing queueing systems, the stock rationing policies always
assign different supply priorities to multiple classes of demands. In the
early literature, the so-called \textit{critical rationing level} was imagined
intuitively, and its existence was further proved, e.g., see Veinott Jr
\cite{Vei:1965} and Topkis \cite{Top:1968}. Once the critical rationing level
is given and the on--hand inventory falls below the critical rationing level,
a low priority demand may be either rejected, back-ordered or discarded such
that the left on--hand inventory will be reserved to supply the future high
priority demands. Thus designing and optimizing the critical rationing levels
becomes a basic management way of inventory rationing across multiple demand
classes. So far, analysis of the critical rationing levels has still been
interesting but difficult and challenging in the study of stock-rationing
queues with multiple demand classes. Therefore, this motivates us in the paper
to be interested in the following questions:

\textbf{(P-a)} Does such a critical rationing level exist?

\textbf{(P-b)} If yes, what are the sufficient (or necessary) conditions for
its existence?

\textbf{(P-c)} If no, which useful characteristics can be further found to
study the optimal rationing policy?

It is interesting and challenging to give a complete answer for the above
three problems . For Problems \textbf{(P-a)} and \textbf{(P-b)}, there have
been some research on the the optimal policy of critical rationing levels in
inventory rationing with multiple demand classes. However, so far quite few
studies have focused on Problem \textbf{(P-c)}, since nobody know from what
point of view and how to approach Problem \textbf{(P-c)}. Fortunately, this
paper proposes and develops an algebraic method to be able to solve Problem
\textbf{(P-c)}. Also, this paper uses the algebraic method to be able to
further deal with Problems \textbf{(P-a)} and \textbf{(P-b)} as a special
example. Therefore, this paper can further sharpen those important results
given in the literature.

Different from the submodular (or supermodular) technique of the Markov
decision processes (MDPs) developed in the inventory rationing literature
(e.g., see Ha \cite{Ha:1997a, Ha:1997b, Ha:2000}), this paper applies the
sensitivity-based optimization to develop an algebraic method in the study of
stock-rationing queues with multiple demand classes. The key of the algebraic
method is to discuss the solution $\mathfrak{P}_{i}^{\left(  \mathbf{d}%
\right)  }$ of the linear equation $G^{\left(  \mathbf{d}\right)  }\left(
i\right)  +b=0$ with respect to a basic economic factor (i.e., price, cost and
reward), where $G^{\left(  \mathbf{d}\right)  }\left(  i\right)  $ is the
perturbation realization factor given in (\ref{Equa-20}). We show that for any
given policy $\mathbf{d}$, the solutions $\mathfrak{P}_{i}^{\left(
\mathbf{d}\right)  }$ for $1\leq i\leq K$ play a key role in solving the above
three problems: \textbf{(P-a)} to \textbf{(P-c)}. See Theorem 10 in Subsection
7.3 for more details.

So far some research has applied the MDPs to discuss inventory rationing (and
stock-rationing queues) across multiple demand classes by means of the
submodular (or supermodular) technique, among which important examples include
Ha \cite{Ha:1997a, Ha:1997b, Ha:2000}, Gayon et al. \cite{Gay:2009}, Benjaafar
and ElHafsi \cite{Ben:2006} and Nadar et al. \cite{Nad:2014}. To do this, it
is a key that the structural properties of the optimal policy need to be
identified by using a set of structured value functions that is preserved
under an optimal operator. Based on this finding, the optimal rationing policy
of the inventory rationing across multiple demand classes can be further
described and expressed by means of the structural properties. In many more
general cases, it is not easy and even very difficult to set up the structural
properties of the optimal policy. For this, some stronger model assumptions
have to be further added to guarantee the existence of structural properties
of the optimal policy. The purpose of improving the applicability of the MDPs
motivates us to propose another algebraic method in this paper to develop a
complete algebraic solution to the optimal policy by means of the
sensitivity-based optimization.

The sensitivity-based optimization may be regarded as a new research branch of
the MDPs, which grows out of infinitesimal perturbation analysis of discrete
event dynamic systems, e.g., see Cao \cite{Cao:1994, Cao:2007}. Note that one
key of the sensitivity-based optimization is to set up and use the so-called
performance difference equation, which is based on the perturbation
realization factor as well as the performance potential related to the Poisson
equation. To the best of our knowledge, this paper is the first to apply the
sensitivity-based optimization to study the stock-rationing queues with
multiple demand classes. On such a research line, there are only a few closely
related works in the recent literature, for example, the energy-efficient data
centers by Ma et al. \cite{Ma:2018} and the group-server queues by Xia et al.
\cite{Xia:2018}. Different from Ma et al. \cite{Ma:2018} and Xia et al.
\cite{Xia:2018}, this paper develops a complete algebraic solution to the
optimal dynamic rationing policy of the stock-rationing queue with multiple
demand classes, and shows that the optimal dynamic rationing policy must be of
transformational threshold type, which can lead to refining three sufficient
conditions under each of which the optimal dynamic rationing policy is of
threshold type (see Theorem 10 in Subsection 7.3). In addition, it is
worthwhile to note that our transformational threshold type results are
sharper than the bang-bang control given in Ma et al. \cite{Ma:2018}, Xia et
al. \cite{Xia:2018} and others.{} Therefore, our algebraic method provides not
only a necessary complement of policy spatial structural integrity but also a
new way of optimality proofs when comparing to the frequently-used submodular
(or supermodular) technique of MDPs. Also, the complete algebraic solution to
the optimal dynamic rationing policy can provide more effective support for
numerical computation of the optimal policy and the optimal profit of this system.

Note that the Poisson equations always play a key role in the study of MDPs.
To the best of our knowledge, this paper provides a new general solution with
two free constants: A potential displacement constant, and another
solution-free constant (see Theorem 1 in Section 5). In addition, it is clear
that such a general solution of the Poisson equations can be further extended
and generalized to a block structure of more general models by means of the
RG-factorizations (see Li \cite{Li:2010}). For the Poisson equations, readers
may refer to, such as Bhulai \cite{Bhu:2002}, Makowski and Shwartz
\cite{Mak:2002}, Asmussen and Bladt \cite{Asm:1994}, Bini et al.
\cite{Bin:2016} and Ma et al. \cite{Ma:2018} for more details.

Based on the above analysis, we summarize the main contributions of this paper
as follows:

\begin{description}
\item[(1)] \textit{A complete algebraic solution:} We study a stock-rationing
queue with two demand classes by means of the sensitivity-based optimization,
and provide a complete algebraic solution to the optimal dynamic rationing
policy. We first show that the optimal dynamic policy must be of
transformational threshold type. Then we refine three sufficient conditions
under each of which the optimal dynamic rationing policy is of threshold type.
See Theorem 10 in Subsection 7.3 for details.

\item[(2)] \textit{A unified computational framework:} To the best of our
knowledge, this paper is the first to apply the sensitivity-based optimization
to analyze the stock-rationing queues with multiple demand classes. Thus it is
necessary and useful to describe the three key steps: (a) Setting up a
policy-based Markov process. (b) Constructing a policy-based Poisson equation,
whose general solution can be used to characterize the monotonicity and
optimality of the long-run average profit of this system. (c) Finding the
optimal dynamic rationing policy in the three different areas of the penalty
cost. In addition, the computational framework can sufficiently support
numerical solution of stock-rationing queues with multiple demand classes
while the submodular (or supermodular) technique of MDPs is very difficult to
deal with more general stock-rationing queues.

\item[(3)] \textit{A difficult area of the penalty cost is solved:} In our
algebraic method, it is a key to set up a linear equation: $G^{\left(
\mathbf{d}\right)  }\left(  i\right)  +b=0$ in the penalty cost $P$ for $1\leq
i\leq K$. We show that the solution $\mathfrak{P}_{i}^{\left(  \mathbf{d}%
\right)  }$ of $G^{\left(  \mathbf{d}\right)  }\left(  i\right)  +b=0$ for
$1\leq i\leq K$ plays a key role in finding the optimal dynamic rationing
policy. By using the solution $\mathfrak{P}_{i}^{\left(  \mathbf{d}\right)  }$
for $1\leq i\leq K$, two key indices $P_{L}\left(  \mathbf{d}\right)  $ and
$P_{H}\left(  \mathbf{d}\right)  $ are defined in (\ref{24}) and (\ref{25})
such that the range of the penalty cost $P$ is divided into three different
areas: (a) $P\geq P_{H}\left(  \mathbf{d}\right)  $; (b) $P_{L}\left(
\mathbf{d}\right)  >0$ and $0<P\leq P_{L}\left(  \mathbf{d}\right)  $; and (c)
$P_{L}\left(  \mathbf{d}\right)  <P<P_{H}\left(  \mathbf{d}\right)  $. In the
first two areas, we show that the optimal dynamic rationing policy is of
threshold type; while for the third area, it is more difficult to analyze the
optimal dynamic rationing policy so that the bang-bang control is always
suggested as a rough result for this case in the literature, e.g., see Ma et
al. \cite{Ma:2018} and Xia et al. \cite{Xia:2018}. Unlike those, this paper
provides a detailed analysis for the difficult area: $P_{L}\left(
\mathbf{d}\right)  <P<P_{H}\left(  \mathbf{d}\right)  $, characterizes the
monotonicity and optimality of the long-run average profit of this system, and
further establish some new structural properties of the optimal dynamic
rationing policy by observing any given reference policy. This leads to a
complete algebraic solution to the optimal dynamic rationing policy.

\item[(4)] \textit{Two different methods can sufficiently support each other:}
Note that our algebraic method sets up a complete algebraic solution to the
optimal dynamic rationing policy, thus it can provide not only a necessary
complement of policy spatial structural integrity but also a new way of
optimality proofs when comparing to the frequently-used submodular (or
supermodular) technique of MDPs. On the other hand, since our algebraic method
and the submodular (or supermodular) technique are all important parts of the
MDPs (the former is to use the poisson equations; while the latter is to apply
the optimality equation), it is clear that the two different methods will
sufficiently support each other in the study of stock-rationing queues (and
rationing inventory) with multiple demand classes.
\end{description}

The remainder of this paper is organized as follows. Section 2 provides a
literature review. Section 3 gives model description for the stock-rationing
queue with two demand classes. Section 4 establishes an optimization problem
to find the optimal dynamic rationing policy, in which we set up a
policy-based birth-death process and define a more general reward function.
Section 5 establishes a policy-based Poisson equation and provides its general
solution with two free constants. Section 6 provides an explicit expression
for the perturbation realization factor ${G}^{\left(  \mathbf{d}\right)
}\left(  i\right)  $, and discusses the solution of the linear equation
${G}^{\left(  \mathbf{d}\right)  }\left(  i\right)  +b=0$ in the penalty cost
$P$ for $1\leq i\leq K$. Section 7 discusses the monotonicity and optimality
of the long-run average profit of this system, and finds the optimal dynamic
rationing policy in three different areas of the penalty cost. Section 8
analyzes the stock-rationing queue under a threshold type (statical) rationing
policy. Section 9 uses numerical experiments to demonstrate our theoretical
results of the optimal dynamic rationing policy. Finally, some concluding
remarks are given in Section 10.

\section{Literature Review}

Our current research is related to three literature streams: The first is the
research on stock-rationing queues, critical rationing levels and their MDP
proofs. The second is on the static rationing policy and the dynamic rationing
policy in inventory rationing across multiple demand classes. The third is on
a simple introduction to the sensitivity-based optimization.

The inventory rationing across multiple demand classes was first analyzed by
Veinott Jr \cite{Vei:1965} in the context of inventory control theory. From
then on, some authors have discussed the inventory rationing problems. So far
such inventory rationing has still been interesting and challenging. Readers
may refer to a book by M\"{o}llering \cite{Mol:2007}; survey papers by Kleijn
and Dekker \cite{Kle:1999} and Li et al. \cite{Li:2019}; and a research
classification by Teunter and Haneveld \cite{Teu:2008}, M\"{o}llering and
Thonemann \cite{Mol:2008}, Van Foreest and Wijngaard \cite{Van:2014} and
Alfieri et al. \cite{Alf:2017}.

\vskip                         0.6cm

\textbf{(a) Stock-rationing queues, critical rationing levels and their MDP proofs}

In a rationing inventory system, a critical rationing level was imagined from
early research and practical experience. Veinott Jr \cite{Vei:1965} first
proposed such a critical rationing level; while Topkis \cite{Top:1968} proved
that the critical rationing level really exists and it is an optimal policy.
Similar results were further developed for two demand classes by Evans
\cite{Eva:1968} and Kaplan \cite{Kap:1969}.

It is a most basic problem how to mathematically prove whether a rationing
inventory system has such a critical rationing level. Ha \cite{Ha:1997a} made
a breakthrough by applying the MDPs to analyze the inventory rationing policy
for a stock-rationing queue with exponential production times, Poisson demand
arrivals, lost sales and multiple demand classes. He proved that the optimal
rationing policy is of critical rationing levels, and showed that not only do
the critical rationing levels exist, but also they are monotone and
stationary. Therefore, the optimal rationing policy was characterized as a
monotone constant sequence of critical rationing levels corresponding to the
multiple demand classes.

Since the seminal work of Ha \cite{Ha:1997a}, it has been interesting to
extend and generalize the way to apply the MDPs to deal with the
stock-rationing queues and the rationing inventory systems. Important examples
include the Erlang production times by Ha \cite{Ha:2000} and Gayon et al.
\cite{Gay:2009}; the backorders with two demand classes by Ha \cite{Ha:1997b}
and with multiple demand classes by de V\'{e}ricourt \cite{De:2001, De:2002};
the parallel production channels by Bulut and Fadilo\u{g}lu \cite{Bul:2011};
the batch ordering by Huang and Iravani \cite{Hua:2008}, the batch production
by Pang et al. \cite{Pan:2014}; the utilization of information by Gayon et al.
\cite{GayB:2009} and ElHafsi et al. \cite{ElH:2010}; an assemble-to-order
production system by Benjaafar and ElHafsi \cite{Ben:2006}, ElHafsi
\cite{ElH:2009}, ElHafsi et al. \cite{ElH:2008, ElH:2015}, Benjaafar et al.
\cite{Ben:2011} and Nadar et al. \cite{Nad:2014}; a two-stage tandem
production system by Xu \cite{Xu:2017}, supply chain by Huang and Iravani
\cite{Hua:2006} and van Wijk et al. \cite{van:2019}; periodic review by Frank
et al. \cite{Fra:2003} and Chen et al. \cite{CheS:2010}; dynamic price by Ding
et al. \cite{Din:2006, Din:2016}, Schulte and Pibernik \cite{Sch:2017}; and so forth.

Different from those works in the literature, this paper applies the
sensitivity-based optimization to study a stock-rationing queue with two
demand classes, and provides a complete algebraic solution to the optimal
dynamic rationing policy. To this end, this paper first shows that the optimal
dynamic policy must be of transformational threshold type. Then it refines
three sufficient conditions under each of which the optimal dynamic rationing
policy is of threshold type. In addition, this paper uses the performance
difference equation to characterize the monotonicity and optimality of the
long-run average profit of this system, and establish some new structural
properties of the optimal dynamic rationing policy by observing any given
reference policy.

\vskip                         0.6cm

\textbf{(b) Inventory rationing across multiple demand classes}

In the inventory rationing literature, there exist two kinds of rationing
policies: The static rationing policy, and the dynamic rationing policy. Note
that the dynamic rationing policy allows a threshold rationing level to be
able to change in time, depending on the number and ages of outstanding
orders. In general, the static rationing policy is possible to miss some
chances to further improve system performance, while the dynamic rationing
policy should reflect better by means of various continuously updated
information, the system performance can be improved dynamically. Deshpande et
al. \cite{Des:2003} indicated that the optimal dynamic rationing policy may
significantly reduce the inventory cost compared with the static rationing policy.

If there exist multiple replenishment opportunities, then the ordering
policies are taken as two different types: Periodic review and continuous
review. Therefore, our literature analysis for inventory rationing focuses on
four different classes through combining the rationing policy (static vs.
dynamic) with the inventory review (continuous vs. periodic) as follows:
Static-continuous, static-periodic, dynamic-continuous, and dynamic-periodic.

\textbf{(b-1) }\textit{The static rationing policy (periodic vs. continuous)}

\textit{The periodic review: }Veinott Jr \cite{Vei:1965} is the first to
introduce an inventory rationing across different demand classes and to
propose a critical rationing level (i.e., the static rationing policy) in a
periodic review inventory system with backorders. Subsequent research further
investigated the periodic review inventory system with multiple demand
classes, for example, the $(s,S)$ policy by Cohen et al. \cite{Coh:1988} and
Tempelmeier \cite{Tem:2006}; the $(S-1,S)$ policy by de V\'{e}ricourt
\cite{De:2002} and Ha \cite{Ha:1997a, Ha:1997b}; the lost sales by Dekkeret et
al. \cite{Dek:2002}; the backorders by M\"{o}llering and Thonemann
\cite{Mol:2008}; and the anticipated critical levels by Wang et al.
\cite{Wan:2013}.

\textit{The continuous review:} Nahmias and Demmy \cite{Nah:1981} is the first
to propose and develop a constant critical level $\left(  Q,r,\mathbf{C}%
\right)  $ policy in a continuous review inventory model with multiple demand
classes, where $Q$ is the fixed batch size, $r$ is the reorder point and
$\mathbf{C}=\left(  C_{1},C_{2},\ldots,C_{n-1}\right)  $ is a set of critical
rationing levels for $n$ demand classes. From that time on, some authors have
discussed the constant critical level $\left(  Q,r,\mathbf{C}\right)  $ policy
in continuous review inventory systems. Readers may refer to recent
publications for details, among which are Melchiors et al. \cite{Mel:2000},
Dekkeret et al. \cite{Dek:1998}, Deshpande et al. \cite{Des:2003}, Isotupa
\cite{Iso:2006}, Arslan et al. \cite{Ars:2007}, M\"{o}llering and Thonemann
\cite{Mol:2008, Mol:2010} and Escalona et al. \cite{Esc:2015, Esc:2017}. In
addition, the $\left(  S-1,S,\mathbf{C}\right)  $ inventory system was
discussed by Dekkeret et al. \cite{Dek:2002}, Kranenburg and van Houtum
\cite{Kra:2007} and so on.

\textbf{(b-2) }\textit{The dynamic rationing policy (continuous vs. periodic)}

\textit{The continuous review: }Topkis \cite{Top:1968} is the first to analyze
the dynamic rationing policy, and to indicate that the optimal rationing
policy is a dynamic policy. Evans \cite{Eva:1968} and Kaplan \cite{Kap:1969}
obtained similar results as that in Topkis \cite{Top:1968} for two demand
classes. Melchiors \cite{Mel:2003} considered a dynamic rationing policy in a
$(s,Q)$ inventory system with a key assumption that there was at most one
outstanding order. Teunter and Haneveld \cite{Teu:2008} developed a continuous
time approach to determine the dynamic rationing policy for two Poisson demand
classes, analyzed the marginal cost to determine the optimal remaining time
for each rationing level, and expressed the optimal threshold policy through a
schematic diagram or a lookup table. Fad\i lo\u{g}lu and Bulut \cite{Fad:2010}
proposed a dynamic rationing policy: Rationing with Exponential Replenishment
Flow (RERF), for continuous review inventory systems with either backorders or
lost sales. Wang et al. \cite{Wan:2013b} developed a dynamic threshold
mechanism to allocate backorders when the multiple outstanding orders for
different demand classes exist for the $(Q,R)$ inventory system.

\textit{The periodic review: }For the dynamic rationing policy in a periodic
review inventory system, readers may refer to, such as two demand classes by
Sobel and Zhang \cite{Sob:2001}, Frank et al. \cite{Fra:2003} and Tan et al.
\cite{Tan:2009}; dynamic critical levels and lost sales by Haynsworth and
Price \cite{Hay:1989}; multiple demand classes by Hung and Hsiao
\cite{Hun:2013}; two backorder classes by Chew et al. \cite{Che:2013}; general
demand processes by Hung et al. \cite{Hun:2012}; mixed backorders and lost
sales by Wang and Tang \cite{Wan:2014}; uncertain demand and production rates
by Turgay et al. \cite{Tur:2015}; and incremental upgrading demands by You
\cite{You:2003}.

\vskip                         0.6cm

\textbf{(c) The sensitivity-based optimization}

In the early 1980s, Ho and Cao \cite{Ho:1991} proposed and developed
infinitesimal perturbation method of discrete event dynamic systems (DEDSs),
which is a new research direction for online simulation optimization of DEDSs
since the 1980s. Readers may refer to the excellent books by, such as Cao
\cite{Cao:1994}, Glasserman \cite{Gla:1991} and Cassandras and Lafortune
\cite{Cas:2008}.

Cao et al. \cite{Cao:1996} and Cao and Chen \cite{Cao:1997} published a
seminal work that transforms the infinitesimal perturbation, together with the
MDPs and the reinforcement learning, into the sensitivity-based optimization.
On this research line, the excellent book by Cao \cite{Cao:2007} summarized
the main results in the study of sensitivity-based optimization. Li and Liu
\cite{Li:2004b} and Chapter 11 in Li \cite{Li:2010} further extended and
generalized the sensitivity-based optimization to a more general framework of
perturbed Markov processes. In addition, the sensitivity-based optimization
can be effectively developed by means of the matrix-analytic method by Neuts
\cite{Neu:1981}, Latouche and Ramaswami \cite{Lat:1999}, and the
RG-factorizations of block-structured Markov processes by Li \cite{Li:2010}
and Ma et al. \cite{Ma:2018}.

So far some research has applied the sensitivity-based optimization to analyze
the MDPs of queues and networks, e.g., see Xia and Cao \cite{Xia:2012}, Xia
and Shihada \cite{Xia:2013}, Xia et al. \cite{Xia:2017, Xia:2018}, Ma et al.
\cite{Ma:2018} and a survey paper by Li et al. \cite{Li:2019}.

\vskip                         0.6cm

Finally, to the best of our knowledge, this paper is the first to apply the
sensitivity-based optimization to analyze the stock-rationing queues with
multiple demand classes. Our algebraic method sets up a complete algebraic
solution to the optimal dynamic rationing policy, thus it provides not only a
necessary complement of policy spatial structural integrity but also some new
proofs of monotonicity and optimality. Note that our algebraic method and the
submodular (or supermodular) technique are all important parts of the MDPs.
The former is to mainly use the poisson equations, which are well related to
Markov (reward) processes; while the latter is to focus on the optimality
equation by applying the monotonous operator theory to prove the optimality.
Therefore, it is clear that the two different methods can sufficiently support
each other in the study of stock-rationing queues with multiple demand
classes. We believe that the methodology and results developed in this paper
can be applicable to the study of stock-rationing queues and open a series of
potentially promising research.

\section{Model Description}

In this section, we describe a stock-rationing queue with two demand classes,
in which a single class of products are supplied to stock at a warehouse, and
the two classes of demands come from two retailers with different priorities.
In addition, we provide system structure, operational mode and mathematical notations.

\textbf{A stock-rationing queue:} The warehouse has the maximal capacity $N$
to stock a single class of products, and the warehouse needs to pay a holding
cost $C_{1}$ per product per unit time. There are two classes of demands to
order the products, in which the demands of Class 1 have a higher priority
than that of Class 2, such that the demands of Class 1 can be satisfied in any
non-zero inventory; while the demands of Class 2 may be either satisfied or
refused based on the inventory level of the products. Figure 1 depicts a
simple physical system to understand the stock-rationing queue.

\begin{figure}[pth]
\centering           \includegraphics[width=12cm]{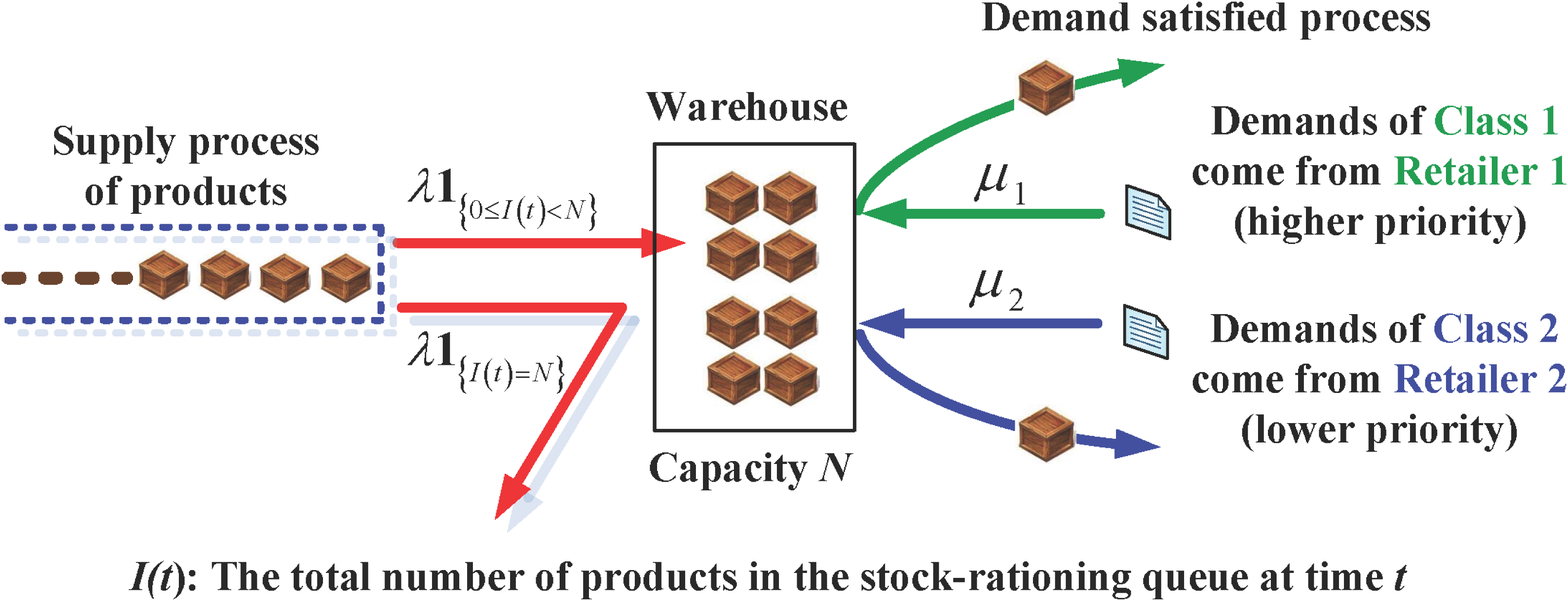}  \caption{A
stock-rationing queue with two demand classes}%
\label{figure:Fig-1}%
\end{figure}

\textbf{The supply process}: The supply stream of the products to the
warehouse is a Poisson process with arrival rate $\lambda$, where the price of
per product is $C_{3}$ paid by the warehouse to the external product supplier.
If the warehouse is full of the products, then any new arriving product has to
be lost. In this case, the warehouse will have an opportunity cost $C_{4}$ per
product rejected into the warehouse.

\textbf{The} \textbf{service processes:} The service times provided by the
warehouse to satisfy the demands of Classes 1 and 2 are i.i.d. and exponential
with service rates $\mu_{1}$ and $\mu_{2}$, respectively. The service
disciplines for the two classes of demands are all First Come First Serve
(FCFS). The warehouse can obtain the service price $R$ when one product is
sold to Retailer 1 or 2. Note that each demand of Class 1 or 2 is satisfied by
one product every time.

\textbf{The stock-rationing rule}: For the two classes of demands, each demand
of Class 1 can always be satisfied in any non-zero inventory; while for
satisfying the demands of Class 2, we need to consider three different cases
as follows:

\textit{Case one}: \textit{The inventory level is zero.} In this case, there
is no product in the warehouse. Thus any new arriving demand has to be
rejected immediately. This leads to the the lost sales cost $C_{2,1}$ (resp.
$C_{2,2}$) per unit time for any lost demand of Class 1 (resp. 2). We assume
that $C_{2,1}>$ $C_{2,2}$, which is used to guarantee the higher priority
service\ for the demands of Class 1 when comparing to the lower priority for
the demands of Class 2.

\textit{Case two}: \textit{The inventory level is low.} In this case, the
number of products in the warehouse is not more than a key threshold $K$,
where the threshold $K$ is subjectively designed by means of some real
experience. Note that the demands of Class 1 have a higher priority to receive
the products than the demands of Class 2. Thus the warehouse will not provide
any product to satisfy the demands of Class 2 under an equal service condition
if the number of products in the warehouse is not more than $K$. Otherwise,
such a service priority is violated (i.e., the demands of Class 2 are
satisfied from a low stock), so that the warehouse must pay a penalty cost $P$
per product supplied to the demands of Class 2 at a low stock. Note that the
penalty cost $P$ measures different priority levels to provide the products
between the two classes of demands.

\textit{Case three}: \textit{The inventory level is high.} In this case, the
number of products in the warehouse is more than the threshold $K$. Thus the
demands of Classes 1 and 2 can be simultaneously satisfied thanks to enough
products in the warehouse.

\textbf{Independence: }We assume that all the random variables defined above
are independent of each other.

In what follows, we use Table 1 to further summarize some above notations.

\begin{table}[th]
\caption{Some costs and prices in the stock-rationing queue}%
\centering                                  \resizebox{\textwidth}  {!}
{
\begin{tabular}
[c]{c|l}\hline
$C_{1}%
$ & The holding cost per unit time per product stored in the
warehouse\\\hline
$C_{2,1}%
$ & The lost sales cost of each lost demand of Class 1\\\hline
$C_{2,2}%
$ & The lost sales cost of each lost demand of Class 2\\\hline
$C_{3}%
$ & The price of per product paid by the warehouse to the external
product supplier\\\hline
$C_{4}%
$ & The opportunity cost per product rejected into the warehouse\\\hline
$P$ & The penalty cost per product supplied to the demands of Class 2 at a low
stock\\\hline
$R$ & The service price of the warehouse paid by each satisfied demand\\\hline
\end{tabular}%

} \end{table}

\section{ Optimization Model Formulation}

In this section, we establish an optimization problem to find the optimal
dynamic rationing policy in the stock-rationing queue. To do this, we set up a
policy-based birth-death process, and define a more general reward function
with respect to both states and policies of the policy-based birth-death process.

To study the stock-rationing queue with two demand classes, we first need to
define both `states' and `policies' to express stochastic dynamics of the
stock-rationing queue.

Let $I(t)$ be the number of products in the warehouse at time $t$. Then it is
regarded as the state of this system at time $t$. Obviously, all the cases of
State $I(t)$ form a state space as follows:%
\[
\mathbf{\Omega}=\{0,1,2,\ldots,N\}.
\]
Also, State $i\in$\ $\mathbf{\Omega}$ is regarded as an inventory level of
this system.

From the states, some policies are defined with a little bit more complexity.
Let $d_{i}$ be a policy related to State $i\in$\ $\mathbf{\Omega}$, and it
expresses whether or not the warehouse prefers to supply some products to the
demands of Class 2 when the inventory level is not more than the threshold $K$
for $0<K\leq N$. Thus, we have%
\begin{equation}
d_{i}=\left\{
\begin{array}
[c]{ll}%
0, & i=0,\\
0,1, & i=1,2,\ldots,K,\\
1, & i=K+1,K+2,\ldots,N,
\end{array}
\right.  \label{1}%
\end{equation}
where $d_{i}=0$ and $1$ represents that the warehouse rejects and satisfies
the demands of Class 2, respectively. Obviously, not only does the policy
$d_{i}$ depend on State $i\in$\ $\mathbf{\Omega}$, but also it is controlled
by the threshold $K$. Of course, for a special case, if $K=N$, then $d_{i}%
\in\left\{  0,1\right\}  $ for $1\leq i\leq N$.

Corresponding to each state in $\mathbf{\Omega}$, we define a time-homogeneous
policy of the stock-rationing queue as%
\[
\mathbf{d}=(d_{0};d_{1},d_{2},\ldots,d_{K};d_{K+1},d_{K+2},\ldots,d_{N}).
\]
It follows from (\ref{1}) that%
\begin{equation}
\mathbf{d}=(0;d_{1},d_{2},\ldots,d_{K};1,1,\ldots,1). \label{2}%
\end{equation}
Thus Policy $\mathbf{d}$ depends on $d_{i}\in\left\{  0,1\right\}  $, which is
related to State $i$\ for $1\leq i\leq K$. Let all the possible policies of
the stock-rationing queue, given in (2), form a policy space as follows:%
\[
\mathcal{D}=\left\{  \mathbf{d}:\mathbf{d}=(0;d_{1},d_{2},\ldots
,d_{K};1,1,\ldots,1),d_{i}\in\left\{  0,1\right\}  ,1\leq i\leq K\right\}  .
\]

\begin{Rem}
In general, the threshold $K$ is subjective and is designed by means of the
real experience of the warehouse manager. If $K=N$, then the policy is
expressed as
\[
\mathbf{d}=(0;d_{1},d_{2},\ldots,d_{N}).
\]
Thus our $K$-based policy $\mathbf{d}=(0;d_{1},d_{2},\ldots,d_{K}%
;1,1,\ldots,1)$ is more general than Policy $\mathbf{d}=(0;d_{1},d_{2}%
,\ldots,d_{N})$.
\end{Rem}

Let $I^{(\mathbf{d})}(t)$ be the state of the stock-rationing queue at time
$t$ under any given policy $\mathbf{d}\in\mathcal{D}$. Then $\left\{
I^{(\mathbf{d})}(t):t\geq0\right\}  $ is a continuous-time policy-based Markov
process on the state space $\mathbf{\Omega}$\ whose state transition relations
are depicted in Figure 2.

\begin{figure}[th]
\centering       \includegraphics[width=\textwidth]{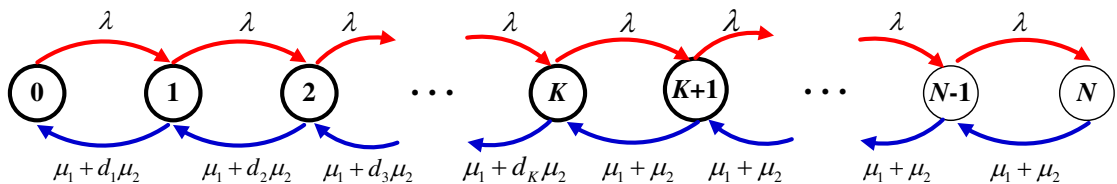}  \caption{State
transition relations of the policy-based Markov process}%
\label{figure:Fig-2}%
\end{figure}

It is easy to see from Figure 2 that $\left\{  I^{(\mathbf{d})}(t):t\geq
0\right\}  $\ is a policy-based birth-death process. Based on this, the
infinitesimal generator of the policy-based birth-death process $\left\{
I^{(\mathbf{d})}(t):t\geq0\right\}  $ is given by%
\begin{equation}
\mathbf{B}^{(\mathbf{d})}=\left(
\begin{array}
[c]{cccccc}%
-\lambda & \lambda &  &  &  & \\
v\left(  d_{1}\right)  & -\left[  \lambda+v\left(  d_{1}\right)  \right]  &
\lambda &  &  & \\
\text{ \ \ \ \ }\ddots & \text{ \ \ \ \ \ \ \ \ \ }\ddots\text{\ } & \text{
\ \ \ \ \ \ \ }\ddots\text{\ } &  &  & \\
& v\left(  d_{K}\right)  & -\left[  \lambda+v\left(  d_{K}\right)  \right]  &
\lambda &  & \\
&  & v\left(  1\right)  & -\left[  \lambda+v\left(  1\right)  \right]  &
\lambda & \\
&  & \text{ \ \ \ \ \ \ \ \ }\ddots & \text{ \ \ \ \ \ \ \ }\ddots & \text{
\ \ \ \ \ \ \ }\ddots & \\
&  &  & v\left(  1\right)  & -\left[  \lambda+v\left(  1\right)  \right]  &
\lambda\\
&  &  &  & v\left(  1\right)  & -v\left(  1\right)
\end{array}
\right)  , \label{3}%
\end{equation}
where $v\left(  d_{i}\right)  =\mu_{1}+d_{i}\mu_{2}$ for $i=1,2,\ldots,K$, and
$v\left(  1\right)  =\mu_{1}+\mu_{2}.$ It is clear that $v\left(
d_{i}\right)  >0$ for $i=1,2,\ldots,K$. Thus the policy-based birth-death
process $\mathbf{B}^{(\mathbf{d})}$\ must be irreducible, aperiodic and
positive recurrent for any given policy $\mathbf{d}\in\mathcal{D}$. In this
case, we write the stationary probability vector of the policy-based
birth-death process $\left\{  I^{(\mathbf{d})}(t):t\geq0\right\}  $\ as%
\begin{equation}
\mathbf{\pi}^{(\mathbf{d})}=\left(  \pi^{(\mathbf{d})}(0);\pi^{(\mathbf{d}%
)}(1),\ldots,\pi^{(\mathbf{d})}(K);\pi^{(\mathbf{d})}(K+1),\ldots
,\pi^{(\mathbf{d})}(N)\right)  . \label{4}%
\end{equation}
Obviously, the stationary probability vector $\mathbf{\pi}^{(\mathbf{d})}$\ is
the unique solution to the system of linear equations: $\mathbf{\pi
}^{(\mathbf{d})}\mathbf{B}^{(\mathbf{d})}=\mathbf{0}$ and $\mathbf{\pi
}^{(\mathbf{d})}\mathbf{e}=1$, where $\mathbf{e}$\ is a column vector of ones
with a suitable dimension. We write \ \
\begin{align}
\xi_{0}  &  =1,\text{ \ \ \ \ \ \ \ \ \ \ \ \ \ \ \ \ \ \ \ \ \ }%
i=0,\nonumber\\
\xi_{i}^{(\mathbf{d})}  &  =\left\{
\begin{array}
[c]{ll}%
\frac{\lambda^{i}}{\prod\limits_{j=1}^{i}v\left(  d_{j}\right)  }, &
i=1,2,\ldots,K,\\
\frac{\lambda^{i}}{\left(  \mu_{1}+\mu_{2}\right)  ^{i-K}\prod\limits_{j=1}%
^{K}v\left(  d_{j}\right)  }, & i=K+1,K+2,\ldots,N,
\end{array}
\right.  \label{5}%
\end{align}
and%
\[
h^{(\mathbf{d})}=1+\sum\limits_{i=1}^{N}\xi_{i}^{(\mathbf{d})}.
\]
It follows from Subsection 1.1.4 of Chapter 1 in Li \cite{Li:2010} that%
\begin{equation}
\pi^{(\mathbf{d})}\left(  i\right)  =\left\{
\begin{array}
[c]{ll}%
\frac{1}{h^{(\mathbf{d})}}, & i=0\\
\frac{1}{h^{(\mathbf{d})}}\xi_{i}^{(\mathbf{d})}, & i=1,2,\ldots,N.
\end{array}
\right.  \label{37}%
\end{equation}

By using the policy-based birth-death process $\mathbf{B}^{(\mathbf{d})}$, now
we define a more general reward function in the stock-rationing queue. It is
seen from Table 1 that the reward function with respect to both states and
policies is defined as a profit rate (i.e. the total system revenue minus the
total system cost per unit time). By observing the impact of Policy
$\mathbf{d}$ on the profit rate, the reward function at State $i$ under Policy
$\mathbf{d}$ is given by%
\begin{align}
f^{(\mathbf{d})}\left(  i\right)   &  =R\left(  \mu_{1}1_{\left\{
i>0\right\}  }+\mu_{2}d_{i}\right)  -C_{1}i-C_{2,1}\mu_{1}1_{\left\{
i=0\right\}  }-C_{2,2}\mu_{2}\left(  1-d_{i}\right) \nonumber\\
\text{ \ }  &  \text{\ \ \ }-C_{3}\lambda1_{\left\{  i<N\right\}  }%
-C_{4}\lambda1_{\left\{  i=N\right\}  }-P\mu_{2}d_{i}1_{\left\{  1\leq i\leq
K\right\}  }, \label{7}%
\end{align}
where, $1_{\left\{  \cdot\right\}  }$ represents the indicator function whose
value is one when the event occurs; otherwise it is zero. By using the
indicator function, satisfying and rejecting the demands of Class 1 are
expressed as $1_{\{i>0\}}$ and $1_{\{i=0\}}$, respectively; the external
products enter or are lost by the warehouse according to $1_{\{i<N\}}$ and
$1_{\{i=N\}}$, respectively; and a penalty cost paid by the warehouse is
denoted by means of $1_{\left\{  1\leq i\leq K\right\}  }$ due to that the
warehouse supplies the products to the demands of Class 2 at a low stock.

For the convenience of readers, it is necessary and useful to explain the
reward function (\ref{7}) from four different cases as follows:

\textbf{Case (a): }For $i=0$,%
\begin{equation}
f\left(  0\right)  =-C_{2,1}\mu_{1}-C_{2,2}\mu_{2}-C_{3}\lambda. \label{8}%
\end{equation}
In Case (a), there is no product in the warehouse, thus it has to reject any
demand of Classes 1 and 2.

\textbf{Case (b): }For $1\leq i\leq K$,%
\begin{equation}
f^{(\mathbf{d})}\left(  i\right)  =R\left(  \mu_{1}+\mu_{2}d_{i}\right)
-C_{1}i-C_{2,2}\mu_{2}\left(  1-d_{i}\right)  -C_{3}\lambda-P\mu_{2}d_{i}.
\label{9}%
\end{equation}
In Case (b), since the inventory level is low for $1\leq i\leq K$, the penalty
cost is paid by the warehouse when it supplies the products to the demands of
Class 2.

Differently from Cases (a) and (b), the inventory level is high for $K+1\leq
i\leq N$ in Cases (c) and (d), thus it can synchronously satisfy the demands
of Classes 1 and 2.

\textbf{Case (c): }For $K+1\leq i\leq N-1$,%
\begin{equation}
f\left(  i\right)  =R\left(  \mu_{1}+\mu_{2}\right)  -C_{1}i-C_{3}\lambda.
\label{10}%
\end{equation}

\textbf{Case (d): }For $i=N$,%
\begin{equation}
f\left(  N\right)  =R\left(  \mu_{1}+\mu_{2}\right)  -C_{1}N-C_{4}\lambda.
\label{11}%
\end{equation}
Note that $C_{3}$ is the price per product paid by the warehouse to the
external product supplier; while $C_{4}$ is the opportunity cost per product
rejected into the warehouse.

Based on the above analysis, we define an $\left(  N+1\right)  $-dimensional
column vector composed of the elements $f\left(  0\right)  ,$ $f^{(\mathbf{d}%
)}\left(  i\right)  $ for $1\leq i\leq K$,\ and $f\left(  j\right)  $\ for
$K+1\leq j\leq N$\ as follows:%
\begin{equation}
\mathbf{f}^{(\mathbf{d})}=\left(  f\left(  0\right)  ;f^{(\mathbf{d})}\left(
1\right)  ,f^{(\mathbf{d})}\left(  2\right)  ,\ldots,f^{(\mathbf{d})}\left(
K\right)  ;f\left(  K+1\right)  ,f\left(  K+2\right)  ,\ldots,f\left(
N\right)  \right)  ^{T}. \label{12}%
\end{equation}

Now, we consider the long-run average profit of the stock-rationing queue (or
the continuous-time policy-based birth-death process $\left\{  I^{(\mathbf{d}%
)}(t):t\geq0\right\}  $) under any given policy $\mathbf{d}$. Let%
\[
\eta^{\mathbf{d}}=\lim_{T\rightarrow\infty}E\left\{  \frac{1}{T}\int_{0}%
^{T}f^{(\mathbf{d})}\left(  I^{(\mathbf{d})}(t)\right)  dt\right\}  .
\]
Then from Section 1 of Chapter 9 of Li \cite{Li:2010}, we obtain%
\begin{equation}
\eta^{\mathbf{d}}=\mathbf{\pi}^{(\mathbf{d})}\mathbf{f}^{(\mathbf{d})},
\label{13}%
\end{equation}
where $\mathbf{\pi}^{(\mathbf{d})}$ and $\mathbf{f}^{(\mathbf{d})}$ are given
by (\ref{4}) and (\ref{12}), respectively.

To further observe the long-run average profit $\eta^{\mathbf{d}}$, here we
show how $\eta^{\mathbf{d}}$ depends on the penalty cost $P$, and
particularly, $\eta^{\mathbf{d}}$ is linear in $P$. To do this, we write that
for $i=0,$
\[
A_{0}=0,\text{ \ }B_{0}=-C_{2,1}\mu_{1}-C_{2,2}\mu_{2}-C_{3}\lambda;
\]
for $i=1,2,\ldots,K,$
\[
A_{i}^{\left(  \mathbf{d}\right)  }=\mu_{2}d_{i},\text{ \ }B_{i}^{\left(
\mathbf{d}\right)  }=R\left(  \mu_{1}+\mu_{2}d_{i}\right)  -C_{1}i-C_{2,2}%
\mu_{2}\left(  1-d_{i}\right)  -C_{3}\lambda;
\]
for $i=K+1,K+2,\ldots,N-1,$%
\[
A_{i}=0,\ B_{i}=R\left(  \mu_{1}+\mu_{2}\right)  -C_{1}i-C_{3}\lambda;
\]
for $i=N,$%
\[
A_{i}=0,\ B_{N}=R\left(  \mu_{1}+\mu_{2}\right)  -C_{1}N-C_{4}\lambda.
\]
Then it follows from (\ref{8}) to (\ref{11}) that for $i=0$,%
\begin{equation}
f\left(  0\right)  =B_{0}; \label{38}%
\end{equation}
for $i=1,2,\ldots,K,$%
\begin{equation}
f^{\left(  \mathbf{d}\right)  }\left(  i\right)  =B_{i}^{\left(
\mathbf{d}\right)  }-PA_{i}^{\left(  \mathbf{d}\right)  }; \label{39}%
\end{equation}
for $i=K+1,K+2,\ldots,N,$%
\begin{equation}
f\left(  i\right)  =B_{i}. \label{40}%
\end{equation}
It follows from (\ref{37}) and (\ref{38}) to (\ref{40}) that
\begin{align}
\eta^{\mathbf{d}}  &  =\mathbf{\pi}^{(\mathbf{d})}\mathbf{f}^{(\mathbf{d}%
)}\nonumber\\
&  =\pi^{\left(  \mathbf{d}\right)  }\left(  0\right)  f\left(  0\right)
+\sum_{i=1}^{K}\pi^{\left(  \mathbf{d}\right)  }\left(  i\right)  f^{\left(
\mathbf{d}\right)  }\left(  i\right)  +\sum_{i=K+1}^{N}\pi^{\left(
\mathbf{d}\right)  }\left(  i\right)  f\left(  i\right) \nonumber\\
&  ={D}^{\left(  \mathbf{d}\right)  }-P{F}^{\left(  \mathbf{d}\right)  },
\label{Equa-9}%
\end{align}
where%
\[
{D}^{\left(  \mathbf{d}\right)  }=\pi^{\left(  \mathbf{d}\right)  }\left(
0\right)  B_{0}+\sum_{i=1}^{K}\pi^{\left(  \mathbf{d}\right)  }\left(
i\right)  B_{i}^{\left(  \mathbf{d}\right)  }+\sum_{i=K+1}^{N}\pi^{\left(
\mathbf{d}\right)  }\left(  i\right)  B_{i}%
\]
and%
\[
{F}^{\left(  \mathbf{d}\right)  }=\sum_{i=1}^{K}\pi^{\left(  \mathbf{d}%
\right)  }\left(  i\right)  A_{i}^{\left(  \mathbf{d}\right)  }.
\]
Hence the long-run average profit $\eta^{\mathbf{d}}$ is linear in the penalty
cost $P$.

We observe that when the inventory level is low, supplying the products to the
demands of Class 2 leads to that both the total system revenue and the total
system cost increase synchronously, and vice versa. Thus there is a tradeoff
between the total system revenue and the total system cost. This motivates us
to find an optimal dynamic rationing policy such that the warehouse has the
maximal profit. Therefore, our objective is to find an optimal dynamic
rationing policy $\mathbf{d}^{\ast}$ such that the long-run average profit
$\eta^{\mathbf{d}}$ is maximal for $\mathbf{d}=\mathbf{d}^{\ast}$, that is,%
\begin{equation}
\mathbf{d}^{\ast}=\mathop{\arg\max}  _{\mathbf{d}\in\mathcal{D}}\left\{
\eta^{\mathbf{d}}\right\}  . \label{14}%
\end{equation}

In fact, it is more difficult and challenging not only to analyze some
interesting structural properties of the optimal rationing policies
$\mathbf{d}^{\ast}$, but also to provide some effective algorithms for
computing the optimal dynamic rationing policy $\mathbf{d}^{\ast}$.

In the remainder of this paper, we apply the sensitivity-based optimization to
study the optimal policy problem (\ref{14}), where the Poisson equations will
play a key role in the study of MDPs (and sensitivity-based optimization).

\section{A Policy-Based Poisson Equation}

In this section, for the stock-rationing queue, we set up a policy-based
Poisson equation which is derived by means of the law of total probability and
analysis on some stopping times of the policy-based birth-death process
$\left\{  I^{\left(  \mathbf{d}\right)  }\left(  t\right)  ,t\geq0\right\}  $.
It is worth noting that the policy-based Poisson equation provides a useful
relation between the sensitivity-based optimization and the MDPs, see, e.g.
Puterman \cite{Put:1994} and Cao \cite{Cao:2007}.

For any given policy $\mathbf{d}\in\mathcal{D}$, it follows from Chapter 2 in
Cao \cite{Cao:2007} that for the continuous-time policy-based birth-death
process $\left\{  I^{\left(  \mathbf{d}\right)  }\left(  t\right)
,t\geq0\right\}  ,$ we define the performance potential as
\begin{equation}
g^{\left(  \mathbf{d}\right)  }\left(  i\right)  =E\left\{  \left.  \int
_{0}^{+\infty}\left[  f^{(\mathbf{d})}\left(  I^{(\mathbf{d})}(t)\right)
-\eta^{\mathbf{d}}\right]  dt\right|  I^{(\mathbf{d})}\left(  0\right)
=i\right\}  , \label{15}%
\end{equation}
where $\eta^{\mathbf{d}}$ is defined in (13). It is seen from Cao
\cite{Cao:2007} that for Policy $\mathbf{d}\in\mathcal{D}$, $g^{\left(
\mathbf{d}\right)  }\left(  i\right)  $ quantifies the contribution of the
initial State $i$ to the long-run average profit of the stock-rationing queue.
Here, $g^{\left(  \mathbf{d}\right)  }\left(  i\right)  $ is also called the
relative value function or the bias in the traditional MDP theory, see, e.g.
Puterman \cite{Put:1994}. We further define a column vector $\mathbf{g}%
^{\left(  \mathbf{d}\right)  }$ as%
\begin{equation}
\mathbf{g}^{\left(  \mathbf{d}\right)  }=\left(  g^{\left(  \mathbf{d}\right)
}\left(  0\right)  ;g^{\left(  \mathbf{d}\right)  }\left(  1\right)
,\ldots,g^{\left(  \mathbf{d}\right)  }\left(  K\right)  ;g^{\left(
\mathbf{d}\right)  }\left(  K+1\right)  ,\ldots,g^{\left(  \mathbf{d}\right)
}\left(  N\right)  \right)  ^{T}. \label{16}%
\end{equation}

To compute the vector $\mathbf{g}^{\left(  \mathbf{d}\right)  }$, we define
the first departure time of the policy-based birth-death process $\left\{
I^{(\mathbf{d})}(t):t\geq0\right\}  $ beginning from State $i$ as%
\[
\tau=\inf\left\{  t\geq0:I^{\left(  \mathbf{d}\right)  }\left(  t\right)  \neq
i\right\}  ,
\]
where $I^{\left(  \mathbf{d}\right)  }\left(  0\right)  =i$. Clearly, $\tau$
is a stopping time of the policy-based birth-death process $\left\{
I^{(\mathbf{d})}(t):t\geq0\right\}  $. Based on this, if $i=0$, then it is
seen from (3) that State $0$ is a boundary state of the policy-based
birth-death process $\mathbf{B}^{(\mathbf{d})}$, hence $I^{\left(
\mathbf{d}\right)  }\left(  \tau\right)  =1$. Similarly, for each State
$i\in\mathbf{\Omega}$, a basic relation is established as follows:%
\begin{equation}
I^{\left(  \mathbf{d}\right)  }\left(  \tau\right)  =\left\{
\begin{array}
[c]{ll}%
1, & i=0,\\
i-1\ \text{or }i+1, & i=1,2,\ldots,N-1,\\
N-1, & i=N.
\end{array}
\right.  \label{17}%
\end{equation}

To compute the column vector $\mathbf{g}^{\left(  \mathbf{d}\right)  }$, we
derive a policy-based Poisson equation in terms of both the stopping time
$\tau$ and the basic relation (17). By using a similar computation to that in
Li and Cao \cite{Li:2004}, we set up the poisson equation according to four
parts as follows:

\textbf{Part (a): }For $i=0$, we have%
\begin{align*}
g^{\left(  \mathbf{d}\right)  }\left(  0\right)   &  =E\left\{  \left.
\int_{0}^{+\infty}\left[  f^{(\mathbf{d})}\left(  I^{(\mathbf{d})}(t)\right)
-\eta^{\mathbf{d}}\right]  dt\right|  I^{(\mathbf{d})}\left(  0\right)
=0\right\} \\
&  =E\left\{  \tau\left|  I^{\left(  \mathbf{d}\right)  }\left(  0\right)
=0\right.  \right\}  \left[  f\left(  0\right)  -\eta^{\mathbf{d}}\right]
+E\left\{  \left.  \int_{\tau}^{+\infty}\left[  f^{(\mathbf{d})}\left(
I^{(\mathbf{d})}(t)\right)  -\eta^{\mathbf{d}}\right]  dt\right|
I^{(\mathbf{d})}\left(  \tau\right)  \right\}  \ \\
&  =\frac{1}{\lambda}\left[  f\left(  0\right)  -\eta^{\mathbf{d}}\right]
+E\left\{  \left.  \int_{0}^{+\infty}\left[  f^{(\mathbf{d})}\left(
I^{(\mathbf{d})}(t)\right)  -\eta^{\mathbf{d}}\right]  dt\right|
I^{(\mathbf{d})}\left(  0\right)  =1\right\} \\
&  =\frac{1}{\lambda}\left[  f\left(  0\right)  -\eta^{\mathbf{d}}\right]
+g^{\left(  \mathbf{d}\right)  }\left(  1\right)  ,
\end{align*}
where for the policy-based\ birth-death process $\left\{  I^{(\mathbf{d}%
)}(t):t\geq0\right\}  $, it is easy to see from Figure 2 that by using
$I^{\left(  \mathbf{d}\right)  }\left(  t\right)  =0$ for $0\leq t<\tau,$%
\begin{align*}
\int_{0}^{\tau}\left[  f^{(\mathbf{d})}\left(  I^{(\mathbf{d})}(t)\right)
-\eta^{\mathbf{d}}\right]  dt  &  =\tau\left[  f\left(  0\right)
-\eta^{\mathbf{d}}\right]  ,\\
E\left\{  \tau\left|  I^{\left(  \mathbf{d}\right)  }\left(  0\right)
=0\right.  \right\}   &  =\frac{1}{\lambda}.
\end{align*}
We obtain%
\begin{equation}
-\lambda g^{\left(  \mathbf{d}\right)  }\left(  0\right)  +\lambda g^{\left(
\mathbf{d}\right)  }\left(  1\right)  =\eta^{\mathbf{d}}-f\left(  0\right)  .
\label{18}%
\end{equation}

\textbf{Part (b): }For $i=1,2,\ldots,K$, it is easy to see from Figure 2 that%
\begin{align*}
g^{\left(  \mathbf{d}\right)  }\left(  i\right)   &  =E\left\{  \left.
\int_{0}^{+\infty}\left[  f^{(\mathbf{d})}\left(  I^{(\mathbf{d})}(t)\right)
-\eta^{\mathbf{d}}\right]  dt\right|  I^{(\mathbf{d})}\left(  0\right)
=i\right\} \\
&  =E\left\{  \tau\left|  I^{\left(  \mathbf{d}\right)  }\left(  0\right)
=i\right.  \right\}  \left[  f^{(\mathbf{d})}\left(  i\right)  -\eta
^{\mathbf{d}}\right]  +E\left\{  \left.  \int_{\tau}^{+\infty}\left[
f^{(\mathbf{d})}\left(  I^{(\mathbf{d})}(t)\right)  -\eta^{\mathbf{d}}\right]
dt\right|  I^{(\mathbf{d})}\left(  \tau\right)  \right\} \\
&  =\frac{1}{v\left(  d_{i}\right)  +\lambda}\left[  f^{(\mathbf{d})}\left(
i\right)  -\eta^{\mathbf{d}}\right] \\
&  \text{ \ \ \ }+\frac{\lambda}{v\left(  d_{i}\right)  +\lambda}E\left\{
\left.  \int_{0}^{+\infty}\left[  f^{(\mathbf{d})}\left(  I^{(\mathbf{d}%
)}(t)\right)  -\eta^{\mathbf{d}}\right]  dt\right|  I^{(\mathbf{d})}\left(
0\right)  =i+1\right\} \\
\text{ \ }  &  \text{\ \ \ \ }+\frac{v\left(  d_{i}\right)  }{v\left(
d_{i}\right)  +\lambda}E\left\{  \left.  \int_{0}^{+\infty}\left[
f^{(\mathbf{d})}\left(  I^{(\mathbf{d})}(t)\right)  -\eta^{\mathbf{d}}\right]
dt\right|  I^{(\mathbf{d})}\left(  0\right)  =i-1\right\} \\
&  =\frac{1}{v\left(  d_{i}\right)  +\lambda}\left[  f^{(\mathbf{d})}\left(
i\right)  -\eta^{\mathbf{d}}\right]  +\frac{\lambda}{v\left(  d_{i}\right)
+\lambda}g^{\left(  \mathbf{d}\right)  }\left(  i+1\right)  +\frac{v\left(
d_{i}\right)  }{v\left(  d_{i}\right)  +\lambda}g^{\left(  \mathbf{d}\right)
}\left(  i-1\right)  ,
\end{align*}
where%
\[
E\left\{  \tau\left|  I^{\left(  \mathbf{d}\right)  }\left(  0\right)
=i\right.  \right\}  =\frac{1}{v\left(  d_{i}\right)  +\lambda}.
\]
We obtain%
\begin{equation}
v\left(  d_{i}\right)  g^{\left(  \mathbf{d}\right)  }\left(  i-1\right)
-\left[  v\left(  d_{i}\right)  +\lambda\right]  g^{\left(  \mathbf{d}\right)
}\left(  i\right)  +\lambda g^{\left(  \mathbf{d}\right)  }\left(  i+1\right)
=\eta^{\mathbf{d}}-f^{(\mathbf{d})}\left(  i\right)  . \label{19}%
\end{equation}

\textbf{Part (c): }For $i=K+1,K+2,\ldots,N-1$, by using Figure 2 we have%
\begin{align*}
g^{\left(  \mathbf{d}\right)  }\left(  i\right)   &  =E\left\{  \left.
\int_{0}^{+\infty}\left[  f^{(\mathbf{d})}\left(  I^{(\mathbf{d})}(t)\right)
-\eta^{\mathbf{d}}\right]  dt\right|  I^{(\mathbf{d})}\left(  0\right)
=i\right\} \\
&  =E\left\{  \tau\left|  I^{\left(  \mathbf{d}\right)  }\left(  0\right)
=i\right.  \right\}  \left[  f\left(  i\right)  -\eta^{\mathbf{d}}\right]
+E\left\{  \left.  \int_{\tau}^{+\infty}\left[  f^{(\mathbf{d})}\left(
I^{(\mathbf{d})}(t)\right)  -\eta^{\mathbf{d}}\right]  dt\right|
I^{(\mathbf{d})}\left(  \tau\right)  \right\} \\
&  =\frac{1}{\mu_{1}+\mu_{2}+\lambda}\left[  f\left(  i\right)  -\eta
^{\mathbf{d}}\right] \\
\text{ \ }  &  \text{\ \ \ \ }+\frac{\lambda}{\mu_{1}+\mu_{2}+\lambda
}E\left\{  \left.  \int_{0}^{+\infty}\left[  f^{(\mathbf{d})}\left(
I^{(\mathbf{d})}(t)\right)  -\eta^{\mathbf{d}}\right]  dt\right|
I^{(\mathbf{d})}\left(  0\right)  =i+1\right\} \\
&  \text{ \ \ \ }+\frac{\mu_{1}+\mu_{2}}{\mu_{1}+\mu_{2}+\lambda}E\left\{
\left.  \int_{0}^{+\infty}\left[  f^{(\mathbf{d})}\left(  I^{(\mathbf{d}%
)}(t)\right)  -\eta^{\mathbf{d}}\right]  dt\right|  I^{(\mathbf{d})}\left(
0\right)  =i-1\right\} \\
&  =\frac{1}{\mu_{1}+\mu_{2}+\lambda}\left[  f\left(  i\right)  -\eta
^{\mathbf{d}}\right]  +\frac{\lambda}{\mu_{1}+\mu_{2}+\lambda}g^{\left(
\mathbf{d}\right)  }\left(  i+1\right)  +\frac{\mu_{1}+\mu_{2}}{\mu_{1}%
+\mu_{2}+\lambda}g^{\left(  \mathbf{d}\right)  }\left(  i-1\right)  ,
\end{align*}
where%
\[
E\left\{  \tau\left|  I^{\left(  \mathbf{d}\right)  }\left(  t\right)
=i\right.  \right\}  =\frac{1}{\mu_{1}+\mu_{2}+\lambda}.
\]
We obtain%
\begin{equation}
\left(  \mu_{1}+\mu_{2}\right)  g^{\left(  \mathbf{d}\right)  }\left(
i-1\right)  -\left(  \mu_{1}+\mu_{2}+\lambda\right)  g^{\left(  \mathbf{d}%
\right)  }\left(  i\right)  +\lambda g^{\left(  \mathbf{d}\right)  }\left(
i+1\right)  =\eta^{\mathbf{d}}-f\left(  i\right)  . \label{33}%
\end{equation}

\textbf{Part (d): }For $i=N$, by using Figure 2 we have%
\begin{align*}
g^{\left(  \mathbf{d}\right)  }\left(  N\right)   &  =E\left\{  \left.
\int_{0}^{+\infty}\left[  f^{(\mathbf{d})}\left(  I^{(\mathbf{d})}(t)\right)
-\eta^{\mathbf{d}}\right]  dt\right|  I^{(\mathbf{d})}\left(  0\right)
=N\right\} \\
&  =E\left\{  \tau\left|  I^{\left(  \mathbf{d}\right)  }\left(  0\right)
=N\right.  \right\}  \left[  f\left(  N\right)  -\eta^{\mathbf{d}}\right]
+E\left\{  \left.  \int_{\tau}^{+\infty}\left[  f^{(\mathbf{d})}\left(
I^{(\mathbf{d})}(t)\right)  -\eta^{\mathbf{d}}\right]  dt\right|
I^{(\mathbf{d})}\left(  \tau\right)  \right\} \\
&  =\frac{1}{\mu_{1}+\mu_{2}}\left[  f\left(  N\right)  -\eta^{\mathbf{d}%
}\right]  +E\left\{  \left.  \int_{0}^{+\infty}\left[  f^{(\mathbf{d})}\left(
I^{(\mathbf{d})}(t)\right)  -\eta^{\mathbf{d}}\right]  dt\right|
I^{(\mathbf{d})}\left(  0\right)  =N-1\right\} \\
&  =\frac{1}{\mu_{1}+\mu_{2}}\left[  f\left(  N\right)  -\eta^{\mathbf{d}%
}\right]  +g^{\left(  \mathbf{d}\right)  }\left(  N-1\right)  ,
\end{align*}
where%
\[
E\left\{  \tau\left|  I^{\left(  \mathbf{d}\right)  }\left(  t\right)
=N\right.  \right\}  =\frac{1}{\mu_{1}+\mu_{2}}.
\]
We obtain%
\begin{equation}
\left(  \mu_{1}+\mu_{2}\right)  g^{\left(  \mathbf{d}\right)  }\left(
N-1\right)  -\left(  \mu_{1}+\mu_{2}\right)  g^{\left(  \mathbf{d}\right)
}\left(  N\right)  =\eta^{\mathbf{d}}-f\left(  N\right)  . \label{20}%
\end{equation}
Thus it follows from (\ref{18}), (\ref{19}), (\ref{33}) and (\ref{20}) that
\begin{equation}
-\mathbf{B}^{\left(  \mathbf{d}\right)  }\mathbf{g}^{\left(  \mathbf{d}%
\right)  }=\mathbf{f}^{(\mathbf{d})}-\eta^{\mathbf{d}}\mathbf{e}, \label{21}%
\end{equation}
where $\mathbf{B}^{\left(  \mathbf{d}\right)  }$, $\mathbf{f}^{(\mathbf{d})}$
and $\eta^{\mathbf{d}}$\ are given in (\ref{3}), (\ref{12}) and (\ref{13}), respectively.

In what follows we provide an effective method to solve the policy-based
Poisson equation, and show that there exist infinitely-many solutions with two
free constants of additive terms. This leads to a general solution with the
two free constants of the policy-based Poisson equation.

To solve the system of linear equations (\ref{21}), it is easy to see that
rank$\left(  \mathbf{B}^{\left(  \mathbf{d}\right)  }\right)  =N$ and
$\det\left(  \mathbf{B}^{\left(  \mathbf{d}\right)  }\right)  =0$ due to the
fact\ that the size of the matrix $\mathbf{B}^{\left(  \mathbf{d}\right)  }$
is $N+1$. Hence, this system of linear equations (\ref{21}) exists
infinitely-many solutions with a free constant of an additive term.

Let $\mathcal{B}$ be a matrix obtained through omitting the first row and the
first column vectors of the matrix $\mathbf{B}^{\left(  \mathbf{d}\right)  }$.
Then,%
\[
\mathcal{B}=\left(
\begin{array}
[c]{cccccc}%
-\left[  \lambda+\nu\left(  d_{1}\right)  \right]  & \lambda &  &  &  & \\
\nu\left(  d_{2}\right)  & -\left[  \lambda+\nu\left(  d_{2}\right)  \right]
& \lambda &  &  & \\
& \ddots\text{ \ \ \ \ \ \ } & \ddots\text{ \ \ \ \ \ \ \ } & \ddots\text{
\ \ \ \ \ \ } &  & \\
& \nu\left(  d_{K}\right)  & -\left[  \lambda+\nu\left(  d_{K}\right)  \right]
& \lambda &  & \\
&  & \nu\left(  1\right)  & -\left[  \lambda+\nu\left(  1\right)  \right]  &
\lambda & \\
&  & \text{\ \ \ \ \ \ \ }\ddots & \ \ \ \ \ \ \  & \ddots\text{
\ \ \ \ \ \ \ \ \ \ }\ddots & \text{ \ \ \ \ \ }\\
&  &  & \nu\left(  1\right)  & -\left[  \lambda+\nu\left(  1\right)  \right]
& \lambda\\
&  &  &  & \nu\left(  1\right)  & -\nu\left(  1\right)
\end{array}
\right)  .
\]
Obviously, rank$\left(  \mathcal{B}\right)  =N.$ Since the size of the matrix
$\mathcal{B}$ is $N$, the matrix $\mathcal{B}$ is invertible, and $\left(
\mathcal{-B}\right)  ^{-1}>0$.

Let $\mathbf{H}^{\left(  \mathbf{d}\right)  }$ and $\mathbf{\varphi
}^{(\mathbf{d})}$ be two column vectors of size $N$ obtained through omitting
the first elements of the two column vectors $\mathbf{f}^{\left(
\mathbf{d}\right)  }-\eta^{\mathbf{d}}\mathbf{e}$ and $\mathbf{g}%
^{(\mathbf{d})}$ of size $N+1$, respectively. Then,%
\[
\mathbf{H}^{\left(  \mathbf{d}\right)  }=\left(
\begin{array}
[c]{c}%
H_{1}^{\left(  \mathbf{d}\right)  }\\
H_{2}^{\left(  \mathbf{d}\right)  }\\
\vdots\\
H_{K}^{\left(  \mathbf{d}\right)  }\\
H_{K+1}^{\left(  \mathbf{d}\right)  }\\
\vdots\\
H_{N}^{\left(  \mathbf{d}\right)  }%
\end{array}
\right)  =\left(
\begin{array}
[c]{c}%
f^{\left(  \mathbf{d}\right)  }\left(  1\right)  -\eta^{\mathbf{d}}\\
f^{\left(  \mathbf{d}\right)  }\left(  2\right)  -\eta^{\mathbf{d}}\\
\vdots\\
f^{\left(  \mathbf{d}\right)  }\left(  K\right)  -\eta^{\mathbf{d}}\\
f\left(  K+1\right)  -\eta^{\mathbf{d}}\\
\vdots\\
f\left(  N\right)  -\eta^{\mathbf{d}}%
\end{array}
\right)  =\left(
\begin{array}
[c]{c}%
\left[  B_{1}^{\left(  \mathbf{d}\right)  }-{D}^{\left(  \mathbf{d}\right)
}\right]  -P\left[  A_{1}^{\left(  \mathbf{d}\right)  }-{F}^{\left(
\mathbf{d}\right)  }\right] \\
\left[  B_{2}^{\left(  \mathbf{d}\right)  }-{D}^{\left(  \mathbf{d}\right)
}\right]  -P\left[  A_{2}^{\left(  \mathbf{d}\right)  }-{F}^{\left(
\mathbf{d}\right)  }\right] \\
\vdots\\
\left[  B_{K}^{\left(  \mathbf{d}\right)  }-{D}^{\left(  \mathbf{d}\right)
}\right]  -P\left[  A_{K}^{\left(  \mathbf{d}\right)  }-{F}^{\left(
\mathbf{d}\right)  }\right] \\
\left[  B_{K+1}-{D}^{\left(  \mathbf{d}\right)  }\right]  -P\left[
A_{K+1}-{F}^{\left(  \mathbf{d}\right)  }\right] \\
\vdots\\
\left[  B_{N}-{D}^{\left(  \mathbf{d}\right)  }\right]  -P\left[  A_{N}%
-{F}^{\left(  \mathbf{d}\right)  }\right]
\end{array}
\right)
\]
and
\[
\mathbf{\varphi}^{(\mathbf{d})}=\left(  g^{\left(  \mathbf{d}\right)  }\left(
1\right)  ,g^{\left(  \mathbf{d}\right)  }\left(  2\right)  ,\ldots,g^{\left(
\mathbf{d}\right)  }\left(  K\right)  ;g^{\left(  \mathbf{d}\right)  }\left(
K+1\right)  ,g^{\left(  \mathbf{d}\right)  }\left(  K+2\right)  ,\ldots
,g^{\left(  \mathbf{d}\right)  }\left(  N\right)  \right)  ^{T}.
\]
Therefore, it follows from (\ref{21}) that%
\begin{equation}
-\mathcal{B}\mathbf{\varphi}^{(\mathbf{d})}=\mathbf{H}^{\left(  \mathbf{d}%
\right)  }+\nu\left(  d_{1}\right)  \mathbf{e}_{1}\mathbf{g}^{\left(
\mathbf{d}\right)  }\left(  0\right)  , \label{Eq-1}%
\end{equation}
where $\mathbf{e}_{1}$ is a column vector with the first element be one and
all the others be zero. Note that the matrix $-\mathcal{B}$ is invertible and
$\left(  -\mathcal{B}\right)  ^{-1}>0$, thus the system of linear equations
(\ref{Eq-1}) always has one unique solution%
\begin{equation}
\mathbf{\varphi}^{(\mathbf{d})}=\left(  -\mathcal{B}\right)  ^{-1}%
\mathbf{H}^{\left(  \mathbf{d}\right)  }+\nu\left(  d_{1}\right)  \left(
-\mathcal{B}\right)  ^{-1}\mathbf{e}_{1}\cdot\Im, \label{Eq-2}%
\end{equation}
where $g^{\left(  \mathbf{d}\right)  }\left(  0\right)  =\Im$ is any given
constant. Let's take a convention%
\[
\left(
\begin{array}
[c]{c}%
a\\
\mathbf{b}%
\end{array}
\right)  =\left(  a,\mathbf{b}\right)  ^{T},
\]
where $\mathbf{b}$ may be a column vector. Then we have%
\begin{align}
\mathbf{g}^{\left(  \mathbf{d}\right)  }  &  =\left(  g^{\left(
\mathbf{d}\right)  }\left(  0\right)  ,\mathbf{\varphi}^{(\mathbf{d})}\right)
^{T}\nonumber\\
&  =\left(  \Im,\left(  -\mathcal{B}\right)  ^{-1}\mathbf{H}^{\left(
\mathbf{d}\right)  }+\nu\left(  d_{1}\right)  \left(  -\mathcal{B}\right)
^{-1}\mathbf{e}_{1}\cdot\Im\right)  ^{T}\nonumber\\
&  =\left(  0,\left(  -\mathcal{B}\right)  ^{-1}\mathbf{H}^{\left(
\mathbf{d}\right)  }\right)  ^{T}+\left(  1,\nu\left(  d_{1}\right)  \left(
-\mathcal{B}\right)  ^{-1}\mathbf{e}_{1}\right)  ^{T}\Im. \label{Eq-2-1}%
\end{align}
Note that $\mathbf{B}^{\left(  \mathbf{d}\right)  }\mathbf{e}=0,$ thus a
general solution to the policy-based Poisson equation is further given by%
\begin{equation}
\mathbf{g}^{\left(  \mathbf{d}\right)  }=\left(  0,\left(  -\mathcal{B}%
\right)  ^{-1}\mathbf{H}^{\left(  \mathbf{d}\right)  }\right)  ^{T}+\left(
1,\nu\left(  d_{1}\right)  \left(  -\mathcal{B}\right)  ^{-1}\mathbf{e}%
_{1}\right)  ^{T}\Im+\xi\mathbf{e}, \label{Eq-2-2}%
\end{equation}
where $\Im$ and $\xi$ are two free constants.

Based on the above analysis, the following theorem summarizes the general
solution of the policy-based Poisson equation.

\begin{The}
\label{The:Poisson}For the Poisson equation $-\mathbf{B}^{\left(
\mathbf{d}\right)  }\mathbf{g}^{\left(  \mathbf{d}\right)  }=\mathbf{f}%
^{(\mathbf{d})}-\eta^{\mathbf{d}}\mathbf{e}$, there exists a key special
solution $\mathbf{g}_{\text{Sp}}^{\mathbf{d}}=\left(  0,\left(  -\mathcal{B}%
\right)  ^{-1}\mathbf{H}^{\left(  \mathbf{d}\right)  }\right)  ^{T}$, and its
general solution is related to two free constants $\Im$ and $\xi$ such that%
\[
\mathbf{g}^{\left(  \mathbf{d}\right)  }=\mathbf{g}_{\text{Sp}}^{\mathbf{d}%
}+\left(  1,\nu\left(  d_{1}\right)  \left(  -\mathcal{B}\right)
^{-1}\mathbf{e}_{1}\right)  ^{T}\Im+\xi\mathbf{e},
\]
where $\xi$ is a potential displacement constant, and $\Im$ is a solution-free constant.
\end{The}

\begin{Rem}
(1) To our best knowledge, this is the first to provide the general solution
of the Poisson equations in the MDPs by means of two different free constants.

(2) Note that $\mathbf{\pi}^{(\mathbf{d})}\mathbf{g}^{\left(  \mathbf{d}%
\right)  }=\eta^{\mathbf{d}}$ and the matrix $-\mathbf{B}^{\left(
\mathbf{d}\right)  }+\mathbf{e\pi}^{(\mathbf{d})}$ is invertible, thus the
Poisson equation $-\mathbf{B}^{\left(  \mathbf{d}\right)  }\mathbf{g}^{\left(
\mathbf{d}\right)  }=\mathbf{f}^{(\mathbf{d})}-\eta^{\mathbf{d}}\mathbf{e}$
can become%
\[
\left(  -\mathbf{B}^{\left(  \mathbf{d}\right)  }+\mathbf{e\pi}^{(\mathbf{d}%
)}\right)  \mathbf{g}^{\left(  \mathbf{d}\right)  }=\mathbf{f}^{(\mathbf{d}%
)}.
\]
This gives a solution of the Poisson equation as follows:%
\[
\mathbf{g}^{\left(  \mathbf{d}\right)  }=\left(  -\mathbf{B}^{\left(
\mathbf{d}\right)  }+\mathbf{e\pi}^{(\mathbf{d})}\right)  ^{-1}\mathbf{f}%
^{(\mathbf{d})}+\xi\mathbf{e},
\]
which is a special solution of the Poisson equation by comparing with that in
Theorem \ref{The:Poisson}.
\end{Rem}

\section{Impact of the Penalty Cost}

In this section, we provide an explicit expression for the perturbation
realization factor of the policy-based birth-death process. Based on this, we
can set up a linear equation in the penalty cost, which is well related to the
performance difference equation. Furthermore, we discuss some useful
properties of policies in the set $\mathcal{D}$ by means of the solution of
the linear equation in the penalty cost.

\subsection{The perturbation realization factor}

We define a perturbation realization factor as%
\begin{equation}
G^{\left(  \mathbf{d}\right)  }\left(  i\right)  \overset{\text{def}}%
{=}g^{\left(  \mathbf{d}\right)  }\left(  i-1\right)  -g^{\left(
\mathbf{d}\right)  }\left(  i\right)  ,i=1,2,\ldots,N. \label{Equa-20}%
\end{equation}
It is easy to see from Cao \cite{Cao:2007} that $G^{\left(  \mathbf{d}\right)
}\left(  i\right)  $ quantifies the difference among two adjacent performance
potentials $g^{\left(  \mathbf{d}\right)  }\left(  i\right)  $ and $g^{\left(
\mathbf{d}\right)  }\left(  i-1\right)  $, and measures the effect on the
long-run average profit of the stock-rationing queue when the system state is
changed from State $i-1$ to State $i$. By using the policy-based Poisson
equation (\ref{21}), we can derive a new system of linear equations, which can
be used to directly express the perturbation realization factor $G^{\left(
\mathbf{d}\right)  }\left(  i\right)  $ for $i=1,2,\ldots,N$.

By using (\ref{Eq-2-2}), we can directly express the perturbation realization
factor $G^{\left(  \mathbf{d}\right)  }\left(  i\right)  $ for $i=1,2,\ldots
,N$. On the other hand, by observing the special structure of the policy-based
Poisson equation (\ref{21}), we can propose a new method of sequence to set up
an explicit expression for $G^{\left(  \mathbf{d}\right)  }\left(  i\right)  $.

For $i=1$, it follows from (\ref{18}) that
\[
-\lambda\left[  g^{\left(  \mathbf{d}\right)  }\left(  0\right)  -g^{\left(
\mathbf{d}\right)  }\left(  1\right)  \right]  =-\lambda G^{\left(
\mathbf{d}\right)  }\left(  1\right)  ,
\]
we have%
\begin{equation}
\lambda G^{\left(  \mathbf{d}\right)  }\left(  1\right)  =f\left(  0\right)
-\eta^{\mathbf{d}}. \label{42}%
\end{equation}
For $i=2,3,\ldots,K$, it follows from (\ref{19}) that%
\begin{align*}
&  \text{ \ \ \ }v\left(  d_{i}\right)  \left[  g^{\left(  \mathbf{d}\right)
}\left(  i-1\right)  -g^{\left(  \mathbf{d}\right)  }\left(  i\right)
\right]  -\lambda\left[  g^{\left(  \mathbf{d}\right)  }\left(  i\right)
-g^{\left(  \mathbf{d}\right)  }\left(  i+1\right)  \right] \\
&  =v\left(  d_{i}\right)  G^{\left(  \mathbf{d}\right)  }\left(  i\right)
-\lambda G^{\left(  \mathbf{d}\right)  }\left(  i+1\right)  ,
\end{align*}
this gives%
\begin{equation}
\lambda G^{\left(  \mathbf{d}\right)  }\left(  i+1\right)  =v\left(
d_{i}\right)  G^{\left(  \mathbf{d}\right)  }\left(  i\right)  +f^{(\mathbf{d}%
)}\left(  i\right)  -\eta^{\mathbf{d}}. \label{43}%
\end{equation}
For $i=K+1,K+2,\ldots,N-1$, it follows from (\ref{33}) that%
\begin{align*}
&  \text{ \ \ \ }\left(  \mu_{1}+\mu_{2}\right)  \left[  g^{\left(
\mathbf{d}\right)  }\left(  i-1\right)  -g^{\left(  \mathbf{d}\right)
}\left(  i\right)  \right]  -\lambda\left[  g^{\left(  \mathbf{d}\right)
}\left(  i\right)  -g^{\left(  \mathbf{d}\right)  }\left(  i+1\right)  \right]
\\
&  =\left(  \mu_{1}+\mu_{2}\right)  G^{\left(  \mathbf{d}\right)  }\left(
i\right)  -\lambda G^{\left(  \mathbf{d}\right)  }\left(  i+1\right)  ,
\end{align*}
we obtain%
\begin{equation}
\lambda G^{\left(  \mathbf{d}\right)  }\left(  i+1\right)  =\left(  \mu
_{1}+\mu_{2}\right)  G^{\left(  \mathbf{d}\right)  }\left(  i\right)
+f\left(  i\right)  -\eta^{\mathbf{d}}. \label{44}%
\end{equation}
For $i=N$, it follows from (\ref{20}) that%
\begin{equation}
\left(  \mu_{1}+\mu_{2}\right)  G^{\left(  \mathbf{d}\right)  }\left(
N\right)  =\eta^{\mathbf{d}}-f\left(  N\right)  . \label{46}%
\end{equation}
By using (\ref{42}), (\ref{43}), (\ref{44}) and (\ref{46}), we obtain a new
system of linear equations satisfied by $G^{\left(  \mathbf{d}\right)
}\left(  i\right)  $ as follows:%
\begin{equation}
\left\{
\begin{array}
[c]{ll}%
\lambda G^{\left(  \mathbf{d}\right)  }\left(  1\right)  =f\left(  0\right)
-\eta^{\mathbf{d}}, & i=1,\\
\lambda G^{\left(  \mathbf{d}\right)  }\left(  i+1\right)  =v\left(
d_{i}\right)  G^{\left(  \mathbf{d}\right)  }\left(  i\right)  +f^{(\mathbf{d}%
)}\left(  i\right)  -\eta^{\mathbf{d}}, & i=2,3,\ldots,K,\\
\lambda G^{\left(  \mathbf{d}\right)  }\left(  i+1\right)  =\left(  \mu
_{1}+\mu_{2}\right)  G^{\left(  \mathbf{d}\right)  }\left(  i\right)
+f\left(  i\right)  -\eta^{\mathbf{d}}, & i=K+1,K+2,\ldots,N-1,\\
\left(  \mu_{1}+\mu_{2}\right)  G^{\left(  \mathbf{d}\right)  }\left(
N\right)  =\eta^{\mathbf{d}}-f\left(  N\right)  , & i=N.
\end{array}
\right.  \label{45}%
\end{equation}

Fortunately, the following theorem can provide an explicit expression for the
perturbation realization factor $G^{\left(  \mathbf{d}\right)  }\left(
i\right)  $ for $1\leq i\leq N$.

\begin{The}
For any given policy $\mathbf{d}$, the perturbation realization factor
$G^{\left(  \mathbf{d}\right)  }\left(  i\right)  $ is\ given by

(a) for $1\leq i\leq K$,%
\begin{equation}
G^{\left(  \mathbf{d}\right)  }\left(  i\right)  =\lambda^{-i}\left[  f\left(
0\right)  -\eta^{\mathbf{d}}\right]  \prod\limits_{k=1}^{i-1}v\left(
d_{k}\right)  +\sum\limits_{r=1}^{i-1}\lambda^{r-i}\left[  f^{\left(
\mathbf{d}\right)  }\left(  r\right)  -\eta^{\mathbf{d}}\right]
\prod\limits_{k=r+1}^{i-1}v\left(  d_{k}\right)  ;\text{ } \label{22}%
\end{equation}

(b) for $K+1\leq i\leq N$,%
\begin{align*}
G^{\left(  \mathbf{d}\right)  }\left(  i\right)   &  =\lambda^{-i}\left[
f\left(  0\right)  -\eta^{\mathbf{d}}\right]  \prod\limits_{k=1}^{K}v\left(
d_{k}\right)  \left[  v\left(  1\right)  \right]  ^{i-K-1}\\
&  \text{ \ \ }+\sum\limits_{r=1}^{K-1}\lambda^{r-K}\left[  f^{\left(
\mathbf{d}\right)  }\left(  r\right)  -\eta^{\mathbf{d}}\right]
\prod\limits_{k=r+1}^{K}v\left(  d_{k}\right)  +\sum\limits_{r=K}^{i-1}%
\lambda^{r-i}\left[  f\left(  r\right)  -\eta^{\mathbf{d}}\right]  \left[
v\left(  1\right)  \right]  ^{i-r-2}.
\end{align*}
\end{The}

\textbf{Proof: }We only prove (a), since the proof of (b) is similar.

It follows from (\ref{45}) that%
\[
G^{\left(  \mathbf{d}\right)  }\left(  1\right)  =\frac{f\left(  0\right)
-\eta^{\mathbf{d}}}{\lambda}.
\]
Similarly, we obtain%
\[
G^{\left(  \mathbf{d}\right)  }\left(  i+1\right)  =\frac{v\left(
d_{i}\right)  }{\lambda}G^{\left(  \mathbf{d}\right)  }\left(  i\right)
+\frac{f^{\left(  \mathbf{d}\right)  }\left(  i\right)  -\eta^{\mathbf{d}}%
}{\lambda},\text{ }i=1,2,\ldots,K.
\]
By using (1.2.4) in Chapter 1 of Elaydi \cite{Ela:1996}, we can obtain the
explicit expression of the perturbation realization factor as follows:%
\[
G^{\left(  \mathbf{d}\right)  }\left(  i\right)  =\lambda^{-i}\left[  f\left(
0\right)  -\eta^{\mathbf{d}}\right]  \prod\limits_{k=1}^{i-1}v\left(
d_{k}\right)  +\sum\limits_{r=1}^{i-1}\lambda^{r-i}\left[  f^{\left(
\mathbf{d}\right)  }\left(  r\right)  -\eta^{\mathbf{d}}\right]
\prod\limits_{k=r+1}^{i-1}v\left(  d_{k}\right)
\]
for $i=1,2,\ldots,K$. This completes the proof. \textbf{{\rule{0.08in}{0.08in}%
}}

\subsection{The performance difference equation}

For any given policy $\mathbf{d}\in\mathcal{D}$, the long-run average profit
of the stock-rationing queue is given by%
\[
\eta^{\mathbf{d}}=\mathbf{\pi}^{\left(  \mathbf{d}\right)  }\mathbf{f}%
^{\left(  \mathbf{d}\right)  },
\]
and the policy-based Poisson equation is given by%
\[
\mathbf{B}^{\left(  \mathbf{d}\right)  }\mathbf{g}^{\left(  \mathbf{d}\right)
}=\eta^{\mathbf{d}}\mathbf{e}-\mathbf{f}^{\left(  \mathbf{d}\right)  }.
\]

It is seen from (\ref{3}) and (\ref{12}) that Policy $\mathbf{d}$ directly
affects not only the elements of the infinitesimal generator $\mathbf{B}%
^{\left(  \mathbf{d}\right)  }$ but also the reward function $\mathbf{f}%
^{\left(  \mathbf{d}\right)  }$. Based on this, if Policy $\mathbf{d}$ changes
to $\mathbf{d}^{\prime}$, then the infinitesimal generator $\mathbf{B}%
^{\left(  \mathbf{d}\right)  }$ and the reward function $\mathbf{f}^{\left(
\mathbf{d}\right)  }$ can have their corresponding changes $\mathbf{B}%
^{\left(  \mathbf{d}^{\prime}\right)  }$ and $\mathbf{f}^{\left(
\mathbf{d}^{\prime}\right)  }$, respectively.

The following lemma provides a useful equation (called performance difference
equation) for the difference $\eta^{\mathbf{d}^{\prime}}-\eta^{\mathbf{d}}$
corresponding to any two different policies $\mathbf{d},\mathbf{d}^{\prime}%
\in\mathcal{D}$. Here, we only restate the performance difference equation
without proof, readers may refer to Cao \cite{Cao:2007} or Ma et al.
\cite{Ma:2018} for more details.

\begin{Lem}
\label{Lem:Diff}For any two policies $\mathbf{d},\mathbf{d}^{\prime}%
\in\mathcal{D}$, we have%
\begin{equation}
\eta^{\mathbf{d}^{\prime}}-\eta^{\mathbf{d}}=\mathbf{\pi}^{\left(
\mathbf{d}^{\prime}\right)  }\left[  \left(  \mathbf{B}^{\left(
\mathbf{d}^{\prime}\right)  }-\mathbf{B}^{\left(  \mathbf{d}\right)  }\right)
\mathbf{g}^{\left(  \mathbf{d}\right)  }+\left(  \mathbf{f}^{\left(
\mathbf{d}^{\prime}\right)  }-\mathbf{f}^{\left(  \mathbf{d}\right)  }\right)
\right]  . \label{28}%
\end{equation}
\end{Lem}

By using the performance difference equation (\ref{28}), we can set up a
partial order relation for the policies in the policy set $\mathcal{D}$ as
follows. For any two policies $\mathbf{d},\mathbf{d}^{\prime}\in\mathcal{D}$,
we write that $\mathbf{d}^{\prime}\succ\mathbf{d}$ if $\eta^{\mathbf{d}%
^{\prime}}>\eta^{\mathbf{d}}$; $\mathbf{d}^{\prime}\thickapprox\mathbf{d}$ if
$\eta^{\mathbf{d}^{\prime}}=\eta^{\mathbf{d}}$; and$\ \mathbf{d}^{\prime}%
\prec\mathbf{d}$ if $\eta^{\mathbf{d}^{\prime}}<\eta^{\mathbf{d}}$. Also, we
write that$\ \mathbf{d}^{\prime}\succeq\mathbf{d}$ if $\eta^{\mathbf{d}%
^{\prime}}\geq\eta^{\mathbf{d}}$; and $\mathbf{d}^{\prime}\preceq\mathbf{d}$
if $\eta^{\mathbf{d}^{\prime}}\leq\eta^{\mathbf{d}}$.

Under this partial order relation, our research target is to find the optimal
policy $\mathbf{d}^{\ast}\in\mathcal{D}$ such that $\mathbf{d}^{\ast}$
$\succeq\mathbf{d}$ for any policy $\mathbf{d}\in\mathcal{D}$, i.e.,%
\[
\mathbf{d}^{\ast}=\underset{\mathbf{d}\in\mathcal{D}}{\arg\max}\left\{
\eta^{\mathbf{d}}\right\}  .
\]
Note that the policy set $\mathcal{D}$ and the state set $\mathbf{\Omega}$ are
all finite, thus an enumeration method using finite comparisons is feasible
for finding the optimal policy $\mathbf{d}^{\ast}$ in the policy set
$\mathcal{D}$.

To find the optimal policy $\mathbf{d}^{\ast}$, we define two policies
$\mathbf{d}$ and $\mathbf{d}^{\prime}$ with an interrelated structure at
Position $i$ as follows:%
\begin{align*}
\mathbf{d}  &  =\left(  0;d_{1},d_{2},\ldots,d_{i-1},\underline{d_{i}}%
,d_{i+1},\ldots,d_{K};1,1,\ldots,1\right)  ,\\
\mathbf{d}^{\prime}  &  =\left(  0;d_{1},d_{2},\ldots,d_{i-1},\underline
{d_{i}^{\prime}},d_{i+1},\ldots,d_{K};1,1,\ldots,1\right)  ,
\end{align*}
where $d_{i}^{\prime},d_{i}\in\left\{  0,1\right\}  $ with $d_{i}^{\prime}\neq
d_{i}$. Clearly, if the two policies $\mathbf{d}$ and $\mathbf{d}^{\prime}$
have an interrelated structure at Position $i$, then only the difference
between the two policies $\mathbf{d}$ and $\mathbf{d}^{\prime}$ is at their
$i$th elements: $d_{i}$ and $d_{i}^{\prime}$.

\begin{Lem}
\label{Lem:PerDE}For the two policies $\mathbf{d}$ and $\mathbf{d}^{\prime}$
with an interrelated structure at Position $i$: $d_{i}$ and $d_{i}^{\prime}$,
we have%
\begin{equation}
\eta^{\mathbf{d}^{\prime}}-\eta^{\mathbf{d}}=\mu_{2}\pi^{\left(
\mathbf{d}^{\prime}\right)  }\left(  i\right)  \left(  d_{i}^{\prime}%
-d_{i}\right)  \left[  {G}^{\left(  \mathbf{d}\right)  }\left(  i\right)
+{b}\right]  , \label{Equa-2}%
\end{equation}
where $b=R+C_{2,2}-P$.
\end{Lem}

\textbf{Proof: }For the two policies $\mathbf{d}$ and $\mathbf{d}^{\prime}$
with an interrelated structure at Position $i$: $d_{i}$ and $d_{i}^{\prime}$,
we have%
\begin{align*}
\mathbf{d}  &  =\left(  0;d_{1},d_{2},\ldots,d_{i-1},\underline{d_{i}}%
,d_{i+1},\ldots,d_{K};1,1,\ldots,1\right)  ,\\
\mathbf{d}^{\prime}  &  =\left(  0;d_{1},d_{2},\ldots,d_{i-1},\underline
{d_{i}^{\prime}},d_{i+1},\ldots,d_{K};1,1,\ldots,1\right)  .
\end{align*}
It is easy to check from (\ref{3}) that
\begin{equation}
\mathbf{B}^{\left(  \mathbf{d}^{\prime}\right)  }-\mathbf{B}^{\left(
\mathbf{d}\right)  }\mathbf{=}\left(
\begin{array}
[c]{ccccccc}%
0 &  &  &  &  &  & \\
0 & \ddots &  &  &  &  & \\
& \ddots & 0 &  &  &  & \\
&  & \left(  d_{i}^{\prime}-d_{i}\right)  \mu_{2} & -\left(  d_{i}^{\prime
}-d_{i}\right)  \mu_{2} &  &  & \\
&  &  & 0 & 0 &  & \\
&  &  &  & \ddots & \ddots & \\
&  &  &  &  & 0 & 0
\end{array}
\right)  . \label{eq-17}%
\end{equation}
Also, from the reward function (\ref{9}), we obtain%
\[
f^{\left(  \mathbf{d}\right)  }\left(  i\right)  =\left(  R+C_{2,2}-P\right)
\mu_{2}d_{i}+R\mu_{1}-C_{1}i-C_{2,2}\mu_{2}-C_{3}\lambda
\]
and%
\[
f^{\left(  \mathbf{d}^{\prime}\right)  }\left(  i\right)  =\left(
R+C_{2,2}-P\right)  \mu_{2}d_{i}^{\prime}+R\mu_{1}-C_{1}i-C_{2,2}\mu_{2}%
-C_{3}\lambda.
\]
This gives
\begin{equation}
\mathbf{f}^{\left(  \mathbf{d}^{\prime}\right)  }-\mathbf{f}^{\left(
\mathbf{d}\right)  }=\left(  0,0,\ldots,0,b\mu_{2}\left(  d_{i}^{\prime}%
-d_{i}\right)  ,0,\ldots,0\right)  ^{T}. \label{eq-18}%
\end{equation}
Thus, it follows from Lemma \ref{Lem:Diff}, (\ref{eq-17})$\ $and (\ref{eq-18})
that%
\begin{align}
\eta^{\mathbf{d}^{\prime}}-\eta^{\mathbf{d}}  &  =\mathbf{\pi}^{\left(
\mathbf{d}^{\prime}\right)  }\left[  \left(  \mathbf{B}^{\left(
\mathbf{d}^{\prime}\right)  }-\mathbf{B}^{\left(  \mathbf{d}\right)  }\right)
\mathbf{g}^{\left(  \mathbf{d}\right)  }+\left(  \mathbf{f}^{\left(
\mathbf{d}^{\prime}\right)  }-\mathbf{f}^{\left(  \mathbf{d}\right)  }\right)
\right] \nonumber\\
&  =\mu_{2}\pi^{\left(  \mathbf{d}^{\prime}\right)  }\left(  i\right)  \left(
d_{i}^{\prime}-d_{i}\right)  \left[  g^{\left(  \mathbf{d}\right)  }\left(
i-1\right)  -g^{\left(  \mathbf{d}\right)  }\left(  i\right)  +b\right]
\nonumber\\
&  =\mu_{2}\pi^{\left(  \mathbf{d}^{\prime}\right)  }\left(  i\right)  \left(
d_{i}^{\prime}-d_{i}\right)  \left[  {G}^{\left(  \mathbf{d}\right)  }\left(
i\right)  +{b}\right]  . \label{eq-18-1}%
\end{align}
This completes the proof. \textbf{{\rule{0.08in}{0.08in}}}

For $d_{i}^{\prime},d_{i}\in\left\{  0,1\right\}  $ with $d_{i}^{\prime}\neq
d_{i}$, we have%
\[
d_{i}^{\prime}-d_{i}=\left\{
\begin{array}
[c]{cc}%
1, & d_{i}^{\prime}=1,d_{i}=0;\\
-1, & d_{i}^{\prime}=0,d_{i}=1.
\end{array}
\right.
\]
Therefore, it is easy to see from (\ref{Equa-2}) that to compare
$\eta^{\mathbf{d}^{\prime}}$ with $\eta^{\mathbf{d}}$, it is necessary to
further analyze the sign of function ${G}^{\left(  \mathbf{d}\right)  }\left(
i\right)  +{b}$. This will be developed in the next subsection.

\subsection{The sign of ${G}^{\left(  \mathbf{d}\right)  }\left(  i\right)
+{b}$}

As seen from (\ref{eq-18-1}), the sign analysis of the performance difference
$\eta^{\mathbf{d}^{\prime}}-\eta^{\mathbf{d}}$ directly depends on that of
${G}^{\left(  \mathbf{d}\right)  }\left(  i\right)  +{b}$. Thus, this
subsection provides the sign analysis of ${G}^{\left(  \mathbf{d}\right)
}\left(  i\right)  +{b}$ with respect to the penalty cost $P$. \ 

Suppose that the inventory level is low. If the service priority is violated
(i.e. the demands of Class 2 are served at a low stock), then the warehouse
has to pay the penalty cost $P$ for each product supplied to the demands of
Class 2. Based on this, we study the influence of the penalty cost $P$ on the
sign of ${G}^{\left(  \mathbf{d}\right)  }\left(  i\right)  +{b}$.

Substituting (\ref{38}), (\ref{39}), (\ref{40}) and (\ref{Equa-9}) into
(\ref{22}), we obtain that for $1\leq i\leq K,$%
\begin{align}
{G}^{\left(  \mathbf{d}\right)  }\left(  i\right)  +{b}  &  =R+C_{2,2}%
+\lambda^{-i}\left[  B_{0}-{D}^{\left(  \mathbf{d}\right)  }\right]
\prod\limits_{k=1}^{i-1}v\left(  d_{k}\right)  +\sum\limits_{r=1}^{i-1}%
\lambda^{r-i}\left[  B_{r}^{\left(  \mathbf{d}\right)  }-{D}^{\left(
\mathbf{d}\right)  }\right]  \prod\limits_{k=r+1}^{i-1}v\left(  d_{k}\right)
\nonumber\\
&  -P\left\{  1+\lambda^{-i}\left[  A_{0}-{F}^{\left(  \mathbf{d}\right)
}\right]  \prod\limits_{k=1}^{i-1}v\left(  d_{k}\right)  +\sum\limits_{r=1}%
^{i-1}\lambda^{r-i}\left[  A_{r}^{\left(  \mathbf{d}\right)  }-{F}^{\left(
\mathbf{d}\right)  }\right]  \prod\limits_{k=r+1}^{i-1}v\left(  d_{k}\right)
\right\}  , \label{40-1}%
\end{align}
which is linear in the penalty cost $P$.

From ${G}^{\left(  \mathbf{d}\right)  }\left(  i\right)  +{b=0}$, we have
\begin{align}
&  \text{ \ \ \ }P\left\{  1+\lambda^{-i}\left[  A_{0}-{F}^{\left(
\mathbf{d}\right)  }\right]  \prod\limits_{k=1}^{i-1}v\left(  d_{k}\right)
+\sum\limits_{r=1}^{i-1}\lambda^{r-i}\left[  A_{r}^{\left(  \mathbf{d}\right)
}-{F}^{\left(  \mathbf{d}\right)  }\right]  \prod\limits_{k=r+1}^{i-1}v\left(
d_{k}\right)  \right\} \nonumber\\
&  =R+C_{2,2}+\lambda^{-i}\left[  B_{0}-{D}^{\left(  \mathbf{d}\right)
}\right]  \prod\limits_{k=1}^{i-1}v\left(  d_{k}\right)  +\sum\limits_{r=1}%
^{i-1}\lambda^{r-i}\left[  B_{r}^{\left(  \mathbf{d}\right)  }-{D}^{\left(
\mathbf{d}\right)  }\right]  \prod\limits_{k=r+1}^{i-1}v\left(  d_{k}\right)
, \label{41}%
\end{align}
thus, the unique solution of the penalty cost $P$ to Equation (\ref{41}) is
given by
\begin{equation}
\mathfrak{P}_{i}^{\left(  \mathbf{d}\right)  }=\frac{R+C_{2,2}+\lambda
^{-i}\left[  B_{0}-{D}^{\left(  \mathbf{d}\right)  }\right]  \prod
\limits_{k=1}^{i-1}v\left(  d_{k}\right)  +\sum\limits_{r=1}^{i-1}%
\lambda^{r-i}\left[  B_{r}^{\left(  \mathbf{d}\right)  }-{D}^{\left(
\mathbf{d}\right)  }\right]  \prod\limits_{k=r+1}^{i-1}v\left(  d_{k}\right)
}{1+\lambda^{-i}\left[  A_{0}-{F}^{\left(  \mathbf{d}\right)  }\right]
\prod\limits_{k=1}^{i-1}v\left(  d_{k}\right)  +\sum\limits_{r=1}^{i-1}%
\lambda^{r-i}\left[  A_{r}^{\left(  \mathbf{d}\right)  }-{F}^{\left(
\mathbf{d}\right)  }\right]  \prod\limits_{k=r+1}^{i-1}v\left(  d_{k}\right)
}. \label{23}%
\end{equation}
It's easy to see from (\ref{40-1}) that if $\mathfrak{P}_{i}^{\left(
\mathbf{d}\right)  }>0$ and $0\leq P\leq\mathfrak{P}_{i}^{\left(
\mathbf{d}\right)  }$, then $G^{\left(  \mathbf{d}\right)  }\left(  i\right)
+b\geq0$; while if $P\geq\mathfrak{P}_{i}^{\left(  \mathbf{d}\right)  },$ then
$G^{\left(  \mathbf{d}\right)  }\left(  i\right)  +b\leq0$. Note that the
equality can hold only if $P=\mathfrak{P}_{i}^{\left(  \mathbf{d}\right)  }$

To understand the solution $\mathfrak{P}_{i}^{\left(  \mathbf{d}\right)  }$
for $1\leq i\leq K$, we use a numerical example to show the solutions in Table
2. To do this, we take the system parameters: $\lambda=3$, $\mu_{1}=4$,
$\mu_{2}=2$, $C_{1}=1$, $C_{2,1}=4$, $C_{2,2}=1$, $C_{3}=5$ and $C_{4}=1$.
Further, we observe three different policies:%
\begin{align*}
\mathbf{d}_{1}  &  =\left(  0;1,1,1,1,1,1,1,1,1,1;1,1,1,1,1\right)  ,\\
\mathbf{d}_{2}  &  =\left(  0;0,0,0,0,0,0,0,0,0,0;1,1,1,1,1\right)  ,\\
\mathbf{d}_{3}  &  =\left(  0;0,0,0,0,0,1,1,1,1,1;1,1,1,1,1\right)  .
\end{align*}

\begin{table}[pbh]
\caption{Numerical analysis of solutions for three different policies}%
\centering                                     \resizebox{\textwidth}  {!}  {
\begin{tabular}
[c]{c|c|c|c|c|c|c|c|c|c|c|c}\hline
$\mathfrak{P}_{i}^{\left(  \mathbf{d}\right)  }$ & $i = \ \ 0$ & 1 & 2 & 3 &
4 & 5 & 6 & 7 & 8 & 9 & 10\\\hline
$\mathbf{d}_{1}$ & 11 & 11 & 38.3 & 20.7 & 21.4 & 21.3 & 20.9 & 20.7 & 20.4 &
20.1 & 19.8\\\hline
$\mathbf{d}_{2}$ & 11 & 11 & -6.5 & 14.3 & 88.4 & 384.9 & $1.5e3$ & $5.6e3 $ &
$2.1e4$ & $7.5e4$ & $2.6e5$\\\hline
$\mathbf{d}_{3}$ & 11 & 11 & -7.1 & 23.0 & -181.4 & -80.4 & -68.6 & 15.0 &
13.7 & 13.3 & 13.1\\\hline
\end{tabular}
}\end{table}

In the stock-rationing queue, we define two critical values related to the
penalty cost $P$ as%
\begin{equation}
P_{H}\left(  \mathbf{d}\right)  =\text{ }\max_{\mathbf{d}\in\mathcal{D}%
}\left\{  0,\mathfrak{P}_{1}^{\left(  \mathbf{d}\right)  },\mathfrak{P}%
_{2}^{\left(  \mathbf{d}\right)  },\ldots,\mathfrak{P}_{K}^{\left(
\mathbf{d}\right)  }\right\}  , \label{24}%
\end{equation}
and%
\begin{equation}
P_{L}\left(  \mathbf{d}\right)  =\min_{\mathbf{d}\in\mathcal{D}}\left\{
\mathfrak{P}_{1}^{\left(  \mathbf{d}\right)  },\mathfrak{P}_{2}^{\left(
\mathbf{d}\right)  },\ldots,\mathfrak{P}_{K}^{\left(  \mathbf{d}\right)
}\right\}  . \label{25}%
\end{equation}
From Table 2, we see that it is possible to have $P_{L}\left(  \mathbf{d}%
\right)  <0$ for Policy $\mathbf{d=d}_{2}$ or $\mathbf{d=d}_{3}$.

The following proposition uses the two critical values $P_{H}\left(
\mathbf{d}\right)  $ and $P_{L}\left(  \mathbf{d}\right)  $, together with the
penalty cost $P$, to provide some sufficient conditions under which the
function $G^{\left(  \mathbf{d}\right)  }\left(  i\right)  +b$ is either
positive, zero or negative.

\begin{Pro}
\label{Pro:Sign}(1) If $P\geq P_{H}\left(  \mathbf{d}\right)  $ for any given
policy $\mathbf{d}\in\mathcal{D}$, then for each $i=1,2,\ldots,K$,%
\[
G^{\left(  \mathbf{d}\right)  }\left(  i\right)  +b\leq0.
\]

(2) If $P_{L}\left(  \mathbf{d}\right)  >0$ and $0\leq P\leq P_{L}\left(
\mathbf{d}\right)  $ for any given policy $\mathbf{d}\in\mathcal{D}$, then for
each $i=1,2,\ldots,K$,%
\[
G^{\left(  \mathbf{d}\right)  }\left(  i\right)  +b\geq0.
\]
\end{Pro}

\textbf{Proof: }(1) For any given policy $\mathbf{d}\in\mathcal{D}$, if $P\geq
P_{H}\left(  \mathbf{d}\right)  $, then it follows from (\ref{24}) that for
each $i=1,2,\ldots,K$,
\[
P\geq\mathfrak{P}_{i}^{\left(  \mathbf{d}\right)  },
\]
this leads to that $G^{\left(  \mathbf{d}\right)  }\left(  i\right)  +b\leq0$.

(2) For any given policy $\mathbf{d}\in\mathcal{D}$, if $P_{L}\left(
\mathbf{d}\right)  >0$ and $0\leq P\leq P_{L}\left(  \mathbf{d}\right)  $,
then\ it follows from (\ref{25}) that for each $i=1,2,\ldots,K$,
\[
0\leq P\leq\mathfrak{P}_{i}^{\left(  \mathbf{d}\right)  },
\]
this gives that $G^{\left(  \mathbf{d}\right)  }\left(  i\right)  +b\geq0$.
This completes the proof. \textbf{{\rule{0.08in}{0.08in}}}\ 

However, for the case with $P_{L}\left(  \mathbf{d}\right)  <P<P_{H}\left(
\mathbf{d}\right)  $ for any given policy $\mathbf{d}\in\mathcal{D}$, it is a
little bit complicated to determine the sign of $G^{\left(  \mathbf{d}\right)
}\left(  i\right)  +b$ for each $i=1,2,\ldots,K$. For this reason, our
discussion will be left in the next section.

For any two policies $\mathbf{d,c}\in\mathcal{D}$,
\begin{align*}
\mathbf{d}  &  =\left(  0;d_{1},d_{2},\ldots,d_{i-1},d_{i},d_{i+1}%
,\ldots,d_{K};1,1,\ldots,1\right)  ,\\
\mathbf{c}  &  =\left(  0;c_{1},c_{2},\ldots,c_{i-1},c_{i},c_{i+1}%
,\ldots,c_{K};1,1,\ldots,1\right)  .
\end{align*}
we write%
\[
S\left(  \mathbf{d,c}\right)  =\left\{  i:d_{i}\neq c_{i},i=1,2,\ldots
,K-1,K\right\}
\]
and its complementary set%
\[
\overline{S\left(  \mathbf{d,c}\right)  }=\left\{  i:d_{i}=c_{i}%
,i=1,2,\ldots,K-1,K\right\}  .
\]
Then%
\[
S\left(  \mathbf{d,c}\right)  \cup\overline{S\left(  \mathbf{d,c}\right)
}=\left\{  1,2,\ldots,K-1,K\right\}  .
\]

The following lemma sets up a policy sequence such that any two adjacent
policies of them have the difference at the corresponding position of only one
element. The proof is easy and is omitted here.

\begin{Lem}
\label{Lem:Pol}For any two policies $\mathbf{d,c}\in\mathcal{D}$, $S\left(
\mathbf{d,c}\right)  =\left\{  i_{1},i_{2},i_{3},\ldots,i_{n-1},i_{n}\right\}
$, then there exist a policy sequence: $\mathbf{d}^{\left(  k\right)  }$ for
$k=1,2,3,\ldots,n-1,n$, such that%
\begin{align*}
S\left(  \mathbf{d,d}^{\left(  1\right)  }\right)   &  =\left\{
j_{1}\right\}  ,\\
S\left(  \mathbf{d}^{\left(  1\right)  }\mathbf{,d}^{\left(  2\right)
}\right)   &  =\left\{  j_{2}\right\}  ,\\
&  \vdots\\
S\left(  \mathbf{d}^{\left(  n-1\right)  }\mathbf{,d}^{\left(  n\right)
}\right)   &  =\left\{  j_{n}\right\}  ,
\end{align*}
where $\mathbf{d}^{\left(  n\right)  }=\mathbf{c}$, and $\left\{  i_{1}%
,i_{2},i_{3},\ldots,i_{n-1},i_{n}\right\}  =\left\{  j_{1},j_{2},j_{3}%
,\ldots,j_{n-1},j_{n}\right\}  $. Also, for $k=1,2,3,\ldots,n-1,n$, we have%
\[
S\left(  \mathbf{d,d}^{\left(  k\right)  }\right)  =\left\{  j_{1},j_{2}%
,j_{3},\ldots,j_{k}\right\}  .
\]
\end{Lem}

The following theorem provides a class property of the policies in the set
$\mathcal{D}$ by means of the function $G^{\left(  \mathbf{c}\right)  }\left(
i\right)  +b$ for any policy $\mathbf{c}\in\mathcal{D}$ and for each $i\in
S\left(  \mathbf{d,c}\right)  $, where Policy $\mathbf{d}$ is any given
reference policy in the set $\mathcal{D}$. Note that the class property will
play a key role in developing some new structural properties of the optimal
dynamic rationing policy.

\begin{The}
\label{The:NP}(1) If $P\geq P_{H}\left(  \mathbf{d}\right)  $ for any given
policy $\mathbf{d}$, then for any policy $\mathbf{c}\in\mathcal{D}$ and for
each $i\in S\left(  \mathbf{d,c}\right)  $,%
\[
G^{\left(  \mathbf{c}\right)  }\left(  i\right)  +b\leq0.
\]

(2) If $P_{L}\left(  \mathbf{d}\right)  >0$ and $0\leq P\leq P_{L}\left(
\mathbf{d}\right)  $ for any given policy $\mathbf{d}$, then for any policy
$\mathbf{c}\in\mathcal{D}$ and for each $i\in S\left(  \mathbf{d,c}\right)  $,%
\[
G^{\left(  \mathbf{c}\right)  }\left(  i\right)  +b\geq0.
\]
\end{The}

\textbf{Proof: }We only prove (1), while (2) can be proved similarly.

If $P\geq P_{H}\left(  \mathbf{d}\right)  $ for any given policy $\mathbf{d}$,
then it follows from (1) of Proposition \ref{Pro:Sign} that for $i=1,2,\ldots
,K,$%
\[
G^{\left(  \mathbf{d}\right)  }\left(  i\right)  +b\leq0.
\]

From Policy $\mathbf{d}$, we observe any different policy $\mathbf{c}%
\in\mathcal{D}$. If the two policies $\mathbf{d}$ and $\mathbf{c}$ have $n$
different elements: $d_{i_{l}}\neq$ $c_{i_{l}}$ for $l=1,2,\ldots,n$, then
$S\left(  \mathbf{d,c}\right)  =\left\{  i_{l}:l=1,2,\ldots,n\right\}  $.

Note that the performance difference equation (\ref{Equa-2}) can only be
applied to two policies $\mathbf{d}^{\prime}$ and $\mathbf{d}$ with an
interrelated structure at Position $i:$ $d_{i}^{\prime},d_{i}\in\left\{
0,1\right\}  $ with $d_{i}^{\prime}\neq d_{i}$, thus for a policy
$\mathbf{c}\in\mathcal{D}$ with $S\left(  \mathbf{d,c}\right)  =\left\{
i_{1},i_{2},i_{3},\ldots,i_{n-1},i_{n}\right\}  $, our following discussion
needs to use the policy sequence: $\mathbf{d}^{\left(  k\right)  }$ for
$k=1,2,3,\ldots,n-1,n$, given in Lemma \ref{Lem:Pol}. To this end, our further
proof is to use the mathematical induction in the following three steps:

\textit{Step one:} Analyzing the two policies $\mathbf{d}$ and $\mathbf{d}%
^{\left(  1\right)  }$.

For each $j_{1}\in\left\{  i_{1},i_{2},i_{3},\ldots,i_{n-1},i_{n}\right\}  $,
we take $S\left(  \mathbf{d,d}^{\left(  1\right)  }\right)  =\left\{
j_{1}\right\}  $. It follows from the performance difference equation
(\ref{Equa-2}) that%
\begin{equation}
\eta^{\mathbf{d}^{\left(  1\right)  }}-\eta^{\mathbf{d}}=\mu_{2}\pi^{\left(
\mathbf{d}^{\left(  1\right)  }\right)  }\left(  j_{1}\right)  \left(
d_{j_{1}}^{\left(  \mathbf{d}^{\left(  1\right)  }\right)  }-d_{j_{1}}\right)
\left[  {G}^{\left(  \mathbf{d}\right)  }\left(  j_{1}\right)  +{b}\right]  .
\label{Equa-10}%
\end{equation}
Similarly, we have%
\begin{equation}
\eta^{\mathbf{d}}-\eta^{\mathbf{d}^{\left(  1\right)  }}=\mu_{2}\pi^{\left(
\mathbf{d}\right)  }\left(  j_{1}\right)  \left(  d_{j_{1}}-d_{j_{1}}^{\left(
\mathbf{d}^{\left(  1\right)  }\right)  }\right)  \left[  {G}^{\left(
\mathbf{d}^{\left(  1\right)  }\right)  }\left(  j_{1}\right)  +{b}\right]  .
\label{Equa-11}%
\end{equation}
It is easy to see from (\ref{Equa-10}) and (\ref{Equa-11}) that%
\begin{equation}
{G}^{\left(  \mathbf{d}^{\left(  1\right)  }\right)  }\left(  j_{1}\right)
+{b=}\frac{\pi^{\left(  \mathbf{d}^{\left(  1\right)  }\right)  }\left(
j_{1}\right)  }{\pi^{\left(  \mathbf{d}\right)  }\left(  j_{1}\right)
}\left[  {G}^{\left(  \mathbf{d}\right)  }\left(  j_{1}\right)  +{b}\right]
\leq0. \label{Equa-11-1}%
\end{equation}
Therefore, for Policy $\mathbf{d}^{\left(  1\right)  }\in\mathcal{D}$,
${G}^{\left(  \mathbf{d}^{\left(  1\right)  }\right)  }\left(  j_{1}\right)
+{b}\leq0$ for each $j_{1}\in\left\{  i_{1},i_{2},i_{3},\ldots,i_{n-1}%
,i_{n}\right\}  $.

\textit{Step two:} Analyzing the two policies $\mathbf{d}^{\left(  1\right)
}$ and $\mathbf{d}^{\left(  2\right)  }$.

For each $j_{2}\in\left\{  i_{1},i_{2},i_{3},\ldots,i_{n-1},i_{n}\right\}  $,
we take $S\left(  \mathbf{d}^{\left(  1\right)  }\mathbf{,d}^{\left(
2\right)  }\right)  =\left\{  j_{2}\right\}  $. It is easy to see from
(\ref{Equa-11-1}) that%
\[
{G}^{\left(  \mathbf{d}^{\left(  2\right)  }\right)  }\left(  j_{2}\right)
+{b=}\frac{\pi^{\left(  \mathbf{d}^{\left(  2\right)  }\right)  }\left(
j_{2}\right)  }{\pi^{\left(  \mathbf{d}^{\left(  1\right)  }\right)  }\left(
j_{2}\right)  }\left[  {G}^{\left(  \mathbf{d}^{\left(  1\right)  }\right)
}\left(  j_{2}\right)  +{b}\right]  \leq0.
\]
Therefore, for Policy $\mathbf{d}^{\left(  2\right)  }\in\mathcal{D}$,
${G}^{\left(  \mathbf{d}^{\left(  2\right)  }\right)  }\left(  j_{2}\right)
+{b}\leq0$ for each $j_{2}\in\left\{  i_{1},i_{2},i_{3},\ldots,i_{n-1}%
,i_{n}\right\}  $.

\textit{Step three:} Assume that for $l=3,4,\ldots,k-2,k-1$, we have obtained
that for Policy $\mathbf{d}^{\left(  l\right)  }\in\mathcal{D}$ with $S\left(
\mathbf{d}^{\left(  l-1\right)  }\mathbf{,d}^{\left(  l\right)  }\right)
=\left\{  j_{l}\right\}  $, ${G}^{\left(  \mathbf{d}^{\left(  l\right)
}\right)  }\left(  j_{l}\right)  +{b}\leq0$ for each $j_{l}\in\left\{
i_{1},i_{2},i_{3},\ldots,i_{n-1},i_{n}\right\}  $. Now, we prove the next case
with $l=k$.

For each $j_{k}\in\left\{  i_{1},i_{2},i_{3},\ldots,i_{n-1},i_{n}\right\}  $,
we take $S\left(  \mathbf{d}^{\left(  k-1\right)  }\mathbf{,d}^{\left(
k\right)  }\right)  =\left\{  j_{k}\right\}  $. It is easy to see from
(\ref{Equa-11-1}) that%
\[
{G}^{\left(  \mathbf{d}^{\left(  k\right)  }\right)  }\left(  j_{k}\right)
+{b=}\frac{\pi^{\left(  \mathbf{d}^{\left(  k\right)  }\right)  }\left(
j_{k}\right)  }{\pi^{\left(  \mathbf{d}^{\left(  k-1\right)  }\right)
}\left(  j_{k}\right)  }\left[  {G}^{\left(  \mathbf{d}^{\left(  k-1\right)
}\right)  }\left(  j_{k}\right)  +{b}\right]  \leq0.
\]
This gives that for Policy $\mathbf{d}^{\left(  k\right)  }\in\mathcal{D}$,
${G}^{\left(  \mathbf{d}^{\left(  k\right)  }\right)  }\left(  j_{k}\right)
+{b}\leq0$ for each $j_{k}\in\left\{  i_{1},i_{2},i_{3},\ldots,i_{n-1}%
,i_{n}\right\}  $. Thus, this result holds for the case with $l=k$.

Following the above analysis, we can prove by induction that for Policy
$\mathbf{d}^{\left(  n\right)  }\in\mathcal{D}$, ${G}^{\left(  \mathbf{d}%
^{\left(  n\right)  }\right)  }\left(  j_{n}\right)  +{b}\leq0$ for each
$j_{n}\in\left\{  i_{1},i_{2},i_{3},\ldots,i_{n-1},i_{n}\right\}  $. Since
$\mathbf{c=d}^{\left(  n\right)  }$, we obtain that for Policy $\mathbf{c}%
\in\mathcal{D}$, ${G}^{\left(  \mathbf{c}\right)  }\left(  i\right)  +{b}%
\leq0$ for each $i\in\left\{  i_{1},i_{2},i_{3},\ldots,i_{n-1},i_{n}\right\}
$. This completes the proof. \textbf{{\rule{0.08in}{0.08in}}}

\section{Monotonicity and Optimality}

In this section, we analyze the optimal dynamic rationing policy in the three
different areas of the penalty cost: $P\geq P_{H}\left(  \mathbf{d}\right)  $;
$P_{L}\left(  \mathbf{d}\right)  >0$ and $0<P\leq P_{L}\left(  \mathbf{d}%
\right)  $; and $P_{L}\left(  \mathbf{d}\right)  <P<P_{H}\left(
\mathbf{d}\right)  $, which are studied as three different subsections,
respectively. For the three areas, some new structural properties of the
optimal dynamic rationing policy are given by using our algebraic method.
Also, it is easy to see that for the first two areas: $P\geq P_{H}\left(
\mathbf{d}\right)  $; and $P_{L}\left(  \mathbf{d}\right)  >0$ and $0<P\leq
P_{L}\left(  \mathbf{d}\right)  $, the optimal dynamic rationing policy is of
threshold type; while for the third area: $P_{L}\left(  \mathbf{d}\right)
<P<P_{H}\left(  \mathbf{d}\right)  $, it may not be of threshold type but must
be of transformational threshold type.

As seen from Lemma \ref{Lem:PerDE}, to compare $\eta^{\mathbf{d}^{\prime}}$
with $\eta^{\mathbf{d}}$, our aim is to focus on only Position $i$ with
$d_{i}^{\prime}\neq d_{i}$ for $d_{i}^{\prime},d_{i}\in\left\{  0,1\right\}
$. Also, Lemma \ref{Lem:Pol} provides a useful class property of the policies
in the set $\mathcal{D}$ under the function $G^{\left(  \mathbf{c}\right)
}\left(  i\right)  +b$ for any policy $\mathbf{c,d}\in\mathcal{D}$ and for
each $i\in S\left(  \mathbf{d,c}\right)  $. The two lemmas are very useful for
our research in the next subsections.

\subsection{The penalty cost $P\geq P_{H}\left(  \mathbf{d}\right)  $}

In this subsection, for the area of the penalty cost: $P\geq P_{H}\left(
\mathbf{d}\right)  $ for any given policy $\mathbf{d}$, we find the optimal
dynamic rationing policy of the stock-rationing queue, and further compute the
maximal long-run average profit of this system.

The following theorem uses the class property of the policies in the set
$\mathcal{D}$, given in (1) of Theorem \ref{The:NP}, to set up some basic
relations between any two policies. Thus, we find the optimal dynamic
rationing policy of the stock-rationing queue.

\begin{The}
\label{The:left}If $P\geq{P}_{H}\left(  \mathbf{d}\right)  $ for any given
policy $\mathbf{d}$, then the optimal dynamic rationing policy of the
stock-rationing queue is given by
\[
\mathbf{d}^{\ast}=\left(  0;0,0,\ldots,0;1,1,\ldots,1\right)  .
\]
This shows that if the penalty cost is higher with $P\geq{P}_{H}\left(
\mathbf{d}\right)  $ for any given policy $\mathbf{d}$, then the warehouse can
not supply any product to the demands of Class $2$.
\end{The}

\textbf{Proof: }If $P\geq{P}_{H}$ for any given policy $\mathbf{d}$, then our
proof will focus on that for any policy $\mathbf{c}\in\mathcal{D}$, we can
have%
\[
\mathbf{d}^{\ast}\succeq\mathbf{c}.
\]

Based on this, we need to study some useful relations among the three
policies: $\mathbf{d}$, $\mathbf{c}$ and $\mathbf{d}^{\ast}$, where
$\mathbf{d}^{\ast}$ is deterministic with $d_{i}^{\ast}=0$ for each
$i=1,2,\ldots,K-1,K$.

To compare $\eta^{\mathbf{c}}$ with $\eta^{\mathbf{d}^{\ast}}$, let $S\left(
\mathbf{d}^{\ast}\mathbf{,c}\right)  =\left\{  n_{l}:l=1,2,\ldots,n\right\}  $
for $1\leq n\leq K$. Then $c_{n_{l}}=1$ for $l=1,2,\ldots,n$, since
$d_{i}^{\ast}=0$ for each $i=1,2,\ldots,K-1,K$.

For the two policies $\mathbf{d}$ and $\mathbf{c}$, we have $d_{i},c_{i}%
\in\left\{  0,1\right\}  $. Further, for the three elements: $d_{i}$, $c_{i}$
and $d_{i}^{\ast}=0$ for $i\in S\left(  \mathbf{d}^{\ast}\mathbf{,c}\right)
$, we need to consider four different cases as follows:

\textit{Case one: }$d_{i}=c_{i}=d_{i}^{\ast}=0$. Since $c_{i}=d_{i}^{\ast}$,
this case does not require any analysis by using Lemma \ref{Lem:PerDE}.

\textit{Case two:} $d_{i}=1$ and $c_{i}=d_{i}^{\ast}=0$. Since $c_{i}%
=d_{i}^{\ast}$, this case does not require any analysis by using Lemma
\ref{Lem:PerDE}.

\textit{Case three: }$c_{i}=1$ and $d_{i}=d_{i}^{\ast}=0$. Note that
$c_{i}\neq d_{i}$, by using (1) of Theorem \ref{The:NP}, we obtain that
$G^{\left(  \mathbf{c}\right)  }\left(  i\right)  +b\leq0$. On the other hand,
since $c_{i}\neq d_{i}^{\ast}$, it follows from the performance difference
equation (\ref{Equa-2}) that for each $i\in S\left(  \mathbf{d}^{\ast
}\mathbf{,c}\right)  $,
\begin{align*}
\eta^{\mathbf{d}^{\ast}}-\eta^{\mathbf{c}}  &  =\mu_{2}\pi^{\left(
\mathbf{d}^{\ast}\right)  }\left(  i\right)  \left(  d_{i}^{\ast}%
-c_{i}\right)  \left[  {G}^{\left(  \mathbf{c}\right)  }\left(  i\right)
+{b}\right] \\
&  =-\mu_{2}\pi^{\left(  \mathbf{d}^{\ast}\right)  }\left(  i\right)  \left[
{G}^{\left(  \mathbf{c}\right)  }\left(  i\right)  +{b}\right]  \geq0.
\end{align*}
Thus $\eta^{\mathbf{d}^{\ast}}\geq\eta^{\mathbf{c}}$, this gives
$\mathbf{d}^{\ast}\succeq\mathbf{c}$.

\textit{Case four:} $d_{i}=c_{i}=1$ and $d_{i}^{\ast}=0$. Note that
$d_{i}^{\ast}\neq d_{i}$, by using (1) of Theorem \ref{The:NP}, we obtain that
$G^{\left(  \mathbf{d}^{\ast}\right)  }\left(  i\right)  +b\leq0$. On the
other hand, since $c_{i}\neq d_{i}^{\ast}$, it follows from the performance
difference equation (\ref{Equa-2}) that for each $i\in S\left(  \mathbf{d}%
^{\ast}\mathbf{,c}\right)  $,
\begin{align*}
\eta^{\mathbf{c}}-\eta^{\mathbf{d}^{\ast}}  &  =\mu_{2}\pi^{\left(
\mathbf{c}\right)  }\left(  i\right)  \left(  c_{i}-d_{i}^{\ast}\right)
\left[  {G}^{\left(  \mathbf{d}^{\ast}\right)  }\left(  i\right)  +{b}\right]
\\
&  =\mu_{2}\pi^{\left(  \mathbf{c}\right)  }\left(  i\right)  \left[
{G}^{\left(  \mathbf{d}^{\ast}\right)  }\left(  i\right)  +{b}\right]  \leq0.
\end{align*}
Thus $\eta^{\mathbf{d}^{\ast}}\geq\eta^{\mathbf{c}}$, this gives
$\mathbf{d}^{\ast}\succeq\mathbf{c}$.

Based on the above four discussions, we obtain that $\mathbf{d}^{\ast}%
\succeq\mathbf{c}$ for any policy $\mathbf{c}\in\mathcal{D}$. This completes
the proof. \textbf{{\rule{0.08in}{0.08in}}}

For $\mathbf{d}^{\ast}=\left(  0;0,0,\ldots,0;1,1,\ldots,1\right)  $, let
$\mathbf{d}^{\left(  n\right)  }$ be a policy in the policy set $\mathcal{D}$
with
\[
S\left(  \mathbf{d}^{\ast}\mathbf{,d}^{\left(  n\right)  }\right)  =\left\{
i_{l}:l=1,2,\ldots,n\right\}
\]
for $1\leq n\leq K$. To understand Policy $\mathbf{d}^{\left(  n\right)  }$,
we take three examples: $S\left(  \mathbf{d}^{\ast}\mathbf{,d}^{\left(
1\right)  }\right)  =\left\{  i_{1}\right\}  $, $S\left(  \mathbf{d}^{\ast
}\mathbf{,d}^{\left(  2\right)  }\right)  =\left\{  i_{1},i_{2}\right\}  $,
$S\left(  \mathbf{d}^{\ast}\mathbf{,d}^{\left(  3\right)  }\right)  =\left\{
i_{1},i_{2},i_{3}\right\}  $. Also, $S\left(  \mathbf{d}^{\left(  n-1\right)
}\mathbf{,d}^{\left(  n\right)  }\right)  =\left\{  i_{n}\right\}  $ for
$1\leq n\leq K$. Note that%
\[
\mathbf{d}^{\left(  K\right)  }==\left(  0;1,1,\ldots,1;1,1,\ldots,1\right)
.
\]

The following corollary provides a set-structured decreasing monotonicity of
the policies $\mathbf{d}^{\left(  n\right)  }\in\mathcal{D}$ for
$n=1,2,\ldots,K-1,K$. In fact, this monotonicity is guaranteed by the class
property of policies in the set $\mathcal{D}$, given in (1) of Theorem
\ref{The:NP}. The proof is easy by using a similar analysis to that in Theorem
\ref{The:left}, thus it is omitted here.

\begin{Cor}
If $P\geq{P}_{H}\left(  \mathbf{d}\right)  $ for any given policy $\mathbf{d}%
$, then%
\[
\mathbf{d}^{\ast}\succeq\mathbf{d}^{\left(  1\right)  }\succeq\mathbf{d}%
^{\left(  2\right)  }\succeq\mathbf{d}^{\left(  3\right)  }\succeq
\cdots\succeq\mathbf{d}^{\left(  K-1\right)  }\succeq\mathbf{d}^{\left(
K\right)  }.
\]
\end{Cor}

In what follows we compute the maximal long-run average profit of the
stock-rationing queue.

When $P\geq{P}_{H}\left(  \mathbf{d}\right)  $ for any given policy
$\mathbf{d}$, the optimal dynamic rationing policy is given by
\[
\mathbf{d}^{\ast}=\left(  0;0,0,\ldots,0;1,1,\ldots,1\right)  ,
\]
thus it follows from (\ref{5}) that%
\[%
\begin{array}
[c]{l}%
\text{ \ \ \ }\xi_{0}=1,\text{ \ \ \ \ \ \ \ \ \ \ \ \ \ \ }i=0,\\
\xi_{i}^{(\mathbf{d}^{\ast})}=\left\{
\begin{array}
[c]{ll}%
\alpha^{i}, & i=1,2,\ldots,K,\\
\left(  \frac{\alpha}{\beta}\right)  ^{K}\beta^{i}, & i=K+1,K+2,\ldots,N,
\end{array}
\right.
\end{array}
\]
and%
\[
h^{(\mathbf{d}^{\ast})}=1+\sum\limits_{i=1}^{N}\xi_{i}^{(\mathbf{d}^{\ast}%
)}=1+\frac{\alpha\left(  1-\alpha^{K}\right)  }{1-\alpha}+\left(
\frac{\alpha}{\beta}\right)  ^{K}\frac{\beta^{K+1}\left(  1-\beta
^{N-K}\right)  }{1-\beta},
\]
where $\alpha=\lambda/\mu_{1}$ and $\beta=\lambda/\left(  \mu_{1}+\mu
_{2}\right)  .$ It follows from (\ref{37}) that%
\[
\pi^{(\mathbf{d}^{\ast})}\left(  i\right)  =\left\{
\begin{array}
[c]{ll}%
\frac{1}{h^{(\mathbf{d}^{\ast})}}, & i=0,\\
\frac{1}{h^{(\mathbf{d}^{\ast})}}\xi_{i}^{(\mathbf{d}^{\ast})}, &
i=1,2,\ldots,N.
\end{array}
\right.
\]
At the same time, it follows from (\ref{8}) to (\ref{11}) that%
\[%
\begin{array}
[c]{ll}%
f\left(  0\right)  =-C_{2,1}\mu_{1}-C_{2,2}\mu_{2}-C_{3}\lambda, & i=0;\\
f^{(\mathbf{d}^{\ast})}\left(  i\right)  =R\mu_{1}-C_{1}i-C_{2,2}\mu_{2}%
-C_{3}\lambda, & 1\leq i\leq K;\\
f\left(  i\right)  =R\left(  \mu_{1}+\mu_{2}\right)  -C_{1}i-C_{3}%
\lambda1_{\left\{  i<N\right\}  }-C_{4}\lambda1_{\left\{  i=N\right\}  }, &
K+1\leq i\leq N.
\end{array}
\]
Since%
\begin{equation}
\eta^{\mathbf{d}^{\ast}}=\sum\limits_{i=0}^{N}\pi^{(\mathbf{d}^{\ast})}\left(
i\right)  f^{(\mathbf{d}^{\ast})}\left(  i\right)  , \label{35}%
\end{equation}
we obtain
\begin{align*}
\eta^{\mathbf{d}^{\ast}}  &  =\frac{1}{h^{(\mathbf{d}^{\ast})}}\left\{
-\left(  C_{2,1}\mu_{1}+C_{2,2}\mu_{2}+C_{3}\lambda\right)  +\sum
\limits_{i=0}^{K}\left(  R\mu_{1}-C_{1}i-C_{2,2}\mu_{2}-C_{3}\lambda\right)
\alpha^{i}\right. \\
&  \text{ \ \ \ }\left.  +\sum\limits_{i=K+1}^{N}\left[  R\left(  \mu_{1}%
+\mu_{2}\right)  -C_{1}i-C_{3}\lambda1_{\left\{  i<N\right\}  }-C_{4}%
\lambda1_{\left\{  i=N\right\}  }\right]  \left(  \frac{\alpha}{\beta}\right)
^{K}\beta^{i}\right\} \\
&  =\frac{1}{h^{(\mathbf{d}^{\ast})}}\left\{  -\gamma_{1}+\gamma
_{2}\frac{\alpha\left(  1-\alpha^{K}\right)  }{1-\alpha}-C_{1}\left[
\frac{\alpha\left(  1-\alpha^{K}\right)  }{\left(  1-\alpha\right)  ^{2}%
}-\frac{K\alpha^{K+1}}{1-\alpha}\right]  \right. \\
&  \text{ \ \ \ }\left.  +\left(  \frac{\alpha}{\beta}\right)  ^{K}\gamma
_{3}\frac{\beta^{K+1}\left(  1-\beta^{N-K}\right)  }{1-\beta}\right. \\
&  \text{ \ \ \ }\left.  -\left(  \frac{\alpha}{\beta}\right)  ^{K}%
C_{1}\left[  \frac{K\beta^{K+1}-N\beta^{N+1}}{1-\beta}+\frac{\beta
^{K+1}\left(  1-\beta^{N-K}\right)  }{\left(  1-\beta\right)  ^{2}}\right]
\right\}  ,
\end{align*}
where%
\begin{align*}
\gamma_{1}  &  =C_{2,1}\mu_{1}+C_{2,2}\mu_{2}+C_{3}\lambda,\\
\gamma_{2}  &  =R\mu_{1}-C_{2,2}\mu_{2}-C_{3}\lambda,\\
\gamma_{3}  &  =R\left(  \mu_{1}+\mu_{2}\right)  -C_{3}\lambda1_{\left\{
i<N\right\}  }-C_{4}\lambda1_{\left\{  i=N\right\}  }.
\end{align*}

\subsection{The penalty cost $P_{L}\left(  \mathbf{d}\right)  >0$ and $0\leq
P\leq P_{L}\left(  \mathbf{d}\right)  $}

In this subsection, we consider the area of the penalty cost: $P_{L}\left(
\mathbf{d}\right)  >0$ and $0\leq P\leq P_{L}\left(  \mathbf{d}\right)  $ for
any given policy $\mathbf{d}$. We first find the optimal dynamic rationing
policy of the stock-rationing queue. Then we compute the maximal long-run
average profit of this system.

The following theorem finds the optimal dynamic rationing policy of the
stock-rationing queue in the area of the penalty cost: $P_{L}\left(
\mathbf{d}\right)  >0$ and $0\leq P\leq P_{L}\left(  \mathbf{d}\right)  $ for
any given policy $\mathbf{d}$. The proof is similar to that of Theorem
\ref{The:left}.

\begin{The}
\label{The:right}If $P_{L}\left(  \mathbf{d}\right)  >0$ and $0\leq P\leq
P_{L}\left(  \mathbf{d}\right)  $ for any given policy $\mathbf{d}$, then the
optimal dynamic rationing policy of the stock-rationing queue is given by
\[
\mathbf{d}^{\ast}=\left(  0;1,1,\ldots,1;1,1,\ldots,1\right)  .
\]
This shows that if the penalty cost is lower with $P_{L}\left(  \mathbf{d}%
\right)  >0$ and $0\leq P\leq{P}_{L}\left(  \mathbf{d}\right)  $, then the
warehouse would like to supply the products to the demands of Class $2$.
\end{The}

\textbf{Proof: }If $P_{L}\left(  \mathbf{d}\right)  >0$ and $0\leq P\leq
P_{L}\left(  \mathbf{d}\right)  $ for any given policy $\mathbf{d}$, then we
prove that for any policy $\mathbf{c}\in\mathcal{D}$,%
\[
\mathbf{d}^{\ast}\succeq\mathbf{c}.
\]

For this, we need to consider the three policies: $\mathbf{d}$, $\mathbf{c}$
and $\mathbf{d}^{\ast}$, where $\mathbf{d}^{\ast}$ is deterministic with
$d_{i}^{\ast}=1$ for each $i=1,2,\ldots,K-1,K$.

To compare $\eta^{\mathbf{c}}$ with $\eta^{\mathbf{d}^{\ast}}$, let $S\left(
\mathbf{d}^{\ast}\mathbf{,c}\right)  =\left\{  n_{l}:l=1,2,\ldots,m\right\}  $
for $1\leq m\leq K$. Then $c_{n_{l}}=0$ for $l=1,2,\ldots,m$, since
$d_{i}^{\ast}=1$ for each $i=1,2,\ldots,K-1,K$.

For the two policies $\mathbf{d}$ and $\mathbf{c}$, we have $d_{i},c_{i}%
\in\left\{  0,1\right\}  $. Based on this, for the three elements: $d_{i}$,
$c_{i}$ and $d_{i}^{\ast}=1$ for $i\in S\left(  \mathbf{d}^{\ast}%
\mathbf{,c}\right)  $, we need to consider four different cases as follows:

\textit{Case one: }$d_{i}=c_{i}=d_{i}^{\ast}=1$. Since $c_{i}=d_{i}^{\ast}$,
this case does not require any analysis according to Lemma \ref{Lem:PerDE}.

\textit{Case two:} $d_{i}=0$ and $c_{i}=d_{i}^{\ast}=1$. Since $c_{i}%
=d_{i}^{\ast}$, this case does not require any analysis by using Lemma
\ref{Lem:PerDE}.

\textit{Case three:} $c_{i}=0$ and $d_{i}=d_{i}^{\ast}=1$. Note that
$c_{i}\neq d_{i}$, by using (2) of Theorem \ref{The:NP}, we obtain that
$G^{\left(  \mathbf{c}\right)  }\left(  i\right)  +b\geq0$. On the other hand,
since $c_{i}\neq d_{i}^{\ast}$, it follows from the performance difference
equation (\ref{Equa-2}) that for each $i\in S\left(  \mathbf{d}^{\ast
}\mathbf{,c}\right)  $,
\begin{align*}
\eta^{\mathbf{d}^{\ast}}-\eta^{\mathbf{c}}  &  =\mu_{2}\pi^{\left(
\mathbf{d}^{\ast}\right)  }\left(  i\right)  \left(  d_{i}^{\ast}%
-c_{i}\right)  \left[  {G}^{\left(  \mathbf{c}\right)  }\left(  i\right)
+{b}\right] \\
&  =\mu_{2}\pi^{\left(  \mathbf{d}^{\ast}\right)  }\left(  i\right)  \left[
{G}^{\left(  \mathbf{c}\right)  }\left(  i\right)  +{b}\right]  \geq0.
\end{align*}
Thus $\eta^{\mathbf{d}^{\ast}}\geq\eta^{\mathbf{c}}$, this gives
$\mathbf{d}^{\ast}\succeq\mathbf{c}$.

\textit{Case four:} $d_{i}=c_{i}=0$ and $d_{i}^{\ast}=1$. Note that
$d_{i}^{\ast}\neq d_{i}$, by using (2) of Theorem \ref{The:NP}, we obtain that
$G^{\left(  \mathbf{d}^{\ast}\right)  }\left(  i\right)  +b\geq0$. On the
other hand, since $c_{i}\neq d_{i}^{\ast}$, it follows from the performance
difference equation (\ref{Equa-2}) that for each $i\in S\left(  \mathbf{d}%
^{\ast}\mathbf{,c}\right)  $,
\begin{align*}
\eta^{\mathbf{c}}-\eta^{\mathbf{d}^{\ast}}  &  =\mu_{2}\pi^{\left(
\mathbf{c}\right)  }\left(  i\right)  \left(  c_{i}-d_{i}^{\ast}\right)
\left[  {G}^{\left(  \mathbf{d}^{\ast}\right)  }\left(  i\right)  +{b}\right]
\\
&  =-\mu_{2}\pi^{\left(  \mathbf{c}\right)  }\left(  i\right)  \left[
{G}^{\left(  \mathbf{d}^{\ast}\right)  }\left(  i\right)  +{b}\right]  \leq0.
\end{align*}
Thus $\eta^{\mathbf{d}^{\ast}}\geq\eta^{\mathbf{c}}$, this gives
$\mathbf{d}^{\ast}\succeq\mathbf{c}$. This completes the proof.
\textbf{{\rule{0.08in}{0.08in}}}

For $\mathbf{d}^{\ast}=\left(  0;1,1,\ldots,1;1,1,\ldots,1\right)  $, let
$\mathbf{d}^{\left(  n\right)  }$ be a policy in the policy set $\mathcal{D}$
with $S\left(  \mathbf{d}^{\ast}\mathbf{,d}^{\left(  n\right)  }\right)
=\left\{  k_{l}:l=1,2,\ldots,n\right\}  $ for $1\leq n\leq K$, where%
\[
\widetilde{\mathbf{d}}^{\left(  K\right)  }=\left(  0;0,0,\ldots
,0;1,1,\ldots,1\right)  .
\]

The following corollary provides a set-structured decreasing monotonicity of
the policies $\mathbf{d}^{\left(  n\right)  }\in\mathcal{D}$ for
$n=1,2,\ldots,K-1,K$. This monotonicity comes from the class property of the
policies in the set $\mathcal{D}$, given in (2) of Theorem \ref{The:NP}. The
proof is easy and omitted here.

\begin{Cor}
If $P_{L}\left(  \mathbf{d}\right)  >0$ and $0\leq P\leq P_{L}\left(
\mathbf{d}\right)  $ for any given policy $\mathbf{d}$, then%
\[
\mathbf{d}^{\ast}\succeq\mathbf{d}^{\left(  1\right)  }\succeq\mathbf{d}%
^{\left(  2\right)  }\succeq\mathbf{d}^{\left(  3\right)  }\succeq
\cdots\succeq\mathbf{d}^{\left(  K-1\right)  }\succeq\mathbf{d}^{\left(
K\right)  }.
\]
\end{Cor}

If $P_{L}\left(  \mathbf{d}\right)  >0$ and $0\leq P\leq P_{L}\left(
\mathbf{d}\right)  $ for any given policy $\mathbf{d}$, then the optimal
dynamic rationing policy is given by
\[
\mathbf{d}^{\ast}=\left(  0;1,1,\ldots,1;1,1,\ldots,1\right)  .
\]
In this case, we obtain%
\[%
\begin{array}
[c]{cl}%
\xi_{0}=1, & i=0,\\
\xi_{i}^{(\mathbf{d}^{\ast})}=\beta^{i}, & i=1,2,\ldots,N,
\end{array}
\]
and%
\[
h^{(\mathbf{d}^{\ast})}=1+\sum\limits_{i=1}^{N}\xi_{i}^{(\mathbf{d}^{\ast}%
)}=1+\frac{\beta\left(  1-\beta^{N}\right)  }{1-\beta}.
\]
It follows from Subsection 1.1.4 of Chapter 1 in Li \cite{Li:2010} that%
\[
\pi^{(\mathbf{d}^{\ast})}\left(  i\right)  =\left\{
\begin{array}
[c]{ll}%
\frac{1}{h^{(\mathbf{d}^{\ast})}} & i=0,\\
\frac{\beta^{i}}{h^{(\mathbf{d}^{\ast})}}, & i=1,2,\ldots,N,
\end{array}
\right.
\]
At the same time, it follows from (\ref{8}) to (\ref{11}) that
\[%
\begin{array}
[c]{ll}%
f\left(  0\right)  =-C_{2,1}\mu_{1}-C_{2,2}\mu_{2}-C_{3}\lambda, & i=0;\\
f^{(\mathbf{d}^{\ast})}\left(  i\right)  =R\left(  \mu_{1}+\mu_{2}\right)
-C_{1}i-C_{3}\lambda-P\mu_{2}, & 1\leq i\leq K;\\
f^{(\mathbf{d}^{\ast})}\left(  i\right)  =R\left(  \mu_{1}+\mu_{2}\right)
-C_{1}i-C_{3}\lambda1_{\left\{  i<N\right\}  }-C_{4}\lambda1_{\left\{
i=N\right\}  }, & K+1\leq i\leq N.
\end{array}
\]
Thus we obtain%
\begin{align*}
\eta^{\mathbf{d}^{\ast}}  &  =\frac{1}{h^{(\mathbf{d}^{\ast})}}\left\{
-\left(  C_{2,1}\mu_{1}+C_{2,2}\mu_{2}+C_{3}\lambda\right)  +\sum
\limits_{i=1}^{K}\left[  R\left(  \mu_{1}+\mu_{2}\right)  -C_{1}i-C_{3}%
\lambda-P\mu_{2}\right]  \beta^{i}\right. \\
\text{ \ \ }  &  \left.  \text{ \ }+\sum\limits_{i=K+1}^{N}\left[  R\left(
\mu_{1}+\mu_{2}\right)  -C_{1}i-C_{3}\lambda1_{\left\{  i<N\right\}  }%
-C_{4}\lambda1_{\left\{  i=N\right\}  }\right]  \beta^{i}\right\} \\
&  =\frac{1}{h^{(\mathbf{d}^{\ast})}}\left\{  -\gamma_{1}+\gamma
_{4}\frac{\beta\left(  1-\beta^{K}\right)  }{1-\beta}-C_{1}\left[
\frac{\beta\left(  1-\beta^{K}\right)  }{\left(  1-\beta\right)  ^{2}%
}-\frac{K\beta^{K+1}}{1-\beta}\right]  \right. \\
&  \left.  \text{ \ }+\gamma_{3}\frac{\beta^{K+1}\left(  1-\beta^{N-K}\right)
}{1-\beta}-C_{1}\left[  \frac{K\beta^{K+1}-N\beta^{N+1}}{1-\beta}%
+\frac{\beta^{K+1}\left(  1-\beta^{N-K}\right)  }{\left(  1-\beta\right)
^{2}}\right]  \right\}  ,
\end{align*}
where%
\begin{align*}
\gamma_{1}  &  =C_{2,1}\mu_{1}+C_{2,2}\mu_{2}+C_{3}\lambda,\\
\gamma_{2}  &  =R\mu_{1}-C_{2,2}\mu_{2}-C_{3}\lambda,\\
\gamma_{3}  &  =R\left(  \mu_{1}+\mu_{2}\right)  -C_{3}\lambda1_{\left\{
i<N\right\}  }-C_{4}\lambda1_{\left\{  i=N\right\}  },\\
\gamma_{4}  &  =R\left(  \mu_{1}+\mu_{2}\right)  -C_{3}\lambda-P\mu_{2}.
\end{align*}

\subsection{The penalty cost $P_{L}\left(  \mathbf{d}\right)  <P<P_{H}\left(
\mathbf{d}\right)  $}

In this subsection, we discuss the third area of the penalty cost:
$P_{L}<P<P_{H}$ for any given policy $\mathbf{d}$. Note that this analysis is
a little more complicated than those in the previous two areas. To this end,
we propose a new algebraic method to find the optimal dynamic rationing policy
of the stock-rationing queue. Based on this, we show that the optimal dynamic
rationing policy may not be of threshold type, but it must be of
transformational threshold type.

For the convenience of readers, it is necessary and useful to simply recall
several previous results as follows.

For any given policy $\mathbf{d}=\left(  0;d_{1},d_{2},\ldots,\ldots
,d_{K-1},d_{K};1,1,\ldots,1\right)  $, the unique solution of the linear
equation $G^{\left(  \mathbf{d}\right)  }\left(  i\right)  +b=0$ in the
penalty cost $P$ is given by
\[
\mathfrak{P}_{i}^{\left(  \mathbf{d}\right)  }=\frac{R+C_{2,2}+\lambda
^{-i}\left[  B_{0}-{D}^{\left(  \mathbf{d}\right)  }\right]  \prod
\limits_{k=1}^{i-1}v\left(  d_{k}\right)  +\sum\limits_{r=1}^{i-1}%
\lambda^{r-i}\left[  B_{r}^{\left(  \mathbf{d}\right)  }-{D}^{\left(
\mathbf{d}\right)  }\right]  \prod\limits_{k=r+1}^{i-1}v\left(  d_{k}\right)
}{1+\lambda^{-i}\left[  A_{0}-{F}^{\left(  \mathbf{d}\right)  }\right]
\prod\limits_{k=1}^{i-1}v\left(  d_{k}\right)  +\sum\limits_{r=1}^{i-1}%
\lambda^{r-i}\left[  A_{r}^{\left(  \mathbf{d}\right)  }-{F}^{\left(
\mathbf{d}\right)  }\right]  \prod\limits_{k=r+1}^{i-1}v\left(  d_{k}\right)
},
\]
which is a fixed real number for $1\leq i\leq K$.

We introduce a convention: If $\mathfrak{P}_{n-1}^{\left(  \mathbf{d}\right)
}<\mathfrak{P}_{n}^{\left(  \mathbf{d}\right)  }=\mathfrak{P}_{n+1}^{\left(
\mathbf{d}\right)  }=\cdots=\mathfrak{P}_{n+i}^{\left(  \mathbf{d}\right)
}=c$ and $\mathfrak{P}_{n-1}^{\left(  \mathbf{d}\right)  }<P\leq c$, then we
write%
\[
\mathfrak{P}_{n-1}^{\left(  \mathbf{d}\right)  }<P\leq\mathfrak{P}%
_{n}^{\left(  \mathbf{d}\right)  }=\mathfrak{P}_{n+1}^{\left(  \mathbf{d}%
\right)  }=\cdots=\mathfrak{P}_{n+i}^{\left(  \mathbf{d}\right)  },
\]
that is, the penalty cost $P$ is written in front of all the equal elements in
the sequence $\left\{  \mathfrak{P}_{k}^{\left(  \mathbf{d}\right)  }:n\leq
k\leq n+i\right\}  $.

For the sequence $\left\{  \mathfrak{P}_{k}^{\left(  \mathbf{d}\right)
}:1\leq k\leq K\right\}  $, we set up a new permutation from the smallest to
the largest as follows:%
\[
\mathfrak{P}_{i_{1}}^{\left(  \mathbf{d}\right)  }\leq\mathfrak{P}_{i_{2}%
}^{\left(  \mathbf{d}\right)  }\leq\cdots\leq\mathfrak{P}_{i_{K-1}}^{\left(
\mathbf{d}\right)  }\leq\mathfrak{P}_{i_{K}}^{\left(  \mathbf{d}\right)  },
\]
it is clear that $\mathfrak{P}_{i_{1}}^{\left(  \mathbf{d}\right)  }%
=P_{L}\left(  \mathbf{d}\right)  $ and $\mathfrak{P}_{i_{K}}^{\left(
\mathbf{d}\right)  }=P_{H}\left(  \mathbf{d}\right)  $. For the convenience of
description, for the incremental sequence $\left\{  \mathfrak{P}_{i_{j}%
}^{\left(  \mathbf{d}\right)  }:1\leq j\leq K\right\}  $, we write its
subscript vector as $\left(  i_{1},i_{2},\ldots,i_{K-1},i_{K}\right)  $. Note
that the subscript vector $\left(  i_{1},i_{2},\ldots,i_{K-1},i_{K}\right)  $
depends on Policy $\mathbf{d}$.

The following lemma shows how the penalty cost $P$ is distributed in the
sequence $\left\{  \mathfrak{P}_{k}^{\left(  \mathbf{d}\right)  }:1\leq k\leq
K\right\}  $.

\begin{Lem}
If $P_{L}\left(  \mathbf{d}\right)  <P<P_{H}\left(  \mathbf{d}\right)  $ for
any given policy $\mathbf{d}$, then there exists the minimal positive integer
$n_{0}\in\left\{  1,2,\ldots,K\right\}  $ such that either%
\[
\mathfrak{P}_{i_{n_{0}}}^{\left(  \mathbf{d}\right)  }<P=\mathfrak{P}%
_{i_{n_{0}+1}}^{\left(  \mathbf{d}\right)  }%
\]
or%
\[
\mathfrak{P}_{i_{n_{0}}}^{\left(  \mathbf{d}\right)  }<P<\mathfrak{P}%
_{i_{n_{0}+1}}^{\left(  \mathbf{d}\right)  }.
\]
\end{Lem}

\textbf{Proof: }Note that%
\[
P_{H}\left(  \mathbf{d}\right)  =\max_{\mathbf{d}\in\mathcal{D}}\left\{
0,\mathfrak{P}_{1}^{\left(  \mathbf{d}\right)  },\mathfrak{P}_{2}^{\left(
\mathbf{d}\right)  },\ldots,\mathfrak{P}_{K}^{\left(  \mathbf{d}\right)
}\right\}
\]
and%
\[
P_{L}\left(  \mathbf{d}\right)  =\max_{\mathbf{d}\in\mathcal{D}}\left\{
\mathfrak{P}_{1}^{\left(  \mathbf{d}\right)  },\mathfrak{P}_{2}^{\left(
\mathbf{d}\right)  },\ldots,\mathfrak{P}_{K}^{\left(  \mathbf{d}\right)
}\right\}  ,
\]
it is easy to see that $P_{H}\left(  \mathbf{d}\right)  $ and $P_{L}\left(
\mathbf{d}\right)  $ are two fixed real numbers. If $P_{L}<P<P_{H}\left(
\mathbf{d}\right)  $ for Policy $\mathbf{d}$, then there exists the minimal
positive integer $n_{0}\in\left\{  1,2,\ldots,K-1,K\right\}  $ such that%
\[
P_{L}\left(  \mathbf{d}\right)  \leq\mathfrak{P}_{i_{n_{0}}}^{\left(
\mathbf{d}\right)  }<P\leq\mathfrak{P}_{i_{n_{0}+1}}^{\left(  \mathbf{d}%
\right)  }<P_{H}\left(  \mathbf{d}\right)  .
\]
This shows that either for $P=\mathfrak{P}_{i_{n_{0}+1}}^{\left(
\mathbf{d}\right)  },$
\[
\mathfrak{P}_{i_{n_{0}}}^{\left(  \mathbf{d}\right)  }<P=\mathfrak{P}%
_{i_{n_{0}+1}}^{\left(  \mathbf{d}\right)  };
\]
or for $P<\mathfrak{P}_{i_{n_{0}+1}}^{\left(  \mathbf{d}\right)  },$%
\[
\mathfrak{P}_{i_{n_{0}}}^{\left(  \mathbf{d}\right)  }<P<\mathfrak{P}%
_{i_{n_{0}+1}}^{\left(  \mathbf{d}\right)  }.
\]
This completes the proof. \textbf{{\rule{0.08in}{0.08in}}}

Now, our task is to develop a new method for finding the optimal dynamic
rationing policy by means of the following two useful information: (a) The
incremental sequence
\[
P_{L}\left(  \mathbf{d}\right)  =\mathfrak{P}_{i_{1}}^{\left(  \mathbf{d}%
\right)  }\leq\mathfrak{P}_{i_{2}}^{\left(  \mathbf{d}\right)  }\leq\cdots
\leq\mathfrak{P}_{i_{K-1}}^{\left(  \mathbf{d}\right)  }\leq\mathfrak{P}%
_{i_{K}}^{\left(  \mathbf{d}\right)  }=P_{H}\left(  \mathbf{d}\right)  ;
\]
and (b) the penalty cost $P$ has a fixed position: $\mathfrak{P}_{i_{n_{0}}%
}^{\left(  \mathbf{d}\right)  }<P\leq\mathfrak{P}_{i_{n_{0}+1}}^{\left(
\mathbf{d}\right)  }$, where $n_{0}$ is the minimal positive integer in the
set $\left\{  1,2,\ldots,K-1,K\right\}  $.

In what follows we discuss two different cases: A simple case and a general case.

\textbf{Case one: A simple case with}%
\begin{equation}
P_{L}\left(  \mathbf{d}\right)  =\mathfrak{P}_{1}^{\left(  \mathbf{d}\right)
}\leq\mathfrak{P}_{2}^{\left(  \mathbf{d}\right)  }\leq\cdots\leq
\mathfrak{P}_{K-1}^{\left(  \mathbf{d}\right)  }\leq\mathfrak{P}_{K}^{\left(
\mathbf{d}\right)  }=P_{H}\left(  \mathbf{d}\right)  . \label{Equa-5}%
\end{equation}
In this case, the subscript vector is expressed as $\left\{  1,2,3,\ldots
,K-1,K\right\}  $ depending on Policy $\mathbf{d}$.

If $P_{L}\left(  \mathbf{d}\right)  <P<P_{H}\left(  \mathbf{d}\right)  $ for
any given policy $\mathbf{d}$, then there exists the minimal positive integer
$n_{0}\in\left\{  1,2,\ldots,K-1,K\right\}  $ such that%
\[
P_{L}\left(  \mathbf{d}\right)  =\mathfrak{P}_{1}^{\left(  \mathbf{d}\right)
}\leq\cdots\leq\mathfrak{P}_{n_{0}-1}^{\left(  \mathbf{d}\right)  }%
<P\leq\mathfrak{P}_{n_{0}}^{\left(  \mathbf{d}\right)  }\leq\cdots
\leq\mathfrak{P}_{K}^{\left(  \mathbf{d}\right)  }=P_{H}\left(  \mathbf{d}%
\right)  .
\]
Based on this, we take two different sets%
\[
\Lambda_{1}=\left\{  \mathfrak{P}_{1}^{\left(  \mathbf{d}\right)
},\mathfrak{P}_{2}^{\left(  \mathbf{d}\right)  },\ldots,\mathfrak{P}_{n_{0}%
-1}^{\left(  \mathbf{d}\right)  }\right\}
\]
and%
\[
\Lambda_{2}=\left\{  \mathfrak{P}_{n_{0}}^{\left(  \mathbf{d}\right)
},\mathfrak{P}_{n_{0}+1}^{\left(  \mathbf{d}\right)  },\ldots,\mathfrak{P}%
_{K}^{\left(  \mathbf{d}\right)  }\right\}  .
\]
By using the two sets $\Lambda_{1}$ and $\Lambda_{2}$, we write%
\[
\overline{P}_{H}\left(  \mathbf{d;}1\rightarrow n_{0}-1\right)  =\max_{1\leq
i\leq n_{0}-1}\left\{  \mathfrak{P}_{i}^{\left(  \mathbf{d}\right)
}\right\}
\]
and%
\[
\overline{P}_{L}\left(  \mathbf{d;}n_{0}\rightarrow K\right)  =\min_{n_{0}\leq
j\leq K}\left\{  \mathfrak{P}_{j}^{\left(  \mathbf{d}\right)  }\right\}
\text{.}%
\]
It is clear that $\overline{P}_{H}\left(  \mathbf{d;}1\rightarrow
n_{0}-1\right)  =\mathfrak{P}_{n_{0}-1}^{\left(  \mathbf{d}\right)  }$ and
$\overline{P}_{L}\left(  \mathbf{d;}n_{0}\rightarrow K\right)  =\mathfrak{P}%
_{n_{0}}^{\left(  \mathbf{d}\right)  }$.

For this simple case, the following theorem finds the optimal dynamic
rationing policy, which is of threshold type.

\begin{The}
\label{The:Opt}For the simple case with $P_{L}\left(  \mathbf{d}\right)
<P<P_{H}\left(  \mathbf{d}\right)  $ for any given policy $\mathbf{d}$, if
there exists the minimal positive integer $n_{0}\in\left\{  1,2,\ldots
,K-1,K\right\}  $ such that%
\begin{equation}
P_{L}\left(  \mathbf{d}\right)  =\mathfrak{P}_{1}^{\left(  \mathbf{d}\right)
}\leq\cdots\leq\mathfrak{P}_{n_{0}-1}^{\left(  \mathbf{d}\right)  }%
<P\leq\mathfrak{P}_{n_{0}}^{\left(  \mathbf{d}\right)  }\leq\cdots
\leq\mathfrak{P}_{K}^{\left(  \mathbf{d}\right)  }=P_{H}\left(  \mathbf{d}%
\right)  , \label{Equa-6}%
\end{equation}
then the optimal dynamic rationing policy is given by%
\[
\mathbf{d}^{\ast}=\left(  0;\underset{n_{0}-1\text{ zeros}}{\underbrace
{0,0,\ldots,0}},\underset{K-n_{0}+1\text{ ones}}{\underbrace{1,1,\ldots,1}%
};1,1,\ldots,1\right)  .
\]
\end{The}

\textbf{Proof: }The proof follows that in Theorems \ref{The:left} and
\ref{The:right}.

On the one hand, in the set $\Lambda_{1}$, it is easy to see from
(\ref{Equa-6}) that $P>\overline{P}_{H}\left(  \mathbf{d;}1\rightarrow
n_{0}-1\right)  $ for Policy $\mathbf{d}$. Now, our aim is to focus on a
sub-policy%
\[
\widetilde{\mathbf{d}}_{a}=\left(  0;d_{1},d_{2},\ldots,d_{n_{0}-1},\ast
,\ast,\ldots,\ast;1,1,\ldots,1\right)  .
\]
For the sub-policy $\left(  d_{1},d_{2},\ldots,d_{n_{0}-1}\right)  $, it is
easy to see from the set $\Lambda_{1}$ that $P>\overline{P}_{H}\left(
\mathbf{d;}1\rightarrow n_{0}-1\right)  $. Thus it follows from Theorem
\ref{The:left} that the optimal dynamic rationing sub-policy is given by%
\[
\widetilde{\mathbf{d}}_{a}^{\ast}=\left(  0;0,0,\ldots,0,\ast,\ast,\ldots
,\ast;1,1,\ldots,1\right)  .
\]

On the other hand, it is seen from the set $\Lambda_{2}$ that $0\leq
P\leq\overline{P}_{L}\left(  \mathbf{d;}n_{0}\rightarrow K\right)  $ for
Policy $\mathbf{d}$. We consider another sub-policy%
\[
\widetilde{\mathbf{d}}_{b}=\left(  0;\ast,\ast,\ldots,\ast,d_{n_{0}}%
,d_{n_{0}+1},\ldots,d_{K};1,1,\ldots,1\right)  .
\]
For the sub-policy $\left(  d_{n_{0}},d_{n_{0}+1},\ldots,d_{K}\right)  $, it
is easy to see from the set $\Lambda_{2}$ that $0\leq P\leq\overline{P}%
_{L}\left(  \mathbf{d;}n_{0}\rightarrow K\right)  $. Thus it is easy to see
from Theorem \ref{The:right} that the optimal dynamic rationing sub-policy is
given by%
\[
\widetilde{\mathbf{d}}_{b}^{\ast}=\left(  0;\ast,\ast,\ldots,\ast
,1,1,\ldots,1;1,1,\ldots,1\right)  .
\]

Based on the above two discussions, from the total set $\Lambda_{1}\cup
\Lambda_{2}$, by observing the total policy $\left(  d_{1},d_{2}%
,\ldots,d_{n_{0}-1};d_{n_{0}},d_{n_{0}+1},\ldots,d_{K}\right)  $ or Policy
$\mathbf{d}$, the optimal dynamic rationing policy is given by%
\[
\mathbf{d}^{\ast}=\widetilde{\left(  \widetilde{\mathbf{d}}_{a}^{\ast}\right)
}_{b}^{\ast}=\widetilde{\left(  \widetilde{\mathbf{d}}_{b}^{\ast}\right)
}_{a}^{\ast}=\left(  0;\underset{n_{0}-1\text{ zeros}}{\underbrace
{0,0,\ldots,0}},\underset{K-n_{0}+1\text{ ones}}{\underbrace{1,1,\ldots,1}%
};1,1,\ldots,1\right)  .
\]
This completes the proof. \textbf{{\rule{0.08in}{0.08in}}}

\begin{Rem}
It is easy to see that in Theorems \ref{The:left}, \ref{The:right} and
\ref{The:Opt}, the optimal dynamic rationing policy is of threshold type
(i.e., critical rationing level).
\end{Rem}

\textbf{Case two: A general case with}%
\[
P_{L}\left(  \mathbf{d}\right)  =\mathfrak{P}_{i_{1}}^{\left(  \mathbf{d}%
\right)  }\leq\mathfrak{P}_{i_{2}}^{\left(  \mathbf{d}\right)  }\leq\cdots
\leq\mathfrak{P}_{i_{K-1}}^{\left(  \mathbf{d}\right)  }\leq\mathfrak{P}%
_{i_{K}}^{\left(  \mathbf{d}\right)  }=P_{H}\left(  \mathbf{d}\right)  .
\]
For the incremental sequence $\left\{  \mathfrak{P}_{i_{j}}^{\left(
\mathbf{d}\right)  }:j=1,2,\ldots,K\right\}  $, we write its subscript vector
as $\left(  i_{1},i_{2},\ldots,i_{K-1},i_{K}\right)  $, which depends on
Policy $\mathbf{d}$. In the general case, we assume that $\left(  i_{1}%
,i_{2},\ldots,i_{K-1},i_{K}\right)  \neq\left(  1,2,\ldots,K-1,K\right)  $.

If $P_{L}\left(  \mathbf{d}\right)  <P<P_{H}\left(  \mathbf{d}\right)  $ for
any given policy $\mathbf{d}$, then there exists the minimal positive integer
$n_{0}\in\left\{  1,2,\ldots,K-1,K\right\}  $ such that
\[
P_{L}\left(  \mathbf{d}\right)  =\mathfrak{P}_{i_{1}}^{\left(  \mathbf{d}%
\right)  }\leq\cdots\leq\mathfrak{P}_{i_{n_{0}-1}}^{\left(  \mathbf{d}\right)
}<P\leq\mathfrak{P}_{i_{n_{0}}}^{\left(  \mathbf{d}\right)  }\leq\cdots
\leq\mathfrak{P}_{i_{K}}^{\left(  \mathbf{d}\right)  }=P_{H}\left(
\mathbf{d}\right)  .
\]
Based on this, we take two sets%
\[
\Lambda_{1}^{G}=\left\{  \mathfrak{P}_{i_{1}}^{\left(  \mathbf{d}\right)
},\mathfrak{P}_{i_{2}}^{\left(  \mathbf{d}\right)  },\ldots,\mathfrak{P}%
_{i_{n_{0}-1}}^{\left(  \mathbf{d}\right)  }\right\}
\]
and%
\[
\Lambda_{2}^{G}=\left\{  \mathfrak{P}_{i_{n_{0}}}^{\left(  \mathbf{d}\right)
},\mathfrak{P}_{i_{n_{0}+1}}^{\left(  \mathbf{d}\right)  },\ldots
,\mathfrak{P}_{i_{K}}^{\left(  \mathbf{d}\right)  }\right\}  .
\]
For the two sets $\Lambda_{1}^{G}$ and $\Lambda_{2}^{G}$, we write%
\[
\overline{P}_{H}^{G}\left(  \mathbf{d;}1\rightarrow n_{0}-1\right)
=\max_{1\leq k\leq n_{0}-1}\left\{  \mathfrak{P}_{i_{k}}^{\left(
\mathbf{d}\right)  }\right\}
\]
and%
\[
\overline{P}_{L}^{G}\left(  \mathbf{d;}n_{0}\rightarrow K\right)  =\min
_{n_{0}\leq k\leq K}\left\{  \mathfrak{P}_{i_{k}}^{\left(  \mathbf{d}\right)
}\right\}  \text{,}%
\]
It is clear that $\overline{P}_{H}^{G}\left(  \mathbf{d;}1\rightarrow
n_{0}-1\right)  =\mathfrak{P}_{i_{n_{0}-1}}^{\left(  \mathbf{d}\right)  }$ and
$\overline{P}_{L}^{G}\left(  \mathbf{d;}n_{0}\rightarrow K\right)
=\mathfrak{P}_{i_{n_{0}}}^{\left(  \mathbf{d}\right)  }$.

Corresponding to the subscript vector of the incremental sequence $\left\{
\mathfrak{P}_{i_{k}}^{\left(  \mathbf{d}\right)  }:1\leq k\leq K\right\}  $,
we transfer Policy%
\[
\mathbf{d}=\left(  0;d_{1},d_{2},\ldots,d_{n_{0}-1},d_{n_{0}},d_{n_{0}%
+1},\ldots,d_{K};1,1,\ldots,1\right)
\]
into a new transformational policy%
\[
\mathbf{d}\left(  \text{Transfer}\right)  =\left(  0;d_{i_{1}},d_{i_{2}%
},\ldots,d_{i_{n_{0}-1}},d_{i_{n_{0}}},d_{i_{n_{0}+1}},\ldots,d_{i_{K}%
};1,1,\ldots,1\right)  .
\]
Therefore, a transformation of the optimal dynamic policy $\mathbf{d}^{\ast}$
is
\[
\left(  1,2,\ldots,K-1,K\right)  \Rightarrow\left(  i_{1},i_{2},\ldots
,i_{K-1},i_{K}\right)  ;
\]
and an inverse transformation of the optimal transformational dynamic policy
$\mathbf{d}^{\ast}\left(  \text{Transfer}\right)  $ is
\[
\left(  i_{1},i_{2},\ldots,i_{K-1},i_{K}\right)  \Rightarrow\left(
1,2,\ldots,K-1,K\right)  .
\]

For the general case, the following theorem finds the optimal dynamic
rationing policy, which may not be of threshold type, but must be of
transformational threshold type.

\begin{The}
\label{The:Trans}For the general case with $P_{L}\left(  \mathbf{d}\right)
<P<P_{H}\left(  \mathbf{d}\right)  $ for any given policy $\mathbf{d}$, if
there exists the minimal positive integer $n_{0}\in\left\{  1,2,\ldots
,K-1,K\right\}  $ such that
\[
P_{L}\left(  \mathbf{d}\right)  =\mathfrak{P}_{i_{1}}^{\left(  \mathbf{d}%
\right)  }\leq\cdots\leq\mathfrak{P}_{i_{n_{0}-1}}^{\left(  \mathbf{d}\right)
}<P\leq\mathfrak{P}_{i_{n_{0}}}^{\left(  \mathbf{d}\right)  }\leq\cdots
\leq\mathfrak{P}_{i_{K}}^{\left(  \mathbf{d}\right)  }=P_{H}\left(
\mathbf{d}\right)  ,
\]
then the optimal \textit{transformational }dynamic rationing policy is given
by%
\[
\mathbf{d}^{\ast}\left(  \text{Transfer}\right)  =\left(  0;\underset
{n_{0}-1\text{ zeros}}{\underbrace{0,0,\ldots,0}},\underset{K-n_{0}+1\text{
ones}}{\underbrace{1,1,\ldots,1}};1,1,\ldots,1\right)  .
\]
\end{The}

\textbf{Proof: }From the set $\Lambda_{1}^{G}$, it is easy to see that
$P>\overline{P}_{H}^{G}\left(  \mathbf{d;}1\rightarrow n_{0}-1\right)  $.
Hence we consider the transformational\textit{ }sub-policy%
\[
\widetilde{\mathbf{d}}_{a}\left(  \text{Transfer}\right)  =\left(  0;d_{i_{1}%
},d_{i_{2}},\ldots,d_{i_{n_{0}-1}},\ast,\ast,\ldots,\ast;1,1,\ldots,1\right)
.
\]
By observing the transformational\textit{ }sub-policy $\left(  d_{i_{1}%
},d_{i_{2}},\ldots,d_{i_{n_{0}-1}}\right)  $ related to $P>\overline{P}%
_{H}^{G}(\mathbf{d;}1\rightarrow n_{0}-1)$, it is easy to see from the proof
of Theorem \ref{The:left} that the optimal transformational dynamic rationing
sub-policy is given by%
\[
\widetilde{\mathbf{d}}_{a}^{\ast}\left(  \text{Transfer}\right)  =\left(
0;0,0,\ldots,0,\ast,\ast,\ldots,\ast;1,1,\ldots,1\right)  .
\]

Similarly, from $0\leq P\leq\overline{P}_{L}^{G}\left(  \mathbf{d;}%
n_{0}\rightarrow K\right)  $ in the set $\Lambda_{2}^{G}$, we discuss the
transformational sub-policy%
\[
\widetilde{\mathbf{d}}_{b}\left(  \text{Transfer}\right)  =\left(  0;\ast
,\ast,\ldots,\ast,d_{i_{n_{0}}},d_{i_{n_{0}+1}},\ldots,d_{i_{K}}%
;1,1,\ldots,1\right)  .
\]
By observing the transformational sub-policy $\left(  d_{n_{0}},d_{n_{0}%
+1},\ldots,d_{K}\right)  $ related to $0\leq P\leq\overline{P}_{L}^{G}\left(
\mathbf{d;}n_{0}\rightarrow K\right)  $, it is easy to see from the proof of
Theorem \ref{The:right} that the optimal transformational dynamic rationing
sub-policy is given by%
\[
\widetilde{\mathbf{d}}_{b}^{\ast}\left(  \text{Transfer}\right)  =\left(
0;\ast,\ast,\ldots,\ast,1,1,\ldots,1;1,1,\ldots,1\right)  .
\]
Therefore, by observing the total transformational sub-policy $(d_{i_{1}%
},d_{i_{2}},\ldots,d_{i_{n_{0}-1}},d_{i_{n_{0}}},d_{i_{n_{0}+1}}$,
$\ldots,d_{i_{K}})$ in the total set $\Lambda_{1}^{G}\cup\Lambda_{2}^{G}$, the
optimal transformational dynamic rationing policy is given by%
\begin{align*}
\mathbf{d}^{\ast}\left(  \text{Transfer}\right)   &  =\widetilde{\left(
\widetilde{\mathbf{d}}_{a}^{\ast}\left(  \text{Transfer}\right)  \right)
}_{b}^{\ast}\left(  \text{Transfer}\right)  =\widetilde{\left(  \widetilde
{\mathbf{d}}_{b}^{\ast}\left(  \text{Transfer}\right)  \right)  }_{a}^{\ast
}\left(  \text{Transfer}\right) \\
&  =\left(  0;\underset{n_{0}-1\text{ zeros}}{\underbrace{0,0,\ldots,0}%
},\underset{K-n_{0}+1\text{ ones}}{\underbrace{1,1,\ldots,1}};1,1,\ldots
,1\right)  .
\end{align*}
This completes the proof. \textbf{{\rule{0.08in}{0.08in}}}

\begin{Rem}
(1) For the general case, although the optimal dynamic rationing policy is not
of threshold type, we show that it must be of transformational threshold type.
Thus the optimal transformational dynamic policy of the stock-rationing queue
has a beautiful form as follows:%
\[
\mathbf{d}^{\ast}\left(  \text{Transfer}\right)  =\left(  0;\underset
{n_{0}-1\text{ zeros}}{\underbrace{0,0,\ldots,0}},\underset{K-n_{0}\text{
ones}}{\underbrace{1,1,\ldots,1};1,1,\ldots,1}\right)  .
\]

(2) We use an inverse transformation of $\mathbf{d}^{\ast}\left(
\text{Transfer}\right)  $ to be able to restore the original optimal dynamic
policy $\mathbf{d}^{\ast}$, since $\mathbf{d}^{\ast}\left(  \text{Transfer}%
\right)  $ is always obtained easily. To indicate such an inverse process, we
take a simple example:%
\[
\mathfrak{P}_{1}^{\left(  \mathbf{d}\right)  }\leq\mathfrak{P}_{3}^{\left(
\mathbf{d}\right)  }\leq\mathfrak{P}_{4}^{\left(  \mathbf{d}\right)  }%
\leq\mathfrak{P}_{7}^{\left(  \mathbf{d}\right)  }<P\leq\mathfrak{P}%
_{2}^{\left(  \mathbf{d}\right)  }\leq\mathfrak{P}_{5}^{\left(  \mathbf{d}%
\right)  }\leq\mathfrak{P}_{6}^{\left(  \mathbf{d}\right)  }\leq
\mathfrak{P}_{8}^{\left(  \mathbf{d}\right)  },
\]
it is easy to check that%
\[
\mathbf{d}^{\ast}=\left(  0;0,1,0,0,1,1,0,1;1,1,1,1\right)  .
\]
\end{Rem}

\begin{Rem}
The transformational version $\mathbf{d}^{\ast}\left(  \text{Transfer}\right)
$ of the optimal dynamic rationing policy $\mathbf{d}^{\ast}$ plays a key role
in the applications of the sensitivity-based optimization to the study of
stock-rationing queues. On the other hand, it is worthwhile to note that the
$RG$-factorization of block-structured Markov processes can be extended and
generalized to a more general optimal transformational version $\mathbf{d}%
^{\ast}\left(  \text{Transfer}\right)  $ in the study of stock-rationing
block-structured queues. See Li \cite{Li:2010} and \cite{Ma:2018} for more details.
\end{Rem}

\begin{Rem}
The bang-bang control is an effective method to roughly describe the optimal
dynamic policy, e.g., see Xia et al. \cite{Xia:2017, Xia:2018} and Ma et al.
\cite{Ma:2018}. However, our optimal transformational dynamic policy
$\mathbf{d}^{\ast}\left(  \text{Transfer}\right)  $ provides a more detailed
result, and also can restore the original optimal dynamic policy
$\mathbf{d}^{\ast}$ by means of an inverse transformation: $\left(
i_{1},i_{2},\ldots,i_{K-1},i_{K}\right)  \Rightarrow\left(  1,2,\ldots
,K-1,K\right)  $. Therefore, our optimal transformational dynamic rationing
policy is superior to the bang-bang control.
\end{Rem}

The following theorem provides a useful summarization for Theorems
\ref{The:left} to \ref{The:Trans}, this shows that we provide a complete
algebraic solution to the optimal dynamic policy of the stock-rationing queue.
Therefore, Problems \textbf{(P-a)} to \textbf{(P-c)} proposed in Introduction
are completely solved by means of our algebraic method.

\begin{The}
\label{The:Level}For the stock-rationing queue with two demand classes, there
must exist an optimal transformational dynamic rationing policy%
\[
\mathbf{d}^{\ast}\left(  \text{Transfer}\right)  =\left(  0;\underset
{n_{0}-1\text{ zeros}}{\underbrace{0,0,\ldots,0}},\underset{K-n_{0}\text{
ones}}{\underbrace{1,1,\ldots,1};1,1,\ldots,1}\right)  .
\]
Based on this finding, we can achieve the following two useful results:

(a) The optimal dynamic rationing policy $\mathbf{d}^{\ast}$ is of critical
rationing level (i.e., threshold type) under each of the three conditions: (i)
$P\geq P_{H}\left(  \mathbf{d}\right)  $ for any given policy $\mathbf{d}$;
(i\negthinspace i) $P_{L}\left(  \mathbf{d}\right)  >0$ and $0\leq P\leq
P_{L}\left(  \mathbf{d}\right)  $ for any given policy $\mathbf{d}$; and
(i\negthinspace i\negthinspace i) $P_{L}\left(  \mathbf{d}\right)
<P<P_{H}\left(  \mathbf{d}\right)  $ with the subscript vector $\left(
1,2,\ldots,K-1,K\right)  $ depending on Policy $\mathbf{d}$.

(b) The optimal dynamic rationing policy is not of critical rationing level
(i.e., threshold type) if $P_{L}\left(  \mathbf{d}\right)  <P<P_{H}\left(
\mathbf{d}\right)  $\ with the subscript vector $\left(  i_{1},i_{2}%
,\ldots,i_{K-1},i_{K}\right)  \neq\left(  1,2,\ldots,K-1,K\right)  $ depending
on Policy $\mathbf{d}$.
\end{The}

\subsection{A global optimal analysis}

In this subsection, for a fixed penalty cost $P$, we discuss how to find a
global optimal policy of the stock-rationing queue with two demand classes by
means of Theorem \ref{The:Level}. Note that if $\mathbf{d}^{\ast}$ is a global
optimal policy of this system, then $\mathbf{d}^{\ast}\succeq\mathbf{c}$ for
any $\mathbf{c}\in\mathcal{D}$. Also, we provide a simple effective method to
be able to find the global optimal policy from the policy set $\mathcal{D}$.

In the policy set $\mathcal{D}$, we define two key policies:%
\[
\mathbf{d}_{1}=\left(  0;0,0,\ldots,0;1,1,\ldots,1\right)
\]
and%
\[
\mathbf{d}_{2}=\left(  0;1,1,\ldots,1;1,1,\ldots,1\right)  .
\]
Note that there are $2^{k}$ different policies in the set $\mathcal{D}$, we
write%
\[
\mathcal{D}=\left\{  \mathbf{d}_{1},\mathbf{d}_{2};\mathbf{c}_{3}%
,\mathbf{c}_{4},\ldots,\mathbf{c}_{2^{k}-1},\mathbf{c}_{2^{k}}\right\}  .
\]

The following theorem describes a useful characteristics of the two key
policies $\mathbf{d}_{1}$ and $\mathbf{d}_{2}$ by means of the class property
of the policies in the set $\mathcal{D}$, given in Theorem \ref{The:NP}. This
characteristics makes us to be able to find the global optimal policy of the
stock-rationing queue.

\begin{The}
\label{The:bound}(1) If a fixed penalty cost $P\geq P_{H}\left(
\mathbf{d}\right)  $ for any given policy $\mathbf{d}$, then $P\geq
P_{H}\left(  \mathbf{d}_{1}\right)  $.

(2) If a fixed penalty cost $P_{L}\left(  \mathbf{d}\right)  >0$ and $0\leq
P\leq P_{L}\left(  \mathbf{d}\right)  $ for any given policy $\mathbf{d}$,
then $P_{L}\left(  \mathbf{d}_{2}\right)  >0$ and $0\leq P\leq P_{L}\left(
\mathbf{d}_{2}\right)  $.
\end{The}

\textbf{Proof: }We only prove (1), while (2) can be proved similarly.

We assume the penalty cost: $P<P_{H}\left(  \mathbf{d}_{1}\right)  $ for
Policy $\mathbf{d}_{1}=\left(  0;0,0,\ldots,0;1,1,\ldots,1\right)  $. Then
there exists the minimal positive integer $n_{0}\in\left\{  1,2,\ldots
,K-1,K\right\}  $ such that
\[
0<P\leq\mathfrak{P}_{i_{n_{0}}}^{\left(  \mathbf{d}_{1}\right)  }\leq
\cdots\leq\mathfrak{P}_{i_{K}}^{\left(  \mathbf{d}_{1}\right)  }=P_{H}\left(
\mathbf{d}_{1}\right)  ,
\]
and also there exists at least a positive integer $m_{0}\in\left\{
n_{0}+1,n_{0}+2,\ldots,K-1,K\right\}  $ such that%
\begin{equation}
\mathfrak{P}_{i_{m_{0}-1}}^{\left(  \mathbf{d}_{1}\right)  }<\mathfrak{P}%
_{i_{m_{0}}}^{\left(  \mathbf{d}_{1}\right)  }. \label{Equa-24}%
\end{equation}

Let%
\[
\overline{P}_{L}^{G}\left(  \mathbf{d}_{1},n_{0}\rightarrow K\right)
=\min\left\{  \mathfrak{P}_{i_{n_{0}}}^{\left(  \mathbf{d}_{1}\right)
},\mathfrak{P}_{i_{n_{0}+1}}^{\left(  \mathbf{d}_{1}\right)  },\ldots
,\mathfrak{P}_{i_{K-1}}^{\left(  \mathbf{d}_{1}\right)  },\mathfrak{P}_{i_{K}%
}^{\left(  \mathbf{d}_{1}\right)  }\right\}  =\mathfrak{P}_{i_{n_{0}}%
}^{\left(  \mathbf{d}_{1}\right)  }>0.
\]
Then from $0\leq P\leq\overline{P}_{L}^{G}\left(  \mathbf{d}_{1}%
,n_{0}\rightarrow K\right)  $, we discuss the transformational sub-policy%
\[
\widetilde{\left(  \mathbf{d}_{1}\right)  }_{b}\left(  \text{Transfer}\right)
=\left(  0;\ast,\ast,\ldots,\ast,d_{i_{n_{0}}},d_{i_{n_{0}+1}},\ldots
,d_{i_{K}};1,1,\ldots,1\right)  .
\]
By observing the transformational sub-policy $\left(  d_{n_{0}},d_{n_{0}%
+1},\ldots,d_{K}\right)  $ related to $0\leq P\leq\overline{P}_{L}^{G}\left(
\mathbf{d}_{1},n_{0}\rightarrow K\right)  $, it is easy to see from the proof
of Theorem \ref{The:right} that the optimal transformational dynamic rationing
sub-policy is given by%
\[
\widetilde{\left(  \mathbf{d}_{1}\right)  }_{b}^{\ast}\left(  \text{Transfer}%
\right)  =\left(  0;\ast,\ast,\ldots,\ast,1,1,\ldots,1;1,1,\ldots,1\right)  .
\]
This gives%
\begin{equation}
\widetilde{\left(  \mathbf{d}_{1}\right)  }_{b}^{\ast}\left(  \text{Transfer}%
\right)  \succ\mathbf{d}_{1}=\mathbf{d}^{\ast} \label{Equa-21}%
\end{equation}
by using (\ref{Equa-24}), where $\mathbf{d}^{\ast}$ is given in Theorem
\ref{The:left}.

Since for a fixed penalty cost $P\geq P_{H}\left(  \mathbf{d}\right)  $ for
Policy $\mathbf{d}$, it follows from Theorem \ref{The:left} that the optimal
dynamic rationing policy of the stock-rationing queue is given by
\[
\mathbf{d}^{\ast}=\left(  0;0,0,\ldots,0;1,1,\ldots,1\right)  .
\]
By using (\ref{Equa-24}), we obtain%
\begin{equation}
\widetilde{\left(  \mathbf{d}_{1}\right)  }_{b}^{\ast}\left(  \text{Transfer}%
\right)  \prec\mathbf{d}^{\ast}. \label{Equa-22}%
\end{equation}

This makes a contradiction between (\ref{Equa-21}) and (\ref{Equa-22}), thus
our assumption on the penalty cost: $P<P_{H}\left(  \mathbf{d}_{1}\right)  $
should not be correct. This completes the proof.
\textbf{{\rule{0.08in}{0.08in}}}

Theorem \ref{The:bound} shows that to find the optimal dynamic rationing
policy of the stock-rationing queue, our first step is to check whether there
exists (a) the penalty cost $P\geq P_{H}\left(  \mathbf{d}_{1}\right)  $, or
(b) the fixed penalty cost $P_{L}\left(  \mathbf{d}_{2}\right)  >0$ and $0\leq
P\leq P_{L}\left(  \mathbf{d}_{2}\right)  $. Thus, the two special policies
$\mathbf{d}_{1}$ and $\mathbf{d}_{2}$ are chosen as the first observation of
our algebraic method on the the optimal dynamic rationing policy.

The following theorem provides the global optimal solution to the optimal
dynamic rationing policy of the stock-rationing queue.

\begin{The}
In the stock-rationing queue with two demand classes, we have

(1) If a fixed penalty cost $P\geq P_{H}\left(  \mathbf{d}_{1}\right)  $, then
$\mathbf{d}^{\ast}=\mathbf{d}_{1}\succeq\mathbf{c}$ for any $\mathbf{c}%
\in\mathcal{D}$.

(2) If a fixed penalty cost $P_{L}\left(  \mathbf{d}_{2}\right)  >0$ or $0\leq
P\leq P_{L}\left(  \mathbf{d}_{2}\right)  $, then $\mathbf{d}^{\ast
}=\mathbf{d}_{2}\succeq\mathbf{c}$ for any $\mathbf{c}\in\mathcal{D}$.

(3) If a fixed penalty cost $P$ satisfies $P<P_{H}\left(  \mathbf{d}%
_{1}\right)  $ and $P>P_{L}\left(  \mathbf{d}_{2}\right)  $, then%
\[
\mathbf{d}^{\ast}=\max\left\{  \widetilde{\left(  \mathbf{d}_{1}\right)  }%
_{b}^{\ast}\left(  \text{Transfer}\right)  ,\widetilde{\left(  \mathbf{d}%
_{2}\right)  }_{a}^{\ast}\left(  \text{Transfer}\right)  ,\left(
\mathbf{c}_{k}\right)  ^{\ast}\left(  \text{Transfer}\right)  \text{ for
}k=3,4,\ldots,K\right\}
\]
and $\mathbf{d}^{\ast}\succeq\mathbf{c}$ for any $\mathbf{c}\in\mathcal{D}$.
\end{The}

\textbf{Proof: }We only prove (3), while (1) and (2) are provided in those of
Theorem \ref{The:bound}.

If $P<P_{H}\left(  \mathbf{d}_{1}\right)  $ or $P>P_{L}\left(  \mathbf{d}%
_{2}\right)  $, then both $\mathbf{d}_{1}$ and $\mathbf{d}_{2}$ are not the
optimal dynamic rationing policy of the system. In this case, by using Theorem
\ref{The:Level}, we indicate that the optimal dynamic rationing policy must be
of transformational threshold type. Thus we have%
\[
\mathbf{d}^{\ast}=\max\left\{  \widetilde{\left(  \mathbf{d}_{1}\right)  }%
_{b}^{\ast}\left(  \text{Transfer}\right)  ,\widetilde{\left(  \mathbf{d}%
_{2}\right)  }_{a}^{\ast}\left(  \text{Transfer}\right)  ,\left(
\mathbf{c}_{k}\right)  ^{\ast}\left(  \text{Transfer}\right)  \text{ for
}k=3,4,\ldots,K\right\}  ,
\]
which is of transformational threshold type, since $K$ is a finite positive
integer. It is clear that $\mathbf{d}^{\ast}\succeq\mathbf{c}$ for any
$\mathbf{c}\in\mathcal{D}$. This completes the proof.
\textbf{{\rule{0.08in}{0.08in}}}

\section{The Static Rationing Policies}

In this section, we analyze the static (i.e., threshold type) rationing
policies of the stock-rationing queue with two demand classes, and discuss the
optimality of the static rationing policies. Furthermore, we provide a
necessary condition under which a static rationing policy is optimal. Based on
this, we can intuitively understand some differences between the optimal
static and dynamic rationing policies.

To study static rationing policy, we define a static policy subset of the
policy set $\mathcal{D}$ as follows. For $\theta=1,2,\ldots,K,K+1$, we write
$\mathbf{d}_{\triangle,\theta}$ as a static rationing policy $\mathbf{d}$ with
$d_{i}=0$ if $1\leq i\leq\theta-1$ and $d_{i}=1$ if $\theta\leq i\leq K$.
Clearly, if $\theta=1$, then%
\[
\mathbf{d}_{\triangle,1}=\left(  0;1,1,\ldots,1;1,1,\ldots,1\right)  ;
\]
if $\theta=K$, then%
\[
\mathbf{d}_{\triangle,K}=\left(  0;0,0,\ldots,0,1;1,1,\ldots,1\right)  ;
\]
and if $\theta=K+1$, then%
\[
\mathbf{d}_{\triangle,K+1}=\left(  0;0,0,\ldots,0;1,1,\ldots,1\right)  .
\]
Let%
\[
\mathcal{D}^{\Delta}=\left\{  \mathbf{d}_{\triangle,\theta}:\theta
=1,2,\ldots,K,K+1\right\}  .
\]
Then%
\[
\mathcal{D}^{\Delta}=\left\{  \left(  0;\underset{\theta-1\text{ zeros}%
}{\underbrace{0,0,\ldots,0}},1,1,\ldots,1;1,1,\ldots,1\right)  :\theta
=1,2,\ldots,K,K+1\right\}  .
\]
It is easy to see that the static rationing policy set $\mathcal{D}^{\Delta
}\subset\mathcal{D}$.

For a static rationing policy $\mathbf{d}_{\triangle,\theta}=\left(
0;\underset{\theta-1\text{ zeros}}{\underbrace{0,0,\ldots,0}},1,1,\ldots
,1;1,1,\ldots,1\right)  $ with $\theta=1,2,\ldots,K,K+1$, it follows from
(\ref{5}) that%
\[%
\begin{array}
[c]{l}%
\text{ \ \ \ \ }\xi_{0}=1,\text{ \ \ \ \ \ \ \ \ \ \ \ \ \ \ \ \ \ \ }i=0;\\
\xi_{i}^{(\mathbf{d}_{\triangle,\theta})}=\left\{
\begin{array}
[c]{ll}%
\alpha^{i}, & i=1,2,\ldots,\theta-1;\\
\left(  \frac{\alpha}{\beta}\right)  ^{\theta-1}\beta^{i}, & i=\theta
,\theta+1,\ldots,N.
\end{array}
\right.
\end{array}
\]
and%
\begin{align*}
h^{(\mathbf{d}_{\triangle,\theta})}  &  =1+\sum\limits_{i=1}^{N}\xi
_{i}^{(\mathbf{d}_{\triangle,\theta})}\\
&  =1+\frac{\alpha\left(  1-\alpha^{\theta-1}\right)  }{1-\alpha}+\left(
\frac{\alpha}{\beta}\right)  ^{\theta-1}\frac{\beta^{\theta}\left(
1-\beta^{N-\theta+1}\right)  }{1-\beta}.
\end{align*}
It follows from (\ref{37}) that%
\[
\pi^{(\mathbf{d}_{\triangle,\theta})}\left(  i\right)  =\left\{
\begin{array}
[c]{ll}%
\frac{1}{h^{(\mathbf{d}_{\triangle,\theta})}}, & i=0;\\
\frac{1}{h^{(\mathbf{d}_{\triangle,\theta})}}\alpha^{i}, & i=1,2,\ldots
,\theta-1;\\
\frac{1}{h^{(\mathbf{d}_{\triangle,\theta})}}\left(  \frac{\alpha}{\beta
}\right)  ^{\theta-1}\beta^{i}, & i=\theta,\theta+1,\ldots,N.
\end{array}
\right.
\]

On the other hand, it follows from (\ref{8}) to (\ref{11}) that for $i=0$
\[
f\left(  0\right)  =-C_{2,1}\mu_{1}-C_{2,2}\mu_{2}-C_{3}\lambda;
\]
for $i=1,2,\ldots,\theta-1$,
\[
f^{\left(  \mathbf{d}_{\triangle,\theta}\right)  }\left(  i\right)  =R\mu
_{1}-C_{1}i-C_{2,2}\mu_{2}-C_{3}\lambda;
\]
for $i=\theta,\theta+1,\ldots,K,$
\[
f^{\left(  \mathbf{d}_{\triangle,\theta}\right)  }\left(  i\right)  =R\left(
\mu_{1}+\mu_{2}\right)  -C_{1}i-C_{3}\lambda-P\mu_{2};
\]
and for $i=K+1,K+2,\ldots,N,$%
\[
f\left(  i\right)  =R\left(  \mu_{1}+\mu_{2}\right)  -C_{1}i-C_{3}%
\lambda1_{\left\{  i<N\right\}  }-C_{4}\lambda1_{\left\{  i=N\right\}  }.
\]
Note that%
\begin{align*}
\eta^{\mathbf{d}_{\triangle,\theta}}  &  =\pi^{\left(  \mathbf{d}%
_{\triangle,\theta}\right)  }\left(  0\right)  f\left(  0\right)  +\sum
_{i=1}^{\theta-1}\pi^{\left(  \mathbf{d}_{\triangle,\theta}\right)  }\left(
i\right)  f^{\left(  \mathbf{d}_{\triangle,\theta}\right)  }\left(  i\right)
\\
&  \text{ \ \ \ }+\sum_{i=\theta}^{K}\pi^{\left(  \mathbf{d}_{\triangle
,\theta}\right)  }\left(  i\right)  f^{\left(  \mathbf{d}_{\triangle,\theta
}\right)  }\left(  i\right)  +\sum_{i=K+1}^{N}\pi^{\left(  \mathbf{d}%
_{\triangle,\theta}\right)  }\left(  i\right)  f\left(  i\right)  ,
\end{align*}
we obtain an explicit expression for the long-run average profit of the
stock-rationing queue under the static rationing policy $\mathbf{d}%
_{\triangle,\theta}$ as follows:%
\begin{align*}
\eta^{\mathbf{d}_{\triangle,\theta}}  &  =\frac{1}{h^{(\mathbf{d}%
_{\triangle,\theta})}}\left\{  -\left(  C_{2,1}\mu_{1}+C_{2,2}\mu_{2}%
+C_{3}\lambda\right)  +\sum_{i=1}^{\theta-1}\alpha^{i}\left(  R\mu_{1}%
-C_{1}i-C_{2,2}\mu_{2}-C_{3}\lambda\right)  \right. \\
&  \text{ \ \ \ }+\sum_{i=\theta}^{K}\left(  \frac{\alpha}{\beta}\right)
^{\theta-1}\beta^{i}\left[  R\left(  \mu_{1}+\mu_{2}\right)  -C_{1}%
i-C_{3}\lambda-P\mu_{2}\right] \\
&  \text{ \ \ \ }\left.  +\sum_{i=K+1}^{N}\left(  \frac{\alpha}{\beta}\right)
^{\theta-1}\beta^{i}\left[  R\left(  \mu_{1}+\mu_{2}\right)  -C_{1}%
i-C_{3}\lambda1_{\left\{  i<N\right\}  }-C_{4}\lambda1_{\left\{  i=N\right\}
}\right]  \right\} \\
&  =\frac{1}{h^{(\mathbf{d}_{\triangle,\theta})}}\left\{  -\gamma_{1}%
+\gamma_{2}\frac{\alpha\left(  1-\alpha^{\theta-1}\right)  }{1-\alpha}%
-C_{1}\left[  \frac{\alpha\left(  1-\alpha^{\theta-1}\right)  }{\left(
1-\alpha\right)  ^{2}}-\frac{\left(  \theta-1\right)  \alpha^{\theta}%
}{1-\alpha}\right]  \right. \\
&  \text{ \ \ \ }-\left(  \frac{\alpha}{\beta}\right)  ^{\theta-1}C_{1}\left[
\frac{\left(  \theta-1\right)  \beta^{\theta}-N\beta^{N+1}}{1-\beta
}+\frac{\beta^{\theta}\left(  1-\beta^{N-\theta+1}\right)  }{\left(
1-\beta\right)  ^{2}}\right] \\
&  \left.  \text{ \ \ }+\left(  \frac{\alpha}{\beta}\right)  ^{\theta-1}%
\gamma_{4}\frac{\beta^{\theta}\left(  1-\beta^{K-\theta+1}\right)  }{1-\beta
}+\left(  \frac{\alpha}{\beta}\right)  ^{\theta-1}\gamma_{3}\frac{\beta
^{K+1}\left(  1-\beta^{N-K}\right)  }{1-\beta}\text{\ \ \ \ }\right\}  .
\end{align*}

Let%
\[
\mathbf{d}_{\triangle,\theta}^{\ast}=\underset{\mathbf{d}_{\triangle,\theta
}\in\mathcal{D}^{\Delta}}{\arg\max}\left\{  \eta^{\mathbf{d}_{\triangle
,\theta}}\right\}
\]
and%
\[
\mathbf{d}_{\triangle,\theta^{\ast}}=\underset{1\leq\theta\leq K+1}{\arg\max
}\left\{  \eta^{\mathbf{d}_{\triangle,\theta}}\right\}  .
\]
Then $\mathbf{d}_{\triangle,\theta}^{\ast}=\mathbf{d}_{\triangle,\theta^{\ast
}}$. Hence we call $\mathbf{d}_{\triangle,\theta}^{\ast}$ (or $\mathbf{d}%
_{\triangle,\theta^{\ast}}$) the optimal static rationing policy in the static
rationing policy set $\mathcal{D}^{\Delta}$. Since $\mathcal{D}^{\Delta
}\subset\mathcal{D}$, the partially ordered set $\mathcal{D}$ shows that
$\mathcal{D}^{\Delta}$ is also a partially ordered set. Based on this, it is
easy to see from the two partially ordered sets $\mathcal{D}$ and
$\mathcal{D}^{\Delta}$ that%
\[
\eta^{\mathbf{d}_{\triangle,\theta}^{\ast}}\leq\eta^{\mathbf{d}^{\ast}}\text{,
\ \ or }\mathbf{d}_{\triangle,\theta}^{\ast}\preceq\mathbf{d}^{\ast},
\]
where $\mathbf{d}^{\ast}$ is the optimal dynamic rationing policy in the set
$\mathcal{D}$.

If $\eta^{\mathbf{d}_{\triangle,\theta}^{\ast}}=\eta^{\mathbf{d}^{\ast}}$,
then the optimal static rationing policy $\mathbf{d}_{\triangle,\theta}^{\ast
}$ is also optimal in the policy set $\mathcal{D}$, thus the optimal dynamic
rationing policy is of threshold type. If $\eta^{\mathbf{d}_{\triangle,\theta
}^{\ast}}<\eta^{\mathbf{d}^{\ast}}$, then the optimal static rationing policy
$\mathbf{d}_{\triangle,\theta}^{\ast}$ is not optimal in the static rationing
policy subset $\mathcal{D}^{\Delta}$, i.e., it is also suboptimal in the
policy set $\mathcal{D}$, thus the optimal dynamic rationing policy is not of
threshold type.

Now, we set up some conditions under which the optimal static rationing policy
$\mathbf{d}_{\triangle,\theta}^{\ast}$ is suboptimal in the dynamic rationing
policy set $\mathcal{D}$.

In the static rationing policy subset $\mathcal{D}^{\Delta}$, it is easy to
see that there must exist a minimal positive integer $\theta^{\ast}\in\left\{
1,2,\ldots,K,K+1\right\}  $ such that%
\[
\mathbf{d}_{\triangle,\theta}^{\ast}=\mathbf{d}_{\triangle,\theta^{\ast}%
}=\left(  0;\underset{\mathbf{\theta}^{\ast}-1\text{ zeros}}{\underbrace
{0,0,\ldots,0},}1,1,\ldots,1;1,1,\ldots,1\right)  .
\]

By using the optimal static rationing policy $\mathbf{d}_{\triangle,\theta
}^{\ast}$ (or $\mathbf{d}_{\triangle,\mathbf{\theta}^{\ast}}$), the following
theorem determines the sign of the function $\ G^{\left(  \mathbf{d}%
_{\triangle,\theta}\right)  }\left(  \theta\right)  +b$ in the three different
points: $\theta=\theta^{\ast}-1,\theta^{\ast},\theta^{\ast}+1$. This may be
useful for us to understand how to use Proposition 1 to give the optimal
long-run average profit of this system.

\begin{The}
In the stock-rationing queue, the static rationing policies $\mathbf{d}%
_{\triangle,\mathbf{\theta}^{\ast}-1}$, $\mathbf{d}_{\triangle,\mathbf{\theta
}^{\ast}}$ and $\mathbf{d}_{\triangle,\mathbf{\theta}^{\ast}+1}$ satisfy the
following conditions:%
\[
G^{\left(  \mathbf{d}_{\triangle,\mathbf{\theta}^{\ast}-1}\right)  }\left(
\theta^{\ast}-1\right)  +b\leq0\mathbf{,}\text{ \ }G^{\left(  \mathbf{d}%
_{\triangle,\mathbf{\theta}^{\ast}}\right)  }\left(  \theta^{\ast}-1\right)
+b\leq0,
\]
and%
\[
G^{\left(  \mathbf{d}_{\triangle,\mathbf{\theta}^{\ast}}\right)  }\left(
\theta^{\ast}\right)  +b\geq0\mathbf{,}\text{ \ }G^{\left(  \mathbf{d}%
_{\triangle,\mathbf{\theta}^{\ast}+1}\right)  }\left(  \theta^{\ast}\right)
+b\geq0.
\]
\end{The}

\textbf{Proof:} We consider three static rationing policies with an
interrelated structure as follows:%
\begin{align*}
\mathbf{d}_{\triangle,\mathbf{\theta}^{\ast}-1}  &  =\left(  0;\underset
{\mathbf{\theta}^{\ast}-2\text{ zeros}}{\underbrace{0,0,\ldots,0}%
,}1,1,1,1,\ldots,1;1,1,\ldots,1\right)  ,\\
\mathbf{d}_{\triangle,\mathbf{\theta}^{\ast}}  &  =\left(  0;\underset
{\mathbf{\theta}^{\ast}-1\text{ zeros}}{\underbrace{0,0,\ldots,0,0}%
,}1,1,1,\ldots,1;1,1,\ldots,1\right)  ,\\
\mathbf{d}_{\triangle,\mathbf{\theta}^{\ast}+1}  &  =\left(  0;\underset
{\mathbf{\theta}^{\ast}\text{ zeros}}{\underbrace{0,0,\ldots,0,0,0},}%
1,\ldots,1;1,1,\ldots,1\right)  .
\end{align*}
Note that $\mathbf{d}_{\triangle,\mathbf{\theta}^{\ast}}$ is the optimal
static rationing policy, and $\mathbf{d}_{\triangle,\theta}^{\ast}%
=\mathbf{d}_{\triangle,\mathbf{\theta}^{\ast}}$. It is clear that
$\mathbf{d}_{\triangle,\mathbf{\theta}^{\ast}}\succeq\mathbf{d}_{\triangle
,\mathbf{\theta}^{\ast}-1}$ and $\mathbf{d}_{\triangle,\mathbf{\theta}^{\ast}%
}\succeq\mathbf{d}_{\triangle,\mathbf{\theta}^{\ast}+1}.$ Thus it follows from
(\ref{Equa-2}) that%
\[
\eta^{\mathbf{d}_{\triangle,\theta^{\ast}+1}}-\eta^{\mathbf{d}_{\triangle
,\mathbf{\theta}^{\ast}}}=-\mu_{2}\pi^{\left(  \mathbf{d}_{\triangle
,\mathbf{\theta}^{\ast}+1}\right)  }\left(  \theta^{\ast}\right)  \left[
G^{\left(  \mathbf{d}_{\triangle,\mathbf{\theta}^{\ast}}\right)  }\left(
\theta^{\ast}\right)  +b\right]  ,
\]
which, together with $\eta^{\mathbf{d}_{\triangle,\theta^{\ast}+1}}%
-\eta^{\mathbf{d}_{\triangle,\mathbf{\theta}^{\ast}}}\leq0$, leads to%
\[
G^{\left(  \mathbf{d}_{\triangle,\mathbf{\theta}^{\ast}}\right)  }\left(
\theta^{\ast}\right)  +b\geq0\mathbf{.}%
\]
On the other hand, it follows from (\ref{Equa-2}) that
\[
\eta^{\mathbf{d}_{\triangle,\mathbf{\theta}^{\ast}}}-\eta^{\mathbf{d}%
_{\triangle,\theta^{\ast}+1}}=\mu_{2}\pi^{\left(  \mathbf{d}_{\triangle
,\mathbf{\theta}^{\ast}}\right)  }\left(  \theta^{\ast}\right)  \left[
G^{\left(  \mathbf{d}_{\triangle,\theta^{\ast}+1}\right)  }\left(
\theta^{\ast}\right)  +b\right]  ,
\]
this gives%
\[
G^{\left(  \mathbf{d}_{\triangle,\theta^{\ast}+1}\right)  }\left(
\theta^{\ast}\right)  +b\geq0
\]

Similarly, by using $\eta^{\mathbf{d}_{\triangle,\mathbf{\theta}^{\ast}}}%
\geq\eta^{\mathbf{d}_{\triangle,\theta^{\ast}-1}}$ and
\[
\eta^{\mathbf{d}_{\triangle,\mathbf{\theta}^{\ast}}}-\eta^{\mathbf{d}%
_{\triangle,\theta^{\ast}-1}}=-\mu_{2}\pi^{\left(  \mathbf{d}_{\triangle
,\mathbf{\theta}^{\ast}}\right)  }\left(  \theta^{\ast}-1\right)  \left[
G^{\left(  \mathbf{d}_{\triangle,\theta^{\ast}-1}\right)  }\left(
\theta^{\ast}-1\right)  +b\right]  ,
\]
this gives%
\[
G^{\left(  \mathbf{d}_{\triangle,\theta^{\ast}-1}\right)  }\left(
\theta^{\ast}-1\right)  +b\leq0\mathbf{;}%
\]
and%
\[
\eta^{\mathbf{d}_{\triangle,\theta^{\ast}-1}}-\eta^{\mathbf{d}_{\triangle
,\mathbf{\theta}^{\ast}}}=\mu_{2}\pi^{\left(  \mathbf{d}_{\triangle
,\theta^{\ast}-1}\right)  }\left(  \theta^{\ast}-1\right)  \left[  G^{\left(
\mathbf{d}_{\triangle,\mathbf{\theta}^{\ast}}\right)  }\left(  \theta^{\ast
}-1\right)  +b\right]  ,
\]
we obtain%
\[
G^{\left(  \mathbf{d}_{\triangle,\mathbf{\theta}^{\ast}}\right)  }\left(
\theta^{\ast}-1\right)  +b\leq0.
\]
This completes the proof. \textbf{{\rule{0.08in}{0.08in}}}

\section{Numerical Experiments}

In this section, by observing several different penalty costs, we conduct
numerical experiments to demonstrate our theoretical results and to gain
insights on the optimal dynamic and static rationing policies in the
stock-rationing queue.

In Examples 1 to 4, we take some common parameters in the stock-rationing
queue with two demand classes as follows:%
\[
C_{1}=1,C_{2,1}=4,C_{2,2}=1,C_{3}=5,C_{4}=1,R=15,N=100.
\]

In Examples 1 and 2, we analyze some difference between the optimal static and
dynamic rationing policies, and use the optimal static rationing policy to
show whether or not the the optimal dynamic rationing policy is of threshold type.

\textbf{Example 1}. We give some useful comparisons of the optimal long-run
average profit between two different penalty costs, and further verify how the
optimality depends on the penalty cost in Theorems \ref{The:left} and
\ref{The:right} for the the optimal dynamic rationing policy. To this end, we
further take the system parameters as $\lambda=3$, $\mu_{1}=4$, $\mu_{2}=2$,
$K=15$ and $1\leq i\leq15$.

\textit{Case one: A higher penalty cost}

We take a higher penalty cost $P=10$. If $d_{i}^{\ast}=0$ for $1\leq i\leq15$
such that a possible optimal dynamic rationing policy $\mathbf{d}^{\ast
}=\left(  0;0,0,\ldots,0;1,1,\ldots,1\right)  $, then we obtain $\eta
^{\mathbf{d}^{\ast}}=22.3$. On the other hand, if $d_{i}^{\prime\ast}=1$ for
$1\leq i\leq15$ such that another possible optimal dynamic rationing policy
$\mathbf{d}^{\prime\ast}=\left(  0;1,1,\ldots,1;1,1,\ldots,1\right)  $, then
we get $\eta^{\mathbf{d}^{\prime\ast}}=13$. By comparing $\eta^{\mathbf{d}%
^{\ast}}=22.3$ with $\eta^{\mathbf{d}^{\prime\ast}}=13$, it is easy to see
that the possible optimal dynamic rationing policy should be $\mathbf{d}%
^{\ast}=\left(  0;0,0,\ldots,0;1,1,\ldots,1\right)  $.

\textit{Case two: A lower penalty cost}

We choose a lower penalty cost $P=0.1$. If $d_{i}^{\ast}=1$ for $1\leq
i\leq15$ such that a possible optimal dynamic rationing policy $\mathbf{d}%
^{\ast}=\left(  0;1,1,\ldots,1;1,1,\ldots,1\right)  $, then we obtain
$\eta^{\mathbf{d}^{\ast}}=22.9$. On the other hand, if $d_{i}^{\prime\ast}=0$
for $1\leq i\leq15$ such that another possible optimal dynamic rationing
policy $\mathbf{d}^{\prime\ast}=\left(  0;0,0,\ldots,0;1,1,\ldots,1\right)  $,
then we get $\eta^{\mathbf{d}^{\prime\ast}}=22.3$. Obviously, the possible
optimal dynamic rationing policy should be $\mathbf{d}^{\ast}=\left(
0;1,1,\ldots,1;1,1,\ldots,1\right)  $.

\textbf{Example 2}. We use the numerical example to demonstrate whether or not
the optimal static rationing policy is suboptimal in the policy set
$\mathcal{D}$. If yes, then we show that the optimal dynamic rationing policy
is not of threshold type. To this end, we take some system parameters:
$\lambda=3$, $\mu_{1}=4$, $\mu_{2}=2$, $K=15$, $1\leq\theta\leq15$. These
parameters are the same as those in Example 1.

In what follows our observation is to focus on the higher penalty cost $P=10$
and the\ lower penalty cost $P=0.1$, respectively.

\textit{Case one: A higher penalty cost}

We observe how the optimal long-run average profit $\eta^{\mathbf{d}^{\ast}}$
depends on the threshold from $\theta=1$ to $\theta=15$. From Figure 3, it is
seen that the optimal threshold is $\theta^{\ast}=9$ and $\eta^{\mathbf{d}%
_{\Delta,\theta^{\ast}}}=21.4$. However, from Case one of Example 1,
$\eta^{\mathbf{d}^{\ast}}=22.3$. Thus we obtain that $\eta^{\mathbf{d}%
_{\Delta,\theta^{\ast}}}=21.4<\eta^{\mathbf{d}^{\ast}}=22.3$. This shows that
the optimal static rationing policy is suboptimal in the policy set
$\mathcal{D}$, and the optimal dynamic rationing policy is not of threshold
type. Thus, $\mathbf{d}^{\ast}=\left(  0;0,0,\ldots,0;1,1,\ldots,1\right)  $,
given in Example 1, is not the optimal dynamic rationing policy yet.
\begin{figure}[ptbh]
\centering                          \includegraphics[width=8.5cm]{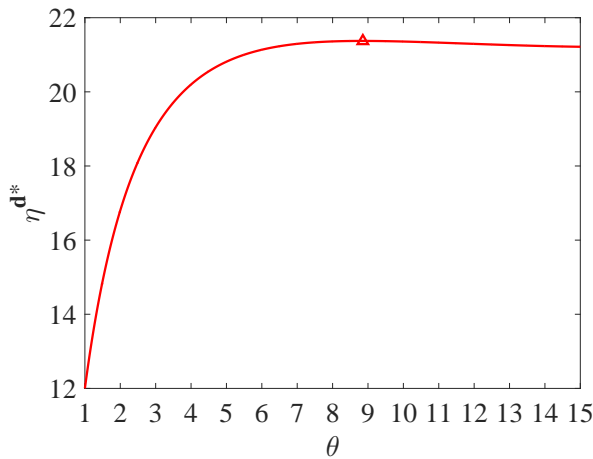}
\caption{The optimal long-run average profit $\eta^{\mathbf{d}^{\ast}}$ vs.
the threshold $\theta$ }%
\label{figure:Fig-5}%
\end{figure}

\textit{Case two: A lower penalty cost}

From Figure 4, it is seen that the optimal threshold is $\theta^{\ast}=3$ and
$\eta^{\mathbf{d}_{\Delta,\theta^{\ast}}}=22.9$. From Case two of\ Example 1,
we obtained $\eta^{\mathbf{d}^{\ast}}=22.9$. This gives that $\eta
^{\mathbf{d}_{\Delta,\theta^{\ast}}}=\eta^{\mathbf{d}^{\ast}}=22.9$.
Therefore, the optimal static rationing policy is the same as the optimal
dynamic rationing policy, and it is optimal in the policy set $\mathcal{D}$,
and the optimal dynamic rationing policy is of threshold type. Thus,
$\mathbf{d}^{\ast}=\left(  0;1,1,\ldots,1;1,1,\ldots,1\right)  $, given in
Example 1, is the optimal dynamic rationing policy. \begin{figure}[ptbh]
\centering                          \includegraphics[width=8.5cm]{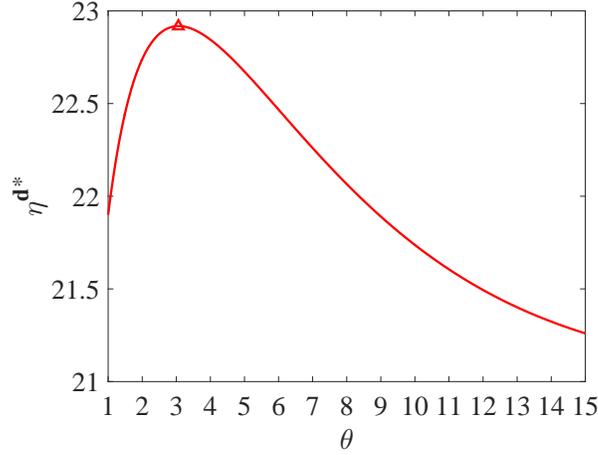}
\caption{The optimal long-run average profit $\eta^{\mathbf{d}^{\ast}}$ vs.
the threshold $\theta$ }%
\label{figure:Fig-6}%
\end{figure}

\textbf{Example 3}. We analyze how the optimal long-run average profit of the
stock-rationing queue depends on the arrival rate. Our observation focuses on
the higher penalty cost $P=10$ and the lower penalty cost $P=0.1$,
respectively. To do this, we further take the system parameters: $\mu_{1}=30$,
$\mu_{2}=40$ and the threshold: $K=5$, $6$, $10$.

\textit{Case one: A higher penalty cost}

Let $P=10$ and $\mathbf{d}^{\ast}=\left(  0;0,0,\ldots,0;1,1,\ldots,1\right)
$. From Figure 5, it is seen that the optimal long-run average profit
$\eta^{\mathbf{d}^{\ast}}$ increases as $\lambda$ increases. In addition, with
the threshold $K$ increases, the optimal long-run average profit
$\eta^{\mathbf{d}^{\ast}}$ increases less slowly as $\lambda$ increases.
\begin{figure}[ptbh]
\centering                          \includegraphics[width=8.5cm]{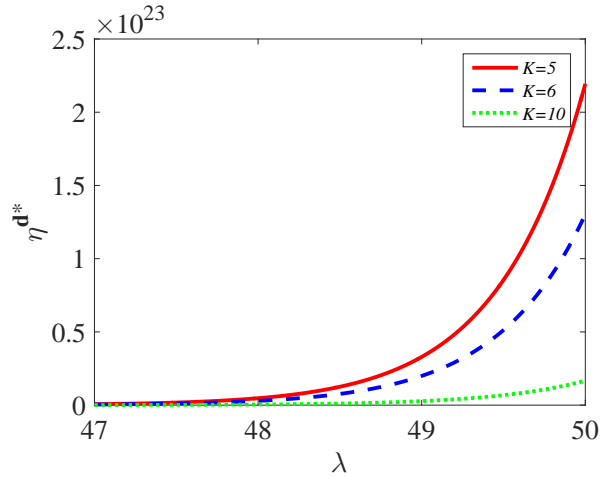}
\caption{$\eta^{\mathbf{d}^{\ast}}$ vs. $\lambda$ under three different
thresholds $K$}%
\label{figure:Fig-3}%
\end{figure}

\textit{Case two: A lower penalty cost}

Let $P=0.1$ and $\mathbf{d}^{\ast}=\left(  0;1,1,\ldots,1;1,1,\ldots,1\right)
$. We discuss how the optimal long-run average profit $\eta^{\mathbf{d}^{\ast
}}$ depends on $\lambda$ for $\lambda\in\left(  65,80\right)  $. From Figure
6, it is seen that the optimal long-run average profit $\eta^{\mathbf{d}%
^{\ast}}$ increases as $\lambda$ increases. In addition, with the threshold
$K$ increases,\ the optimal long-run average profit $\eta^{\mathbf{d}^{\ast}}$
increases less slowly as $\lambda$ increases. \begin{figure}[ptbh]
\centering                          \includegraphics[width=8.5cm]{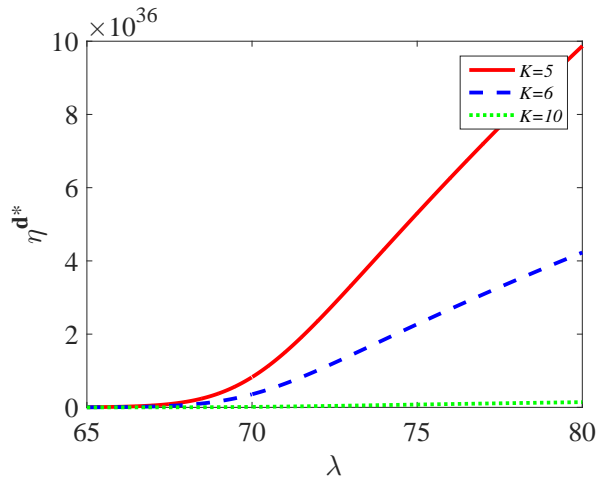}
\caption{$\eta^{\mathbf{d}^{\ast}}$ vs. $\lambda$ under three different
thresholds $K$}%
\label{figure:Fig-4}%
\end{figure}

\textbf{Example 4}. Our observation is to focus on how the penalty cost $P$
influences the long-run average profit $\eta^{\mathbf{d}}$ for any given
policy $\mathbf{d}$. From (\ref{Equa-9}), it is easy to see that for any given
policy $\mathbf{d}$, the long-run average profit $\eta^{\mathbf{d}}$ is linear
in the penalty cost $P$. To show this, we take the system parameters:
$P\in\left(  0,50\right)  $, $\mu_{1}=4$, $\mu_{2}=2$, $\lambda=3$ and $K=15$.
In this case, we observe the special policy: $\mathbf{d}_{1}=$ $\mathbf{d}%
^{\ast}=\left(  0;0,0,\ldots,0;1,1,\ldots,1\right)  $. Figure 7 shows that for
the special policy $\mathbf{d}_{1}=$ $\mathbf{d}^{\ast}$, the long-run average
profit $\eta^{\mathbf{d}^{\ast}}$ linearly decreases as $P$ increases.
\begin{figure}[ptbh]
\centering                          \includegraphics[width=8.5cm]{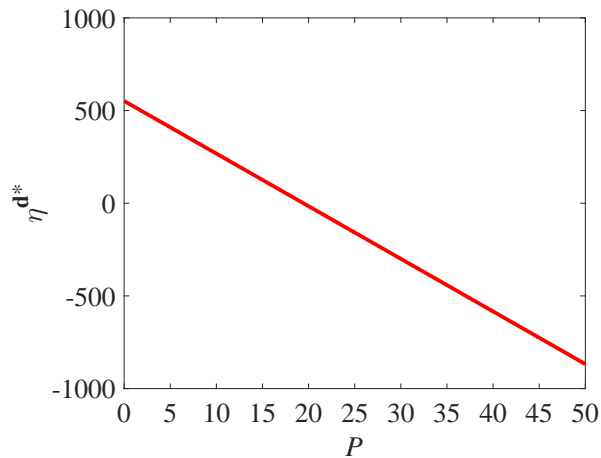}
\caption{The long-run average profit $\eta^{\mathbf{d}^{\ast}}$ vs. the
penalty cost $P$ }%
\label{figure:Fig-7}%
\end{figure}

\section{Concluding Remarks}

In this paper, we highlight intuitive understanding on the optimal dynamic
rationing policy of the stock-rationing queue with two demand classes by means
of the sensitivity-based optimization. To find the optimal dynamic rationing
policy, we establish a policy-based birth-death process and a more general
reward function such that the long-run average profit of the stock-rationing
queue is expressed explicitly. Furthermore, we set up a policy-based Poisson
equation and provide an explicit expression for its solution. Based on this,
we derive a performance difference equation between any two policies such that
we can find the optimal dynamic rationing policy and compute the maximal
long-run average profit from three different areas of the penalty costs.
Therefore, we provide an algebraic method to set up a complete algebraic
solution to the optimal dynamic rationing policy. We show that the optimal
dynamic policy must be of transformational threshold type, which leads to
refining three simple sufficient conditions under each of which the optimal
dynamic policy is of threshold type. In addition, we develop some new
structural properties (e.g., set-structured monotonicity, and class property
of policies) of the optimal dynamic rationing policy. Therefore, we believe
that the methodology and results developed in this paper can be applicable to
analyzing more general stock-rationing queues, and open a series of
potentially promising research.

Along such a line, there are a number of interesting directions for potential
future research, for example:

$\bullet$ Extending to the stock-rationing queues with multiple demand
classes, multiple types of products, backorders, batch order, batch
production, and so on.

$\bullet$ Analyzing non-Poisson input, such as Markovian arrival processes
(MAPs); and/or non-exponential service times, e.g. the PH distributions.

$\bullet$ Discussing how the long-run profit can be influenced by some concave
or convex reward functions.

$\bullet$ Studying individual or social optimization for stock-rationing
queues from a perspective of game theory by means of the sensitivity-based optimization.

\section*{Acknowledgements}

The authors thank the associate editor and two anonymous reviewers for their
many valuable comments to improve the presentation of this paper. Quan-Lin Li
thanks Shaohui Zheng at Hong Kong University of Science and Technology (HKUST)
for providing the inventory control problem during his visiting Professor
Zheng in 2007. Furthermore, the authors gratitude Xi-Ren Cao at HKUST, Li Xia
at Sun Yat-sen University and Xiaole Wu at Fudan University for their valuable
discussion and suggestions. Quan-Lin Li was supported by the National Natural
Science Foundation of China under grants No. 71671158 and 71932002 and by the
Beijing Social Science Foundation Research Base Project under grant No. 19JDGLA004.

\end{document}